# ON THE RATIONALIZATION OF THE $K(n)$-LOCAL SPHERE

TOBIAS BARTHEL, TOMER M. SCHLANK, NATHANIEL STAPLETON, AND JARED WEINSTEIN

*Dedicated to the memory of Jack Morava (1944–2025)*

Abstract. We compute the rational homotopy groups of the $K(n)$-local sphere for all heights $n$ and all primes $p$, verifying a prediction that goes back to the pioneering work of Morava in the early 1970s. More precisely, we show that the inclusion of the Witt vectors into the Lubin–Tate ring induces a split injection on continuous stabilizer cohomology with torsion cokernel of bounded exponent, thereby proving Hopkins' chromatic splitting conjecture and the vanishing conjecture of Beaudry–Goerss–Henn rationally. The key ingredients are the equivalence between the Lubin–Tate tower and the Drinfeld tower due to Faltings and Scholze–Weinstein, integral $p$-adic Hodge theory, and an integral refinement of a theorem of Tate on the Galois cohomology of non-archimedean fields.

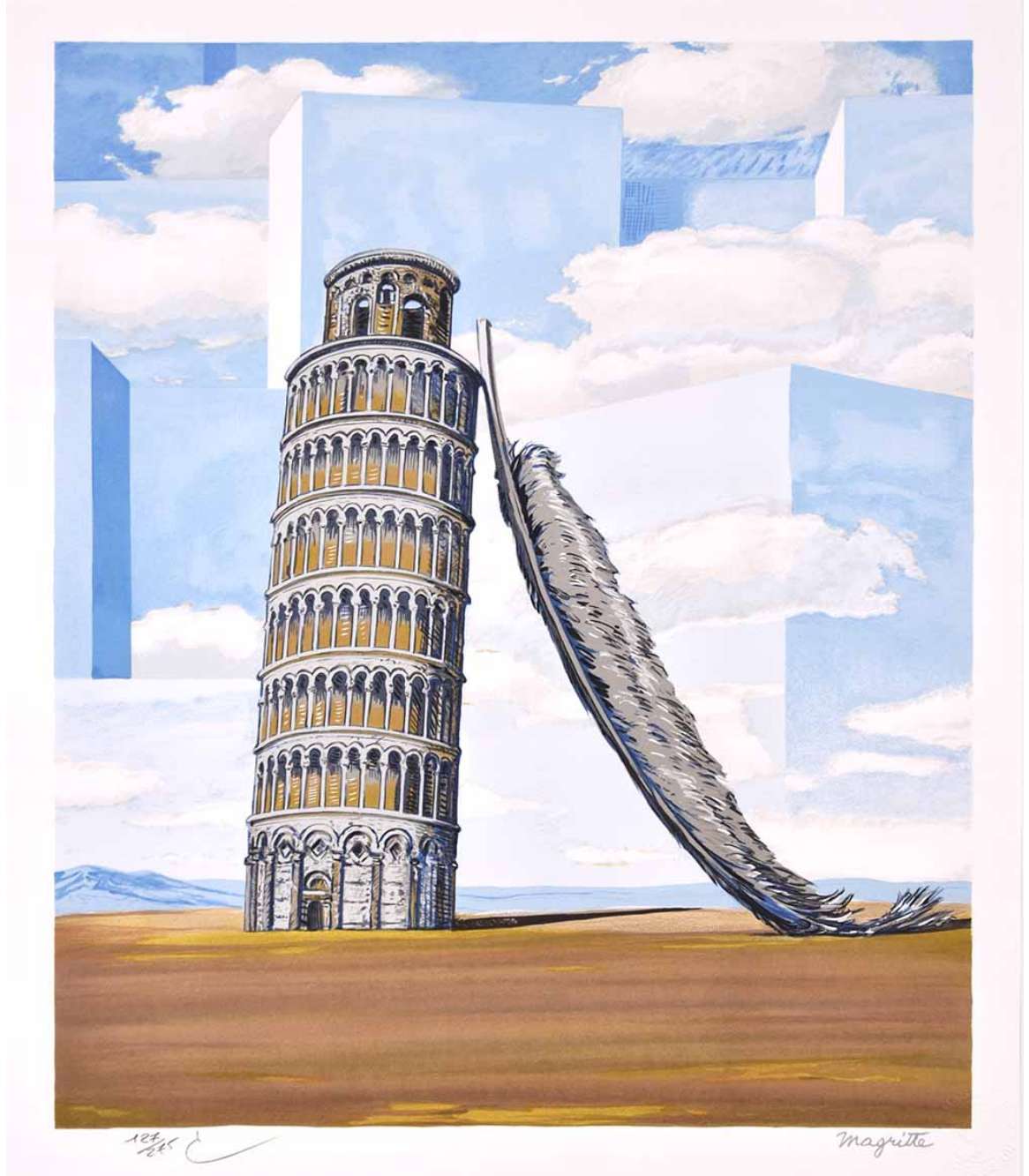

René Magritte, *Memory of a Journey*, 1955.

## 1. Introduction

A central problem in homotopy theory is to understand the homotopy groups of spheres $\pi_{k+d}S^k$, i.e., the group of continuous pointed maps $S^{k+d} \to S^k$ up to homotopy. It follows from the Freudenthal suspension theorem that $\pi_{k+d}S^k$ stabilizes for $k > d+1$, so one may first seek to determine the stable homotopy groups of spheres, $\pi_d S^0 := \lim_{k\to\infty} \pi_{k+d}S^k$. These are abelian groups. We have $\pi_d S^0 = 0$ for $d < 0$ and an isomorphism $\pi_0 S^0 \cong \mathbb{Z}$ encoding the degree of a map. For $d > 0$, Serre showed that $\pi_d S^0$ is finite.

Early attempts at understanding $\pi_* S^0$ through explicit calculations in small degrees only provided limited information about the large-scale structure. Chromatic homotopy theory begins with the deep observation that the elements of $\pi_* S^0$ may be organized into certain periodic families of increasing periodicity, which depend on a prime number $p$. From a modern perspective, these periodic families are the homotopy groups of localizations $L_n S^0$ of $S^0$. These localizations provide successive approximations to the sphere spectrum

$$S^0 \to \ldots \to L_n S^0 \to \ldots \to L_1 S^0 \to L_0 S^0 \simeq S^0_{\mathbb{Q}}$$

whose associated graded pieces are controlled by the $K(n)$-local spheres $L_{K(n)}S^0$, implicitly depending on $p$. The homotopy limit of this tower recovers $S^0_{(p)}$, the $p$-localization of $S^0$. Since $\pi_* S^0_{(p)}$ is the $p$-localization of $\pi_* S^0$, and since these taken together determine $\pi_* S^0$, a fundamental problem in the field is to understand the homotopy groups $\pi_* L_{K(n)}S^0$.

In the case $n = 0$, we have $L_{K(0)}S^0 \cong S^0_{\mathbb{Q}}$, the rational sphere spectrum. As an immediate consequence of the aforementioned theorem of Serre, $\pi_* S^0_{\mathbb{Q}} \cong \mathbb{Q} \otimes \pi_* S^0$ is $\mathbb{Q}$ in degree 0 and 0 otherwise. In the case $n = 1$, $\pi_* L_{K(1)}S^0$ was calculated in the 1970s by Adams–Baird (unpublished) and Ravenel [Rav84]. The case $n = 2$ took many years of work by many people and was only recently resolved (see [SY95, SW02b, SW02a, GHMR05, Beh12, Koh13, BGH22] for an incomplete list); even stating the answer is very involved. Consequently, a full computation of $\pi_* L_{K(n)}S^0$ for $n > 2$ seems to be out of reach.

In light of this, attention over the last few decades has gradually turned towards understanding structural features of $\pi_* L_{K(n)}S^0$. Since the 1970s and motivated by the work of Lazard and Morava, a guiding problem has been to determine the rank of the torsion-free part of $\pi_* L_{K(n)}S^0$, that is, the dimension of the $\mathbb{Q}_p$-vector space $\mathbb{Q} \otimes \pi_* L_{K(n)}S^0$. Through the full force of the computations mentioned above, this dimension is known for $n \leq 2$ and all primes $p$. In this paper, we resolve this question completely for all $n$ and all primes $p$:

**Theorem A.** *There is an isomorphism of graded $\mathbb{Q}$-algebras*

$$\mathbb{Q} \otimes \pi_* L_{K(n)}S^0 \cong \Lambda_{\mathbb{Q}_p}(\zeta_1, \zeta_2, \ldots, \zeta_n),$$

*where the latter is the exterior $\mathbb{Q}_p$-algebra on generators $\zeta_i$ in degree $1 - 2i$.*

In particular, this result confirms the rational part of Hopkins' chromatic splitting conjecture [Hov95] for all primes $p$ and all heights $n$. Previously, this was only known for $n \leq 2$ through the explicit computation in the works listed above. As such, Theorem A constitutes the first general result in the direction of this conjecture since the construction of the class $\zeta \in \pi_{-1} L_{K(n)} S^0$ by Devinatz and Hopkins [DH04] in the early 2000s. We will give a more thorough explanation of the homotopical context for our results in Section 2.

In order to explain our approach, we recall that the homotopy groups of $L_{K(n)} S^0$ can be approached algebraically through Lubin and Tate's deformation theory of formal groups. Let $\Gamma_n$ be a formal group of dimension 1 and height $n$ over $\overline{\mathbb{F}}_p$, and let $\mathbb{G}_n$ be the so-called Morava stabilizer group, defined as the group of automorphisms of $\Gamma_n$ as a functor on $\mathbb{F}_p$-algebras. Then $\mathbb{G}_n$ is an extension of $\mathrm{Gal}(\overline{\mathbb{F}}_p/\mathbb{F}_p) \cong \widehat{\mathbb{Z}}$ by $\mathrm{Aut}_{\overline{\mathbb{F}}_p}(\Gamma_n)$. Let $W = W(\overline{\mathbb{F}}_p)$ be the ring of $p$-typical Witt vectors. By Lubin–Tate theory, there is a complete local ring $A \cong W[\![u_1, \ldots, u_{n-1}]\!]$, corepresenting deformations of $\Gamma_n$, which admits a continuous action by $\mathbb{G}_n$. Further, the invariant differentials $\omega$ of the universal deformation of $\Gamma_n$ form an invertible $A$-module, and the natural actions of $\mathbb{G}_n$ on $A$ and $\omega$ extend to an action of $\mathbb{G}_n$ on the graded ring $A_* = \bigoplus_{t \in 2\mathbb{Z}} \omega^{\otimes t/2}$ (i.e., $A_*$ is evenly concentrated). The Devinatz–Hopkins spectral sequence [DH04] takes the form

$$H^s_{\mathrm{cts}}(\mathbb{G}_n, A_t) \implies \pi_{t-s} L_{K(n)} S^0,$$

where $H^s_{\mathrm{cts}}$ refers to continuous cohomology. Thus understanding the bi-graded ring $H^s_{\mathrm{cts}}(\mathbb{G}_n, A_t)$ is of great importance in chromatic homotopy theory. However, the action of $\mathbb{G}_n$ on $A_*$ is difficult to describe, see [DH95].

Consider instead the problem of computing $H^*_{\mathrm{cts}}(\mathbb{G}_n, W)$, where $\mathbb{G}_n$ acts on $W$ through its quotient $\mathrm{Gal}(\overline{\mathbb{F}}_p/\mathbb{F}_p)$. A classical theorem of Lazard [Laz65] states that the cohomology of a $p$-adic Lie group $G$ with $\mathbb{Q}_p$-coefficients can be computed in terms of Lie algebra cohomology. Applied to the $p$-adic Lie group $\mathrm{Aut}_{\overline{\mathbb{F}}_p}(\Gamma_n)$, Lazard's theorem provides an isomorphism of graded $\mathbb{Q}_p$-algebras, see Lemma 3.8.5:

$$H^*_{\mathrm{cts}}(\mathbb{G}_n, W) \otimes_{\mathbb{Z}_p} \mathbb{Q}_p \cong \Lambda_{\mathbb{Q}_p}(x_1, x_2, \ldots, x_n) \tag{1.0.1}$$

Here, the right hand side is the exterior $\mathbb{Q}_p$-algebra on generators $x_i$ of degree $2i-1$.

Remarkably, and verified by extensive calculations for heights $n \leq 2$ over the last 40 years, work of Morava [Mor85] from the early 1970s suggests that the natural map of $\mathbb{Z}_p$-modules

$$H^*_{\mathrm{cts}}(\mathbb{G}_n, W) \longrightarrow H^*_{\mathrm{cts}}(\mathbb{G}_n, A) \tag{1.0.2}$$

is a rational isomorphism; i.e., it becomes an isomorphism after inverting $p$. The main result of this paper establishes a refinement of this conjecture:

**Theorem B.** *For every integer $s \geq 0$, the natural map $W \hookrightarrow A$ induces a split injection*

$$H^s_{\mathrm{cts}}(\mathbb{G}_n, W) \hookrightarrow H^s_{\mathrm{cts}}(\mathbb{G}_n, A)$$

*whose complement is killed by a power of $p$. In particular,*

$$H^s_{\mathrm{cts}}(\mathbb{G}_n, W) \otimes_{\mathbb{Z}_p} \mathbb{Q}_p \to H^s_{\mathrm{cts}}(\mathbb{G}_n, A) \otimes_{\mathbb{Z}_p} \mathbb{Q}_p$$

*is an isomorphism.*

We have not attempted to make explicit the power of $p$ which kills the complement of $H^s_{\mathrm{cts}}(\mathbb{G}_n, W)$ in $H^s_{\mathrm{cts}}(\mathbb{G}_n, A)$, though this should be possible by our methods. In fact, at heights $n = 1, 2$ this complement is zero, suggesting that this might be the case in general; this is known as the "chromatic vanishing conjecture".

The proof of Theorem B is summarized in Section 3.9. It relies upon recent advances in $p$-adic geometry. Ultimately we draw much of our power from the isomorphism, due to Faltings [Fal02a] and clarified by Scholze and Weinstein [SW13], between the Lubin–Tate and Drinfeld towers. Faltings' isomorphism may be regarded as an equivalence of diamond stacks:

$$[\mathrm{LT}_K^\diamond / \mathbb{G}_n] \simeq [\mathcal{H}^\diamond / \mathrm{GL}_n(\mathbb{Z}_p)]. \tag{1.0.3}$$

(See Theorem 3.9.1 for explanation of the notation.) The cohomology of the stack on the left accesses $H^*_{\mathrm{cts}}(\mathbb{G}_n, A)$. The main idea is to use this isomorphism to replace the opaque action of $\mathbb{G}_n$ on Lubin–Tate space LT with the far more transparent action of the group $\mathrm{GL}_n(\mathbb{Z}_p)$ on Drinfeld's symmetric space $\mathcal{H}$. Theorem A is then readily deduced from Theorem B via the Devinatz–Hopkins spectral sequence.

The groups appearing in (1.0.3) are profinite, and accordingly the stacks appearing there must be construed as living on the pro-étale topology $X_{\text{proét}}$ on rigid-analytic spaces. Thereby the pro-étale cohomology of rigid-analytic spaces, first considered in [Sch13a] and expanded in [BMS18], enters the picture as an indispensable tool. Much of this article is concerned with controlling the pro-étale cohomology $H^*(X_{\text{proét}}, \hat{\mathcal{O}}^+)$ of the completed structure sheaf $\hat{\mathcal{O}}^+$ on $X_{\text{proét}}$ for a rigid-analytic space $X$ over a local field of characteristic $(0, p)$, for example $X = \mathrm{LT}_K$ or $X = \mathcal{H}$.

The results we obtain are new even for the case when $X = \operatorname{Spa} K$ is a single point. In that case, $H^*(X_{\text{proét}}, \hat{\mathcal{O}}^+)$ is the continuous cohomology $H^*_{\mathrm{cts}}(\mathrm{Gal}(\overline{K}/K), \mathcal{O}_C)$, where $C$ is the completion of an algebraic closure $\overline{K}$ of $K$. One of our results (Theorem 4.0.4) is a refinement of a classical theorem of Tate.

**Theorem C.** *Let $K$ be a local field of characteristic $(0, p)$, and let $C$ be the completion of an algebraic closure $\overline{K}$ of $K$. Let $\mathcal{O}_K \subset K$ and $\mathcal{O}_C \subset C$ denote the valuation rings. There is an isomorphism of graded $\mathcal{O}_K$-modules:*

$$H^*_{\mathrm{cts}}(\mathrm{Gal}(\overline{K}/K), \mathcal{O}_C) \cong \mathcal{O}_K[\varepsilon] \oplus T,$$

*where $\varepsilon$ is a degree 1 element with $\varepsilon^2 = 0$ and $T$ is $p^N$-torsion with $N$ a constant which does not depend on $p$ or $K$ (in fact, $N = 5$ suffices).*

Our results on the pro-étale cohomology of rigid-analytic spaces $X$ apply under certain hypotheses regarding the existence of semi-stable models for $X$ after extension of scalars. The cleanest result we have (Corollary 5.6.9) applies when $X$ admits a semistable model over $\mathcal{O}_K$.

**Theorem D.** *Let $K$ be a local field of characteristic $(0, p)$, and let $\mathfrak{X}$ be a quasi-separated semistable formal scheme over $\mathcal{O}_K$, with rigid-analytic generic fiber $X$. There is a natural isomorphism of graded $K$-vector spaces:*

$$H^*(\mathfrak{X}, \mathcal{O}) \otimes_{\mathcal{O}_K} K[\varepsilon] \cong H^*(X_{\text{proét}}, \hat{\mathcal{O}}^+) \otimes_{\mathcal{O}_K} K$$

*Here $H^*(\mathfrak{X}, \mathcal{O})$ is coherent cohomology, and $\varepsilon$ is a degree 1 element with $\varepsilon^2 = 0$.*

**Outline of the document.** Since this paper is written with an audience of both arithmetic geometers and homotopy theorists in mind, we have chosen to include additional background material that might be familiar for one group but not necessarily the other. In that spirit, we begin in Section 2 with a rapid review of the key players in chromatic homotopy theory, working towards stating the chromatic splitting conjecture. We then construct a splitting on the level of cohomology via power operations and deduce Theorem A from Theorem B. Section 3 then collects preliminary material from arithmetic geometry, including the notion of an adic space, a short treatment of continuous cohomology from the condensed perspective, and the pro-étale topology. The section concludes with an outline of the proof of Theorem B, see Section 3.9, which is then fleshed out in the rest of the paper. In Section 4, we establish integral refinements of a theorem of Tate, including Theorem C, by determining bounds on the torsion exponents in the Galois cohomology of $\mathcal{O}_C$ and its Tate twists, where $C$ is the completion of an algebraic closure of a local field $K$ of mixed characteristic. We then globalize this result in Section 5 to the pro-étale cohomology of the generic fiber of a semistable formal scheme over $\mathcal{O}_K$ with coefficients in the sheaf of bounded functions, and prove Theorem D. Finally, in Section 6, we put all the pieces together to prove Theorem B, by first applying our methods separately to the Drinfeld tower and the Lubin–Tate tower, and then deducing our main theorem via the isomorphism of towers.

**Acknowledgements.** We thank Agnès Beaudry, Pierre Colmez, Paul Goerss, Mike Hopkins, Kiran Kedlaya, Niko Naumann, Eric Peterson, Jay Pottharst, Joaquín Rodrigues Jacinto, Juan Esteban Rodríguez Camargo, Peter Scholze, Ehud de Shalit, and Padmavathi Srinivasan for many helpful discussions. We are also grateful to Pierre Colmez, Jack Morava, Niko Naumann, Ningchuan Zhang, and especially Peter Scholze for comments on a preliminary draft of this paper. We also thank the anonymous referees for their valuable feedback on an earlier version.

TB was supported by the European Research Council (ERC) under Horizon Europe (grant No. 101042990) and would like to thank the Max Planck Institute for its hospitality. NS and TS were supported by the US-Israel Binational Science Foundation under grant 2018389. NS was supported by a Sloan research fellowship, NSF Grants DMS-2304781 and DMS-1906236, and a grant from the Simons Foundation (MP-TSM-00002836, NS). JW was supported by NSF Grant DMS-2401472.

## 2. Chromatic homotopy theory

The goal of this section is to place Theorem B in the context of stable homotopy theory and deduce Theorem A from Theorem B. To this end, we begin with a rapid review of some relevant material from chromatic homotopy theory, before turning to the applications. We refer the reader interested in a more thorough introduction to the subject to the following sources: [Lur10, BB20, BGH22]. The readers familiar with chromatic homotopy theory can safely skip ahead to the new results, beginning in Section 2.5.

### 2.1. Chromatic characteristics.

Our starting point is the chromatic perspective on the category of spectra as envisioned by Morava [Mor85] and Ravenel [Rav84] and established by Devinatz, Hopkins, and Smith [DHS88, HS98]. This story has been told by many, and we take a revisionistic approach following [BB20], also freely using the language of higher algebra as developed by Lurie in [Lur09, Lur17].

In order to motivate the homotopical constructions, let us first recast some familiar concepts from algebra in more category-theoretic language; the resulting definitions can then be transported more easily to higher algebra. Let $\mathrm{Mod}(\mathbb{Z})$ be the symmetric monoidal abelian category of abelian groups. A non-trivial unital associative ring $A \in \mathrm{Mod}(\mathbb{Z})$ is said to be a division algebra (or skew field) if any module over $A$ is free. Two division algebras $A$ and $B$ are said to be of the same characteristic if $A \otimes_{\mathbb{Z}} B \neq 0$. It is straightforward to verify that this notion induces an equivalence relation on the collection of all division algebras in $\mathrm{Mod}(\mathbb{Z})$.

It turns out that we can classify all such characteristics: Indeed, the minimal representative of the equivalence classes of characteristics of division algebras in $\mathrm{Mod}(\mathbb{Z})$ are given by the prime fields $\mathbb{F}_p$ for primes $p$ and $\mathbb{Q}$. This fact is essentially a translation of the basic classification of prime fields in classical algebra.

Stable homotopy theory is the study of the category Sp of spectra, which forms a higher analogue of the category of abelian groups; Waldhausen and May coined the term 'brave new algebra.' The role of the integers is then played by the sphere spectrum $S^0$, and the tensor product is replaced by the smash product written as $\otimes$ or $\otimes_{S^0}$ for emphasis. Equipped with this structure, Sp forms a symmetric monoidal stable $\infty$-category; we may therefore speak of rings and their modules in this setting. Interpreted in this context, formally there is an identification $\mathrm{Sp} = \mathrm{Mod}(S^0)$. The next definition is then the natural higher algebraic counterpart to the concept of characteristic as discussed above:

**Definition 2.1.1.** A *division algebra* in Sp is a unital associative ring spectrum $A$ such that every $A$-module $M$ splits into a direct sum of shifts of free rank 1 modules. Two division algebras $A, B \in \mathrm{Sp}$ are of the same *(chromatic) characteristic* if and only if $A \otimes_{S^0} B \neq 0$.

We again obtain an equivalence relation on the collection of all division algebras in Sp, so we are naturally led to ask if we can understand the equivalence classes. This question has been answered completely in the aforementioned seminal work of Devinatz, Hopkins, and Smith. Stating their classification in the form we want requires a short detour. Let $p$ be a prime, let $n \in \mathbb{N} \cup \{\infty\}$, let $\kappa$ be a perfect field of characteristic $p$, and finally let $\Gamma$ be a 1-dimensional height $n$ commutative formal group over $\kappa$. It is an insight of Morava, based on earlier work of Quillen, that this data lifts to Sp: there exists a multiplicative cohomology theory $K(\Gamma, \kappa)^*$ with the following properties:

(1) The value of $K(\Gamma, \kappa)^*$ on a point is given by

$$K(\Gamma, \kappa)^*(\mathrm{pt}) \cong \begin{cases} \mathbb{Q} & \text{if } n = 0 \\ \kappa[v_n^{\pm 1}] & \text{if } 0 < n < \infty \\ \kappa & \text{if } n = \infty, \end{cases}$$

where $v_n$ is a formal variable in degree $2p^n - 2$.

(2) $K(\Gamma, \kappa)^*$ is complex oriented, and the $K(\Gamma, \kappa)^*$-cohomology of complex projective space represents the formal group $\Gamma$:

$$\mathrm{Spf}(K(\Gamma, \kappa)^*(\mathbf{CP}^\infty)) \cong \Gamma.$$

(3) $K(\Gamma, \kappa)^*$ satisfies the Künneth formula for any two spectra $X, Y$:

$$K(\Gamma, \kappa)^*(X \otimes Y) \cong K(\Gamma, \kappa)^*(X) \otimes_{K(\Gamma,\kappa)^*} K(\Gamma, \kappa)^*(Y).$$

By Brown representability, $K(\Gamma, \kappa)^*$ is represented in the category of spectra by a spectrum $K(\Gamma, \kappa)$, which admits the structure[1] of a unital and associative ring spectrum, known as *Morava K-theory* (at height $n$ and over the field $\kappa$). Since the ring of coefficients $K(\Gamma, \kappa)^*(\mathrm{pt})$ is a graded field, $K(\Gamma, \kappa)$ itself must be a division algebra in the sense of Definition 2.1.1. For example, if $\hat{\mathbb{G}}_a$ is the formal additive group, then $K(\hat{\mathbb{G}}_a, \mathbb{Q})$ and $K(\hat{\mathbb{G}}_a, \mathbb{F}_p)$ represent singular cohomology with coefficients in $\mathbb{Q}$ and $\mathbb{F}_p$, respectively. If $\hat{\mathbb{G}}_m$ is the formal multiplicative group, then $K(\hat{\mathbb{G}}_m, \mathbb{F}_p)$ is a summand of mod $p$ complex topological $K$-theory. Generalizing the last example, for any prime $p$ and any height $n \in \mathbb{N}$, there exists a formal group $\Gamma_n$ of height $n$ over $\mathbb{F}_p$. Mildly abusing notation, we set $K(n, p) \coloneqq K(\Gamma_n, \mathbb{F}_p)$, keeping the choice of field implicit. (After base change to an algebraically closed field, $\Gamma_n$ is unique up to isomorphism by a theorem of Lazard [Laz75].) We also set $K(0, p) = K(\hat{\mathbb{G}}_a, \mathbb{Q})$ and $K(\infty, p) = K(\hat{\mathbb{G}}_a, \mathbb{F}_p)$.

Armed with a good collection of division algebras, we can now return to the classification of characteristics in Sp to state:

**Theorem 2.1.2** (Devinatz–Hopkins–Smith)**.** *The collection of Morava K-theories $K(n, p)$ for $p$ ranging through the primes and $n \in \mathbb{N} \cup \{\infty\}$ forms a complete and pairwise distinct set of representatives for the characteristics of division algebras in* Sp*. Moreover, the $K(n, p)$s are minimal in the sense that any division algebra of the same characteristic as $K(n, p)$ has the structure of a module over $K(n, p)$.*

In other, more plain terms: the Morava $K$-theories $K(n, p)$ provide precisely the prime fields of the category of spectra. A couple of remarks are in order.

[1] In fact, the ring structure is not unique, but none of what follows depends on the choice.

- (Non-commutativity) In contrast to the situation in classical algebra, the Morava $K$-theories in intermediate characteristic $0 < n < \infty$ cannot be made commutative. In fact, they do not even afford the structure of an $\mathbb{E}_2$-ring spectrum, see for instance [ACB19]. This is the main reason to work with division algebras in the definition of characteristic.
- (Interrelation) A finitely generated abelian group $M$ with $M \otimes \mathbb{Q} \neq 0$ also satisfies $M \otimes \mathbb{F}_p \neq 0$ for all primes $p$. This statement has a chromatic refinement: If $X$ is a finite spectrum, i.e., a compact object in Sp, then $K(n,p)^*(X) \neq 0$ implies $K(n+1,p)^*(X) \neq 0$. Both the classical algebraic statement and the chromatic statement are false in general for non-compact objects.

To access and isolate the part of Sp that is visible to a fixed Morava $K$-theory $K(n,p)$, we need another important tool from stable homotopy theory, namely the theory of Bousfield localization. These form a suitable generalization of localizations and (derived) completions familiar from commutative algebra.

Fix an arbitrary spectrum $M \in \mathrm{Sp}$. A spectrum $X$ is said to be *$M$-acyclic* if $M_*(X) = 0$; a spectrum $Y$ is then called *$M$-local* if any map $X \to Y$ from an $M$-acyclic spectrum $X$ is null, i.e., factors through a zero object. Intuitively, we wish to quotient Sp by the ideal of all $M$-acyclic spectra to focus on those spectra which are "seen" by $M$. Bousfield [Bou79] rigorously proved that this works, thereby constructing a localization functor $L_M \colon \mathrm{Sp} \to \mathrm{Sp}$ with the following properties:

(1) $L_M X = 0$ if and only if $X$ is $M$-acyclic;
(2) $L_M$ is idempotent and has essential image spanned by the $M$-local spectra;
(3) for any $Z \in \mathrm{Sp}$, there is a natural map $Z \to L_M Z$ which exhibits $L_M Z$ as the initial map out of $Z$ to an $M$-local spectrum.

These properties characterize $L_M$ uniquely up to homotopy. Finally, we denote the full subcategory of $M$-local spectra by $\mathrm{Sp}_M$; alternatively, this category is obtained from Sp by inverting all $M$-equivalences, i.e., those maps which induce isomorphisms in $M_*$-homology. Via localization, $\mathrm{Sp}_M$ inherits a structure of symmetric monoidal $\infty$-category from Sp, with tensor product given by the *localized* smash product.

Some examples might be illuminating. If $M = \mathbb{Q}$, then $L_{\mathbb{Q}}$ is rationalization, whose effect on homotopy groups of any spectrum is tensoring with $\mathbb{Q}$. The element $p \in \mathbb{Z} \cong \pi_0 S^0$ is represented by a map $p \colon S^0 \to S^0$, whose cofiber we denote by $S^0/p$, the *mod $p$ Moore spectrum*. On the one hand, the local category $\mathrm{Sp}_{S^0/p}$ is the category of $p$-complete spectra, and localization at $S^0/p$ has the effect of derived $p$-completion on homotopy groups. It is customary to write $X_p \coloneqq L_{S^0/p}X$ and $\mathrm{Sp}_p \coloneqq \mathrm{Sp}_{S^0/p}$. On the other hand, localizing at $M = S^0[1/p]$, the colimit over multiplication by $p$ on $S^0$, has the effect of inverting $p$ on $S^0$. Similarly, we can construct spectral analogues of $p$-localization, by inverting all primes but $p$, to obtain the category $\mathrm{Sp}_{(p)}$ of $p$-local spectra. For more information about various localizations on Sp, we refer to [Bou79].

The main example of interest to us is $\mathrm{Sp}_{K(n,p)}$, the so-called $K(n,p)$-local category. In light of Theorem 2.1.2, it forms an irreducible piece of the category of spectra, as we will explain momentarily, and it is one of the key objects of study in chromatic homotopy theory. Note that the functor $L_{K(n,p)}$ and thus $\mathrm{Sp}_{K(n,p)}$ only depend on the characteristic and are in particular independent of the choice of $\Gamma_n$.

2.2. **Chromatic divide and conquer.** From now on, we will restrict attention to the category of $p$-local spectra for a fixed prime $p$, and usually drop the prime from the notation. Intuitively speaking, the idea of the chromatic approach to stable homotopy theory is to filter the category of $p$-local spectra $\mathrm{Sp}_{(p)}$ by its subcategories of mixed chromatic characteristics $(0,1,\ldots,n)$ for $n \to \infty$. Here we say that a spectrum has mixed chromatic characteristic $(0,1,\ldots,n)$ if it is local with respect to the direct sum $K(0) \oplus K(1) \oplus \ldots \oplus K(n)$, and we write $L_n$ for the associated

Bousfield localization. There is then a sequence of Bousfield localization functors and natural transformations,

$$\ldots \to L_n \to L_{n-1} \to \ldots \to L_1 \to L_0 \tag{2.2.1}$$

the so-called *chromatic tower*. Note that the bottom layer is rationalization $L_0 = L_{\mathbb{Q}}$. Also note that the infinite height has been omitted, a curiosity justified by the chromatic convergence theorem of Hopkins–Ravenel [Rav92]: If $X$ is a $p$-local finite spectrum, then it can be recovered from its chromatic tower (2.2.1):

$$X \simeq \lim_n L_n X.$$

Given $X$, it is then sensible to ask for the graded pieces of its chromatic filtration, i.e., the difference between $L_nX$ and $L_{n-1}X$. This is captured by the *chromatic fracture square*, taking the form of a (homotopy) pullback square that exists for any height $n > 0$ and an arbitrary spectrum $X$ (see e.g., [ACB22]):

$$\begin{array}{ccc} L_nX & \longrightarrow & L_{K(n)}X \\ \downarrow & & \downarrow \\ L_{n-1}X & \longrightarrow & L_{n-1}L_{K(n)}X. \end{array} \tag{2.2.2}$$

Geometrically, one should think of this square as being analogous to the gluing square for a sheaf over an open-closed decomposition of a space. In this picture, $L_{K(n)}X$ corresponds to the sheaf over a formal neighborhood of a point, while $L_{n-1}X$ is the restriction to the open complement of the point. The term $L_{n-1}L_{K(n)}X$ along with the maps pointing to it then control the gluing process.

The weakest form of Hopkins' *chromatic splitting conjecture*, recorded in [Hov95], stipulates that the bottom horizontal map in (2.2.2) is split for $X = S^0$ (and hence for any finite spectrum $X$), so that the chromatic assembly process takes a particularly simple form. The strong form of the conjecture (Conjecture 2.4.2 below) gives a complete description of $L_{n-1}L_{K(n)}S^0$ in terms of the $L_iS^0$ for $i \leq n$. Assuming it, one could inductively reduce the study of the $L_nS^0$ to that of the $L_{K(n)}S^0$, whose homotopy groups are identified after inverting $p$ by our Theorem A.

2.3. **Morava $E$-theory and the cohomology of the stabilizer group.** Just as a height $n$ formal group over $\mathbb{F}_p$ gives rise to the spectrum $K(n)$, the Lubin–Tate ring also admits a spectral incarnation. Let $\Gamma$ be a height $n$ formal group over a perfect field $\kappa$ of characteristic $p$. Let $A(\Gamma,\kappa)$ denote its ring of deformations, so that $A(\Gamma,\kappa) \cong W(\kappa)[\![u_1,\ldots,u_{n-1}]\!]$. An unpublished theorem of Goerss–Hopkins–Miller [GH04, HM14], revisited and extended by Lurie in [Lur18], lifts this data to a commutative algebra in $\mathrm{Sp}_{K(n)}$, called $E(\Gamma,\kappa)$, with the property that

$$\pi_*E(\Gamma,\kappa) \cong A(\Gamma,\kappa)[\beta,\beta^{-1}],\ |\beta| = 2.$$

The commutative ring spectrum $E(\Gamma,\kappa)$ is known as *Morava $E$-theory* or *Lubin–Tate theory.*

In fact, Goerss, Hopkins, and Miller prove something stronger. Consider the 1-category of formal groups over perfect fields FG. The objects of FG are given by pairs $(\Gamma,\kappa)$ as above, and a morphism $(\Gamma,\kappa) \to (\Gamma',\kappa')$ in FG consists of a ring map $i\colon \kappa \to \kappa'$ together with an isomorphism of formal groups $i^*\Gamma \xrightarrow{\sim} \Gamma'$. Goerss, Hopkins, and Miller produce a fully faithful functor $E(-,-)\colon \mathrm{FG} \to \mathrm{CAlg}(\mathrm{Sp})$. It is important to note that the source is a 1-category and that the target is an $\infty$-category. Thus this theorem identifies a very rigid portion of CAlg(Sp), in which the mapping spaces are homotopy equivalent to sets.

The underlying spectrum of $E(\Gamma,\kappa)$ is easily constructed as a consequence of the Landweber exact functor theorem. That theorem produces cohomology theories out of the complex cobordism spectrum with specified formal groups, so long as these satisfy a certain hypothesis. The universal deformation $\overline{\Gamma}$ of $\Gamma$ over $A(\Gamma,\kappa)$ satisfies the hypothesis, and this gives rise to Morava $E$-theory. However, producing Morava $E$-theory as a commutative algebra in $\mathrm{Sp}_{K(n)}$ is quite a bit more difficult and requires either obstruction theory or a derived deformation theory of formal groups—this is the content of the Goerss–Hopkins–Miller theorem.

Let $\Gamma_n$ be any formal group of height $n$ over $\overline{\mathbb{F}}_p$; then $\Gamma_n$ is unique up to isomorphism. We write $E_n = E(\Gamma, \overline{\mathbb{F}}_p)$. We let $\mathbb{G}_n = \mathrm{Aut}(\Gamma_n, \overline{\mathbb{F}}_p) \cong \mathrm{Aut}\, E_n$, the *Morava stabilizer group*. As a topological group, $\mathbb{G}_n$ may be identified with the profinite completion $\widehat{D^\times}$ of $D^\times$, where $D$ is the central simple algebra of invariant $\frac{1}{n}$ over $\mathbb{Q}_p$, with ring of integers $\mathcal{O}_D = \mathrm{End}_{\overline{\mathbb{F}}_p}(\Gamma_n)$.

It turns out that $L_{E_n} = L_n$ by [Rav84, Theorem 2.1(d)], so $E_n$ does not provide us with a new localization functor. Furthermore, there is a close relationship between $E_n$ and $K(n)$, akin to the one between a local ring and its residue field. Since the ideal $I_n = (p, u_1, \ldots, u_{n-1})$ of $\pi_0 E_n$ is generated by a regular sequence, we may form a (not necessarily commutative) ring spectrum $E_n/I_n$ by iterated cofibers. This will have the property that $\pi_*(E_n/I_n) \cong \pi_*(E_n)/I_n = \kappa[\beta, \beta^{-1}]$. The ring spectrum obtained in this way is equivalent as a spectrum to a finite direct sum of suspensions of $K(\Gamma_n, \kappa)$:

$$E_n/I_n \simeq \bigoplus_{0 \le i \le p^n - 2} \Sigma^{2i} K(\Gamma_n, \kappa).$$

Here, the direct sum accounts for the fact that $E_n$ is 2-periodic, while the periodicity of $K(\Gamma_n, \kappa)$ is $2p^n - 2$. We refer the reader to [BB20] and the references given therein for more details on this material.

The unit in $\mathrm{Sp}_{K(n)}$ is the $K(n)$-local sphere, denoted $L_{K(n)}S^0$. Being the unit, $L_{K(n)}S^0$ is the initial object in $\mathrm{CAlg}(\mathrm{Sp}_{K(n)})$, so there is a canonical map of commutative algebras $L_{K(n)}S^0 \to E_n$. Since $\mathbb{G}_n$ acts on $E_n$ through commutative algebra automorphisms, this map is $\mathbb{G}_n$-equivariant for the trivial action on $L_{K(n)}S^0$. A result of Devinatz and Hopkins [DH04], reinterpreted in Rognes' framework of spectral Galois extensions [Rog08], says that the unit map $L_{K(n)}S^0 \to E_n$ exhibits $E_n$ as a pro-Galois extension of $L_{K(n)}S^0$ in $\mathrm{Sp}_{K(n)}$, with Galois group $\mathbb{G}_n$. Concretely, this means that we have canonical equivalences of commutative ring spectra

$$L_{K(n)}S^0 \simeq E_n^{h\mathbb{G}_n} \quad \text{and} \quad L_{K(n)}(E_n \otimes E_n) \simeq C_{\mathrm{cts}}(\mathbb{G}_n, E_n), \tag{2.3.1}$$

where $C_{\mathrm{cts}}(\mathbb{G}_n, E_n)$ denotes the ring spectrum of continuous functions on $\mathbb{G}_n$ with coefficients in $E_n$. This enables us to run Galois descent along $L_{K(n)}S^0 \to E_n$. Form the associated $K(n)$-local cosimplicial Amitsur complex

$$L_{K(n)}S^0 \to E_n^{\hat{\otimes}\bullet+1} \coloneqq L_{K(n)}\left(E_n \rightrightarrows E_n \otimes E_n \,\substack{\rightarrow\\[-1ex]\rightarrow\\[-1ex]\rightarrow}\, E_n^{\otimes 3} \ldots\right), \tag{2.3.2}$$

where we have omitted the degeneracy maps from the display. Applying the homotopy groups to the resolution (2.3.2) and using (2.3.1) to identify the abutment and $E_2$-page, we obtain a Bousfield–Kan spectral sequence of signature

$$E_2^{s,t} \cong H^s_{\mathrm{cts}}(\mathbb{G}_n, \pi_t E_n) \implies \pi_{t-s} L_{K(n)}S^0, \tag{2.3.3}$$

where $H^*_{\mathrm{cts}}$ denotes cohomology with continuous cocycles.

Since $\mathbb{G}_n$ has finite virtual cohomological dimension ([Mor85, Section 2.2.0]), this spectral sequence provides an excellent approximation to the homotopy groups of $L_{K(n)}S^0$: it converges strongly with a finite horizontal vanishing line on some finite page, i.e., there exists $r \geq 2$ and some $N > 0$ such that $E_r^{s,t} = 0$ for all $s > N$. In fact, if $p$ is odd and $2(p-1) \geq n^2$, then we may take $r = 2$ and $N = n^2$. The spectral sequence (2.3.3), often simply referred to as "the" descent spectral sequence in chromatic homotopy theory, provides a gateway between stable homotopy theory and $p$-adic geometry.

2.4. **The chromatic splitting conjecture and the vanishing conjecture.** Computational evidence at low heights $n \leq 2$ suggests that the continuous cohomology of the action of $\mathbb{G}_n$ on $W \subseteq \pi_0 E_n$ largely controls the behavior of $\pi_* L_{K(n)} S^0$, as we shall now explain. Recall the isomorphism (1.0.1)

$$H^*_{\mathrm{cts}}(\mathbb{G}_n, W) \otimes_{\mathbb{Z}_p} \mathbb{Q}_p \cong \Lambda_{\mathbb{Q}_p}(x_1, x_2, \ldots, x_n),$$

where the latter is the exterior $\mathbb{Q}_p$-algebra on generators $x_i$ in degree $2i-1$. The relevance of these classes for stable homotopy theory was first realized by Morava, who also exhibited explicit cocycle representatives for the $x_i$, see [Mor85, Remark 2.2.5].

Each of the classes $x_i$ can be lifted[2] to a class $\tilde{x}_i$ in the integral cohomology ring $H^*_{\mathrm{cts}}(\mathbb{G}_n, W)$. Let

$$\varphi\colon H^*_{\mathrm{cts}}(\mathbb{G}_n, W) \to H^*_{\mathrm{cts}}(\mathbb{G}_n, \pi_0 E_n) \cong E_2^{*,0} \tag{2.4.1}$$

be the natural map induced from the inclusion $W \hookrightarrow A \cong \pi_0 E_n$.

Hopkins' *chromatic splitting conjecture*, recorded in [Hov95], predicts that the classes $\tilde{x}_i$ carry all the relevant information about the homotopy groups of $L_{K(n)} S^0$; more precisely:

**Conjecture 2.4.2** (Chromatic splitting conjecture)**.** *Let $n \geq 0$. For $p$ sufficiently large with respect to $n$ and each $i = 1, \ldots, n$, we have the following behavior.*

*(1) For each $i = 1, \ldots, n$, the class $\varphi(\tilde{x}_i)$ survives the spectral sequence* (2.3.3)*, thereby giving rise to a homotopy class*

$$e_i \in \pi_{1-2i} L_{K(n)} S^0,$$

*or all the same, a map $e_i\colon S_p^{1-2i} \to L_{K(n)} S^0$.*

*(2) Let $S_p^d$ be the $p$-complete sphere spectrum shifted in degree $d \in \mathbb{Z}$. The composition $S_p^{1-2i} \xrightarrow{e_i} L_{K(n)} S^0 \to L_{n-1} L_{K(n)} S^0$ factors through a map $\overline{e}_i\colon L_{n-i} S_p^{1-2i} \to L_{n-1} L_{K(n)} S^0$.*

*(3) The $\overline{e}_i$ induce an equivalence of spectra*

$$L_{n-1} L_{K(n)} S^0 \cong \bigoplus_{\substack{0 \leq j \leq n \\ 1 \leq i_1 < \ldots < i_j \leq n}} \left( \bigotimes_{k=1}^{j} L_{n-i_k} S_p^{1-2i_k} \right),$$

*where the tensor product is taken in the category of $L_{n-1} S_p^0$-modules. Here, the right hand side is indexed on the $\mathbb{Z}_p$-module generators of the exterior algebra $\Lambda_{\mathbb{Z}_p}(\overline{e}_1, \overline{e}_2, \ldots, \overline{e}_n)$. This becomes an equivalence of ring spectra when the right side is endowed with the ring structure in which the classes $\overline{e}_i$ constructed in* (2) *form exterior generators.*

A more refined formulation of the conjecture as well as the necessary modifications for the prime 2 can be found in [BGH22]. It has been verified by explicit computation of both sides for heights $n \leq 2$ and all primes $p$. After tensoring with $\mathbb{Q}$, the chromatic splitting conjecture predicts that

$$\mathbb{Q} \otimes \pi_* L_{K(n)} S^0 \cong \Lambda_{\mathbb{Q}_p}(\overline{e}_1, \overline{e}_2, \ldots, \overline{e}_n), \ |e_i| = 1 - 2i. \tag{2.4.3}$$

Indeed, there are natural equivalences $\mathbb{Q} \otimes L_i X \simeq \mathbb{Q} \otimes X$ for all spectra $X$ and all $i \geq 0$, so (3) of Conjecture 2.4.2 rationalizes to (2.4.3). Further chromatic consequences of Conjecture 2.4.2 can be found in [BB20] and [Mor14].

A related question concerns the map $\varphi$ appearing in (2.4.1). The following conjecture[3] was formulated by Beaudry, Goerss, and Henn in [BGH22, Page 3].

**Conjecture 2.4.4** (Vanishing conjecture)**.** *The inclusion of coefficients $W \hookrightarrow \pi_0 E_n$ induces an isomorphism $H^*_{\mathrm{cts}}(\mathbb{G}_n, W) \cong H^*_{\mathrm{cts}}(\mathbb{G}_n, \pi_0 E_n)$.*

[2] In the chromatic practice, there are certain preferred choices of lifts, but these will not matter for our purposes here. For details, we refer to [Hov95].

[3] The name comes from the equivalent formulation that $H^*_{\mathrm{cts}}(\mathbb{G}_n, \pi_0 E_n / W) = 0$.

Conjecture 2.4.4 has been verified for all $n \leq 2$ and for all primes $p$; moreover, in these cases, it does not need to be modified for $p = 2$, but rather accounts for the additional complications witnessed there. If correct, this conjecture would substantially simplify the task of understanding $\pi_* L_{K(n)} S^0$. Our goal in this paper is to prove (2.4.3) for all heights $n$ and all primes $p$.

2.5. **Power operations and the splitting of $W \to A$.** Fix a prime $p$ and a height $n \geq 1$, and let $E = E_n$ be Morava $E$-theory. Thus $A = E^0 \cong W[\![u_1, \ldots, u_{n-1}]\!]$, where $W = W(\overline{\mathbb{F}}_p)$. For the first time in the article, we introduce new material.

Our proof of Theorem A begins with the following proposition.

**Proposition 2.5.1.** *The inclusion $W \hookrightarrow A$ admits a continuous $\mathbb{G}_n$-equivariant (additive) splitting. In other words, there is a $\mathbb{G}_n$-equivariant decomposition of topological abelian groups*

$$A \cong W \oplus A^c.$$

Our proof of Proposition 2.5.1 uses input from homotopy theory, namely the transfer maps and power operations on Morava $E$-theory. Many general results regarding power operations can be found in [BMMS86]. We briefly recall this theory here.

Let $F$ be an $\mathbb{E}_\infty$-ring. For each $m \geq 1$, the $\mathbb{E}_\infty$-ring structure on $F$ induces an $m$-fold multiplication map

$$(F^{\otimes m})_{h\Sigma_m} \to F, \tag{2.5.2}$$

where $(-)_{h\Sigma_m}$ denotes homotopy orbits with respect to the permutation action of the symmetric group $\Sigma_m$ on $F^{\otimes m}$.

For $m \geq 1$, the $m$th power operation on $F^0 = \pi_0 F$ is the multiplicative map of rings

$$P^m \colon F^0 \to F^0(B\Sigma_m)$$

defined as the composite

$$[S^0, F] \to [S^0_{h\Sigma_m}, F^{\otimes m}_{h\Sigma_m}] \to [S^0_{h\Sigma_m}, F] = F^0(B\Sigma_m).$$

Here, the first map inputs a map $S^0 \to F$, applies the $m$th tensor power to obtain a map $S^0 \cong (S^0)^{\otimes m} \to F^{\otimes m}$ (recall that $S^0$ is the unit in spectra), and then passes to the homotopy orbits for the $\Sigma_m$ action on either side. The second map is given by postcomposition with the multiplication map (2.5.2). By convention, $P^0$ is the constant function with value 1, and $P^1$ is identity map on $F^0$.

In the situation that $F$ is also $K(n)$-local, the power operation $P^m \colon F^0 \to F^0(B\Sigma_m)$ can be upgraded to a continuous map between topological abelian groups. In general, if $X$ is a $K(n)$-local spectrum, then the abelian group $\pi_0 X$ admits a topology, as studied in [HS99, Chapter 11]. The key point is that the unit object $L_{K(n)} S^0$ in the $K(n)$-local category can be expressed as a sequential colimit

$$L_{K(n)} S^0 \simeq \operatorname{colim}_{k \in \mathbb{N}} M_k, \tag{2.5.3}$$

where each $M_k$ is a compact object in the $K(n)$-local category. For instance, one could take $M_k$ to be the Spanier–Whitehead dual of the $K(n)$-localized generalized Moore spectrum $S^0/(p^{\ell_k}, v_1^{\ell_k}, \ldots, v_{n-1}^{\ell_k})$, for appropriately chosen $\ell_k$.

Now if $X$ is a $K(n)$-local spectrum, we can construct a map

$$\pi_0 X = [S^0, X] \cong [L_{K(n)} S^0, X] \to \lim_{k \in \mathbb{N}} [M_k, X].$$

We give $\pi_0 X$ the coarsest topology which makes the map continuous. Any map of $K(n)$-local spectra $X \to Y$ induces a continuous map $\pi_0 X \to \pi_0 Y$ [HS99, Proposition 11.1.(b)].

Returning to the setup of power operations, suppose $F$ is a $K(n)$-local $\mathbb{E}_\infty$-ring. Then $F^{B\Sigma_m}$ is also $K(n)$-local, and so both $\pi_0 F$ and $\pi_0(F^{B\Sigma_m}) = F^0(B\Sigma_m)$ are topological abelian groups.

**Lemma 2.5.4.** *Let $F$ be a $K(n)$-local $\mathbb{E}_\infty$-ring, and let $m \geq 0$. The power operation $P^m \colon F^0 \to F^0(B\Sigma_m)$ is continuous.*

*Proof.* The case $m = 0$ is trivial, so assume $m \geq 1$. The continuity of $P^m$ can be expressed this way: for every $k \in \mathbb{N}$ there exists $j \in \mathbb{N}$ and a map $f_k \colon [M_j, F] \to [M_k, F^{B\Sigma_m}]$ such that

$$\begin{array}{ccc}
[L_{K(n)}S^0, F] & \xrightarrow{P^m} & [L_{K(n)}S^0, F^{B\Sigma_m}] \\
\downarrow & & \downarrow \\
[M_j, F] & \xrightarrow{f_k} & [M_k, F^{B\Sigma_m}]
\end{array}$$

commutes.

We begin by giving a presentation of $(L_{K(n)}S^0)^{\otimes m}_{h\Sigma_m}$ as a colimit:

$$\begin{aligned}
(L_{K(n)}S^0)^{\otimes m}_{h\Sigma_m} &\simeq ((\operatorname{colim}_{k\in\mathbb{N}} M_k)^{\otimes m})_{h\Sigma_m} \\
&\simeq (\operatorname{colim}_{(k_1,\dots,k_m)\in\mathbb{N}^m}(M_{k_1} \otimes \dots \otimes M_{k_m}))_{h\Sigma_m} \\
&\simeq (\operatorname{colim}_{j\in\mathbb{N}} M_j^{\otimes m})_{h\Sigma_m} \\
&\simeq \operatorname{colim}_{j\in\mathbb{N}}((M_j^{\otimes m})_{h\Sigma_m}).
\end{aligned}$$

The first equivalence is the $m$th tensor power of (2.5.3). The second equivalence follows from the fact that the tensor product commutes with colimits in each variable. The third equivalence follows from the fact that the diagram $\mathbb{N}$ is filtered, which implies that the diagonal map $\mathbb{N} \to \mathbb{N}^m$ is cofinal. The fourth equivalence follows from the fact that colimits commute.

The key point now is that each $(M_k)_{h\Sigma_m}$ is compact. Indeed, it is equivalent to $M_k \otimes (L_{K(n)}S^0)_{h\Sigma_m}$, where the tensor product is taken in $\mathrm{Sp}_{K(n)}$. We apply the facts that $M_k$ is compact by construction, $(L_{K(n)}S^0)_{h\Sigma_m} \simeq L_{K(n)}B\Sigma_m$ is dualizable (see [HS99, Cor. 8.7]), and the tensor product of a compact object with a dualizable object is compact.

Therefore the map $(M_k)_{h\Sigma_m} \to (L_{K(n)}S^0)_{h\Sigma_m} \xrightarrow{\sim} (L_{K(n)}S^0)^{\otimes m}_{h\Sigma_m}$ must factor through $(M_j^{\otimes m})_{h\Sigma_m}$ for some $j$, which is to say that there is a map $(M_k)_{h\Sigma_m} \to (M_j)^{\otimes m}_{h\Sigma_m}$ making the diagram

$$\begin{array}{ccc}
(M_k)_{h\Sigma_m} & \longrightarrow & (M_j)^{\otimes m}_{h\Sigma_m} \\
\downarrow & & \downarrow \\
(L_{K(n)}S^0)_{h\Sigma_m} & \xrightarrow{\simeq} & (L_{K(n)}S^0)^{\otimes m}_{h\Sigma_m}
\end{array}$$

commute.

We now have a commutative diagram

$$\begin{array}{ccccccc}
[L_{K(n)}S^0, F] & \longrightarrow & [(L_{K(n)}S^0)^{\otimes m}_{h\Sigma_m}, F^{\otimes m}_{h\Sigma_m}] & \xrightarrow{\cong} & [(L_{K(n)}S^0)_{h\Sigma_m}, F^{\otimes m}_{h\Sigma_m}] & \longrightarrow & [(L_{K(n)}S^0)_{h\Sigma_m}, F] \\
\downarrow & & \downarrow & & \downarrow & & \downarrow \\
[M_j, F] & \longrightarrow & [(M_j)^{\otimes m}_{h\Sigma_m}, F^{\otimes m}_{h\Sigma_m}] & \longrightarrow & [(M_k)_{h\Sigma_m}, F^{\otimes m}_{h\Sigma_m}] & \longrightarrow & [(M_k)_{h\Sigma_m}, F]
\end{array}$$

in which the composition of the upper horizontal arrows in $P^m$. The desired map $f_k$ can now be defined as the composition of the lower horizontal arrows.

∎

We return our focus to Morava $E$-theory $E = E_n$. Goerss, Hopkins, and Miller [GH04, Chapter 7] show that $E$ admits an essentially unique $\mathbb{E}_\infty$-ring structure. Therefore there are power operations $P_m \colon E^0 \to E^0(B\Sigma_m)$.

**Corollary 2.5.5.** *The power operations on $E^0$ are continuous with respect to the $I$-adic topology, where $I \subset E^0$ is the maximal ideal.*

*Proof.* By [HS99, Proposition 11.9], if $X$ is a dualizable object of $\mathrm{Sp}_{K(n)}$ and $E$ is height $n$ Morava $E$-theory, then the topology on $\pi_0(E^X) = E^0(X)$ is the $I$-adic topology, where $I \subset \pi_0 E$ is the maximal ideal. Applied to the case that $X$ is the suspension spectrum of $BG$ for $G$ trivial and $G = \Sigma_m$, respectively, we find that the topologies on $E^0$ and $E^0(B\Sigma_m)$ are the $I$-adic ones. The corollary is now a case of Lemma 2.5.4. □

*Remark* 2.5.6. Alternatively, one could deduce Lemma 2.5.4 by endowing $E_n$ with the structure of a condensed commutative ring spectrum, as for example in [Mor23, Section 2.3].

Let $G$ be a finite group, and let $H \subseteq G$ be a subgroup. The induced map between classifying spaces $BH \to BG$ is a finite cover, so for any cohomology theory $F$, there is a homomorphism of graded abelian groups $\mathrm{Tr}_H^G \colon F^*(BH) \to F^*(BG)$, known as the transfer map. If $F$ is $K(n)$-local for some $n \geq 1$, the transfer map can be made more general: any homomorphism of finite groups $H \to G$ induces a map of dualizable $K(n)$-local spectra

$$L_{K(n)}\Sigma^\infty_+ BG \to L_{K(n)}\Sigma^\infty_+ BH,$$

which in turn induces a transfer map $\mathrm{Tr}_H^G \colon F^*(BH) \to F^*(BG)$ [GS96, Theorem 1.1 and Corollary 1.2].

The transfer map is transitive with respect to composable morphisms between finite groups. In degree 0, the transfer map $F^0(BH) \to F^0(BG)$ is continuous, since it arises by taking $\pi_0$ of a map of $K(n)$-local spectra. Specializing further to the case $F = E$, the transfer map $\mathrm{Tr}_H^G \colon E^0(BH) \to E^0(BG)$ is continuous with respect to the $I$-adic topology on either side.

We are at last ready to prove Proposition 2.5.1.

*Proof of Proposition 2.5.1.* Recall that $A = E^0$. Let $P^m \colon E^0 \to E^0(B\Sigma_m)$ be the $m$th power operation, determined by the $\mathbb{E}_\infty$-ring structure on $E$. We have the following identity [BMMS86, Chapter VIII, Proposition 1.1.(vi)]:

$$P^m(x+y) = \sum_{i+j=m} \mathrm{Tr}_{\Sigma_i \times \Sigma_j}^{\Sigma_m}(P^i(x)P^j(y)) \tag{2.5.7}$$

for $x, y \in E^0$.

Let

$$\beta_m \colon E^0 \xrightarrow{P^m} E^0(B\Sigma_m) \xrightarrow{\mathrm{Tr}_{\Sigma_m}^e} E^0$$

be the composite of the $m$th power operation with the $K(n)$-local transfer map along the surjection from $\Sigma_m$ to the trivial group. Since $P^m$ is continuous (Corollary 2.5.5) and so is $\mathrm{Tr}_{\Sigma_m}^e$, the map $\beta_m$ is continuous.

It follows from the formula for $P^m(x+y)$ in (2.5.7) together with the transitivity of the transfer that $\beta_m(x+y) = \sum_{i+j=m} \beta_i(x)\beta_j(y)$. This implies that the formal sum $\beta(x) = \sum_{m \geq 0} \beta_m x^m$, considered as map $E^0 \to E^0[\![x]\!]$, satisfies $\beta(x+y) = \beta(x)\beta(y)$. Since $\beta_0(a) = 1$ and $\beta_1(a) = a$ for all $a \in E^0$, the map $\beta$ factors through a homomorphism from the additive group $E^0$ to the subgroup $1 + xE^0[\![x]\!]$ of $E^0[\![x]\!]^\times$. Now we may quotient the target by the maximal ideal in $E^0$ to obtain a map $E^0 \to \overline{\mathbb{F}}_p[\![x]\!]$ that sends addition to multiplication. The big Witt vectors $W_{\mathrm{big}}(\overline{\mathbb{F}}_p)$ may be canonically identified (additively) with the abelian group of units in $\overline{\mathbb{F}}_p[\![x]\!]$ with constant coefficient 1 under multiplication. Further, the $p$-typical Witt vectors $W = W(\overline{\mathbb{F}}_p)$ splits off of the big Witt vectors. The quotient map $(1 + x\overline{\mathbb{F}}_p[\![x]\!]) \cong W_{\mathrm{big}}(\overline{\mathbb{F}}_p) \to \overline{\mathbb{F}}_p$ is given by reading off the coefficient of $x$, and this quotient map factors through $W$. The maps constructed so far fit into a diagram:

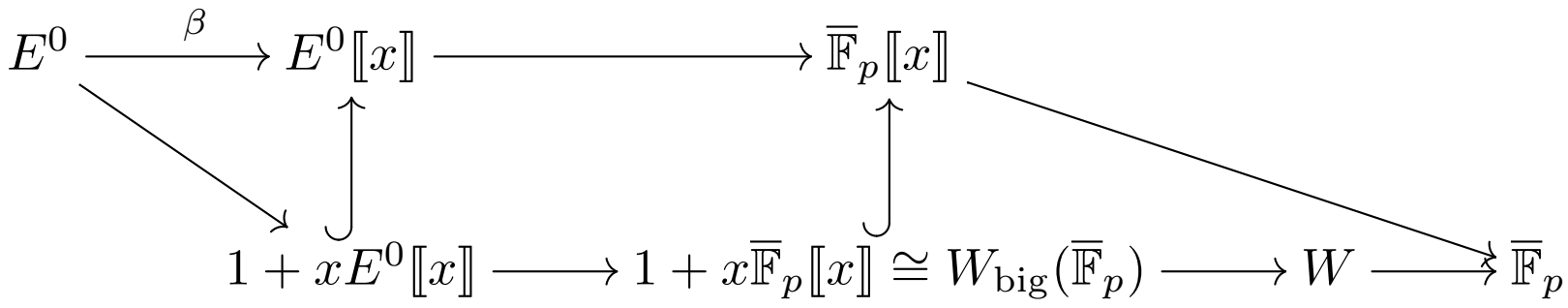

Let $\gamma\colon E^0 \to W$ be the composition of the maps appearing in the diagram. Precomposing with the inclusion $W \to E^0$, we obtain an additive endomorphism $f$ of $W$. This is a map between $p$-adically complete modules and it is the identity modulo $p$; therefore $f$ is an automorphism. The map $\alpha := f^{-1} \circ \gamma$ is therefore a section of $W \to E^0$.

Further, the maps that go into the construction of $\alpha$ are $\mathbb{G}_n$-equivariant. For $P_m$, this follows from the fact that $\mathbb{G}_n$ acts on $E$ via $\mathbb{E}_\infty$-ring maps. The transfer map is $\mathbb{G}_n$-equivariant as it is given by restriction along a map of spectra (alternatively by the formula for this transfer and the action of $\mathbb{G}_n$ on the level of characters). □

*Remark* 2.5.8. We do not know if the map $f\colon W \to W$ appearing in the proof above is the identity.

**Corollary 2.5.9.** *The inclusion $W \hookrightarrow A$ induces a split injection $H^*_{\mathrm{cts}}(\mathbb{G}_n, W) \to H^*_{\mathrm{cts}}(\mathbb{G}_n, A)$ with cokernel $H^*_{\mathrm{cts}}(\mathbb{G}_n, A^c)$.*

Theorem B has thus been reduced to the statement that $H^*_{\mathrm{cts}}(\mathbb{G}_n, A^c)$ is $p$-power torsion. This will be established in the course of the next sections.

2.6. **The proof of Theorem A assuming Theorem B.** We finish this section by explaining how to deduce Theorem A from Theorem B. The key point is that, rationally, the cohomology of the stabilizer group action on the homotopy groups of Morava $E$-theory simplifies dramatically in non-zero degrees:

**Lemma 2.6.1.** *For all $t \neq 0$ and all $s \in \mathbb{Z}$, we have $H^s_{\mathrm{cts}}(\mathbb{G}_n, \mathbb{Q} \otimes \pi_t E_n) = 0$.*

*Proof.* Recall that the (extended) Morava stabilizer group $\mathbb{G}_n$ can be described naturally as a semidirect product

$$1 \longrightarrow \mathcal{O}_D^\times \longrightarrow \mathbb{G}_n \simeq \mathcal{O}_D^\times \rtimes \mathrm{Gal}(\overline{\mathbb{F}}_p/\mathbb{F}_p) \longrightarrow \mathrm{Gal}(\overline{\mathbb{F}}_p/\mathbb{F}_p) \longrightarrow 1, \tag{2.6.2}$$

where $\mathcal{O}_D^\times$ is isomorphic to the automorphism group of our chosen formal group law $\Gamma_n$ over $\overline{\mathbb{F}}_p$. The center of $\mathcal{O}_D^\times$ is isomorphic to $\mathbb{Z}_p^\times$ and we may consider the central subgroup $\mathbb{Z}_p \subset \mathbb{Z}_p^\times \trianglelefteq \mathcal{O}_D^\times$, which we can take to be generated by the element $1 + p \in \mathbb{Z}_p^\times$. Fixing some integer $t$, the associated convergent Lyndon–Hochschild–Serre spectral sequence for continuous cohomology (Proposition 3.3.6 below) has signature

$$H^p_{\mathrm{cts}}(\mathcal{O}_D^\times/\mathbb{Z}_p, H^q_{\mathrm{cts}}(\mathbb{Z}_p, \mathbb{Q} \otimes \pi_t E_n)) \implies H^{p+q}_{\mathrm{cts}}(\mathcal{O}_D^\times, \mathbb{Q} \otimes \pi_t E_n).$$

Since the groups we take cohomology of are compact, we may move the rationalization out of the coefficients. It is thus enough to show that

$$H^q_{\mathrm{cts}}(\mathbb{Z}_p, \mathbb{Q} \otimes \pi_t E_n) = 0$$

for $t \neq 0$. To this end, we use that the generator of $\mathbb{Z}_p$ acts by multiplication by $(1 + p)^t$ on $\pi_{2t} E_n$, see for example [BB20, Section 3.2.2(c)], while the odd homotopy groups are zero. The continuous $\mathbb{Z}_p$-cohomology of $\mathbb{Q} \otimes \pi_t E_n$ is thus computed via the complex

$$\mathbb{Q} \otimes \pi_{2t} E_n \xrightarrow{(1+p)^t - 1} \mathbb{Q} \otimes \pi_{2t} E_n.$$

Since $\mathbb{Q} \otimes \pi_{2t} E_n$ is a $\mathbb{Q}$-vector space, when $t \neq 0$ the action by $(1+p)^t - 1$ is invertible, so the complex is acyclic. We then conclude by another application of the Lyndon–Hochschild–Serre spectral sequence, this time for the extension (2.6.2). □

**Proposition 2.6.3.** *Theorem B implies Theorem A.*

*Proof.* We use the Devinatz–Hopkins spectral sequence (2.3.3):

$$E_2^{s,t} \cong H^s_{\mathrm{cts}}(\mathbb{G}_n, \pi_t E_n) \implies \pi_{t-s} L_{K(n)} S^0.$$

This spectral sequence converges strongly and collapses on a finite page with a horizontal vanishing line. On the one hand, it follows that rationalization yields another strongly convergent spectral sequence

$$\mathbb{Q} \otimes E_2^{s,t} \cong H^s_{\mathrm{cts}}(\mathbb{G}_n, \mathbb{Q} \otimes \pi_t E_n) \implies \mathbb{Q} \otimes \pi_{t-s} L_{K(n)} S^0.$$

Here, the identification of the $E_2$-term uses that $\mathbb{G}_n$ is a compact group so that rationalization commutes with taking continuous cohomology; see for example [BP21, Corollary 12.9]. On the other hand, Lemma 2.6.1 implies that the rationalized Devinatz–Hopkins spectral sequence collapses on the $E_2$-page, resulting in a graded isomorphism

$$H^*_{\mathrm{cts}}(\mathbb{G}_n, \mathbb{Q} \otimes \pi_0 E_n) \cong \mathbb{Q} \otimes \pi_{-*} L_{K(n)} S^0.$$

Since the spectral sequence is multiplicative, this is in fact an isomorphism of graded rings. Theorem B combined with Lazard's theorem [Laz65] as stated in Lemma 3.8.5 identifies the left hand side as

$$H^*_{\mathrm{cts}}(\mathbb{G}_n, \mathbb{Q} \otimes \pi_0 E_n) \cong H^*_{\mathrm{cts}}(\mathbb{G}_n, \mathbb{Q} \otimes W) \cong \Lambda_{\mathbb{Q}_p}(x_1, x_2, \ldots, x_n),$$

with $x_i$ in cohomological degree $2i-1$. This gives Theorem A. □

## 3. Arithmetic prerequisites

To complete the proof of Theorem B, we must show that the complement of the split injection

$$H^*_{\mathrm{cts}}(\mathbb{G}_n, W) \to H^*_{\mathrm{cts}}(\mathbb{G}_n, A)$$

is $p^N$-torsion for some $N \geq 0$. Our proof lies entirely within the domain of $p$-adic geometry. In this section, we review some fundamentals of $p$-adic geometry and condensed mathematics and how these interact with each other. At the end of the section, we offer a summary of the proof of Theorem B.

3.1. **Adic spaces.** We offer the reader a brief summary of the necessary techniques from nonarchimedean analytic geometry, starting with Huber's category of *adic spaces* [Hub94]. This category contains all formal schemes, rigid-analytic varieties, and perfectoid spaces. For a more leisurely exposition, see the last named author's chapter in [BCKW19].

A topological ring $A$ is a *Huber ring* if it contains an open subring $A_0$ whose topology is induced by a finitely generated ideal $I \subset A_0$. (Only the topology is part of the structure of a Huber ring, not $A_0$ or $I$.) A subset $S$ of a Huber ring $A$ is *bounded* if for every open subset $U \subset A$ containing the origin, there exists another such subset $V \subset A$ such that $VS \subset U$. A single element $f \in A$ is *power-bounded* if $\{f^n\}_{n \geq 1}$ is bounded. A *Huber pair* is a pair $(A, A^+)$ consisting of a Huber ring $A$ and an open and integrally closed subring $A^+ \subset A$ whose elements are power-bounded. A *continuous valuation* on a Huber ring $A$ is a continuous multiplicative function $|\;| : A \to H \cup \{0\}$, where $H$ is a totally ordered abelian group (written multiplicatively).

The adic spectrum $\mathrm{Spa}(A, A^+)$ is the set of equivalence classes of continuous valuations satisfying $|A^+| \leq 1$. It is endowed with the topology generated by *rational subsets* of the form

$$U = U\left(\frac{f_1, \dots, f_r}{g}\right) = \left\{ |\ | \in \mathrm{Spa}(A, A^+) \;\middle|\; |f_i| \leq |g| \neq 0,\ i = 1, \dots, r \right\}$$

for $f_1, \dots, f_r \in A$ generating an open ideal and $g \in A$. Then $\mathrm{Spa}(A, A^+)$ is quasi-compact.

One defines presheaves of rings $\mathcal{O}_X^+ \subset \mathcal{O}_X$ on $X = \mathrm{Spa}(A, A^+)$ as follows. For the rational subset $U$ above, we declare that $\mathcal{O}_X(U)$ is the completion of $A[f_i/g]$ with respect to the topology in which $A_0[f_i/g]$ (with its $I$-adic topology) is an open subring, and $\mathcal{O}_X^+(U)$ is the completion of the integral closure of $A^+[f_i/g]$ in $A[f_i/g]$. Then $(\mathcal{O}_X(U), \mathcal{O}_X^+(U))$ is another Huber pair. For each $x \in \mathrm{Spa}(A, A^+)$, one gets a continuous valuation $|\ |_x$ on the local ring $\mathcal{O}_{X,x}$.

For certain classes of Huber pairs $(A, A^+)$ (including all which are considered in this article), the presheaves $\mathcal{O}_X$ and $\mathcal{O}_X^+$ are sheaves. Consider the category of triples $(X, \mathcal{O}_X, \{\ \}_x)$, where $X$ is a topological space, $\mathcal{O}_X$ is a sheaf of topological rings, and $\{|\ |_x\}$ is a collection of continuous valuations $|\ |_x$ on $\mathcal{O}_{X,x}$ for each $x \in X$. An *adic space* is an object in this category which is locally isomorphic to one arising from a Huber pair. The *integral structure sheaf* $\mathcal{O}_X^+$ is the subsheaf of $\mathcal{O}_X$ whose sections over $U$ consisting of those $f \in \mathcal{O}_X(U)$ for which $|f|_x \leq 1$ for all $x \in U$.

As a basic example, let $X = \mathrm{Spa}(\mathcal{O}_K, \mathcal{O}_K)$, where $K$ is a complete metric field containing $\mathbb{Q}_p$, and $\mathcal{O}_K$ is its ring of integers. Assuming that the residue field of $\mathcal{O}_K$ is algebraic over $\mathbb{F}_p$, the ringed space $(X, \mathcal{O}_X)$ is isomorphic to $\mathrm{Spec}\, \mathcal{O}_K$. In particular it has two points: a generic point lying in $\mathrm{Spa}(K, \mathcal{O}_K)$ (corresponding to the given metric on $K$) and a special point (which satisfies $|p| = 0$).

A *rigid-analytic space* over $K$ is an adic space over $\mathrm{Spa}(K, \mathcal{O}_K)$ that is locally isomorphic to an affinoid adic space of the form $\mathrm{Spa}(A, A^+)$, where $(A, A^+)$ obeys a certain finiteness condition. Namely, $A$ is isomorphic to a ring quotient of a Tate algebra $K\langle T_1, \dots, T_d\rangle$, and $A^+$ is equal to the subring of power-bounded elements of $A$. A special case is $A = K\langle T_1, \dots, T_d\rangle$, in which case $\mathrm{Spa}(A, A^+)$ is the closed ball of radius 1. (Chronologically, Tate's theory of rigid-analytic spaces long predates Huber's theory of adic spaces; rigid-analytic spaces as Tate defined them sit inside of adic spaces over $\mathrm{Spa}(K, \mathcal{O}_K)$ as a full subcategory, so there is no harm in thinking of them this way.)

Fix a continuous real-valued valuation $|\ |$ representing the sole point of $\mathrm{Spa}(K, \mathcal{O}_K)$, so as to fix a value of $|p|$.

**Example 3.1.1** (The rigid-analytic open ball)**.** Let $A = \mathcal{O}_K[\![T_1, \dots, T_d]\!]$, and let $\mathfrak{B}^d = \mathrm{Spa}(A, A)$, the formal $d$-dimensional unit ball over $\mathrm{Spa}(\mathcal{O}_K, \mathcal{O}_K)$. Let $B^{d,\circ}$ be the fiber of $\mathfrak{B}^d$ over the generic point of $\mathrm{Spa}(\mathcal{O}_K, \mathcal{O}_K)$; i.e., $B^{d,\circ}$ is the locus where $|p| \neq 0$. Let us observe that $B^{d,\circ}$ is exhausted by affinoid (closed) balls over $K$. For each real number $r$ of the form $r = |p|^{1/n}$ with $n = 1, 2, 3, \dots$, let $B_r^d$ be the following rational subset of $\mathfrak{B}^d$:

$$B_r^d = U\left(\frac{T_1^n, \dots, T_d^n, p}{p}\right).$$

Then $B_r^d$ is an affinoid rigid-analytic space over $K$. We have

$$B^{d,\circ} = \varinjlim_{r<1} B_r^d,$$

since for each continuous valuation on $A$ with $|p| \neq 0$, we must have $|T_i|^n \leq |p|$ for $n$ sufficiently large. Therefore $B^{d,\circ}$ is a rigid-analytic space over $K$.

For a rigid-analytic space $X/K$, there is an important distinction between $H^0(X, \mathcal{O}_X)$ and $H^0(X, \mathcal{O}_X^+)[1/p]$. When $X$ is affinoid, these agree, but in general they are quite different. In the situation of the rigid-analytic open ball, we have

$$H^0(B^{d,\circ}, \mathcal{O}_X^+)[1/p] = \mathcal{O}_K[\![T_1,\dots,T_d]\!][1/p]$$

is the ring of power series in $K[\![T_1,\dots,T_d]\!]$ which are bounded on the open unit ball. Whereas, $H^0(B^{d,\circ}, \mathcal{O}_X)$ is the much larger ring of power series in $K[\![T_1,\dots,T_d]\!]$ which converge on the open unit ball.

3.2. **Condensed mathematics.** When $G$ is a topological group and $M$ is a topological abelian group admitting a continuous $G$-action, continuous cohomology groups $H^i_{\mathrm{cts}}(G,M)$ are defined using the complex of continuous cocycles. At this level of generality, however, one has no abelian category of topological abelian groups, and so $M \mapsto H^i_{\mathrm{cts}}(G,M)$ has no interpretation as the $i$th derived functor of fixed points $M \mapsto M^G$.

The language of condensed mathematics [CS] is well-suited to resolve this issue in all contexts arising in this article. We quickly review the main points of the theory, ignoring set-theoretic issues throughout.[4] The pro-étale site of a point $*_{\text{proét}}$ is the category of profinite sets with jointly surjective continuous maps as the covers. A *condensed set* is a sheaf of sets on $*_{\text{proét}}$. Similarly there are condensed groups, rings, etc. There is a functor $X \mapsto \underline{X}$ from topological spaces/groups/rings to condensed sets/groups/rings, via $\underline{X}(S) = C_{\mathrm{cts}}(S,X)$, meaning continuous maps $S \to X$. This functor is fully faithful when restricted to compactly generated topological spaces [CS, Proposition 1.7].

Let Cond(Ab) be the category of condensed abelian groups. Then Cond(Ab) is an abelian category containing all limits and colimits [CS, Theorem 1.10]. It has a symmetric monoidal tensor product $\mathcal{M} \otimes \mathcal{N}$ and an internal Hom-functor $\underline{\mathrm{Hom}}(\mathcal{M},\mathcal{N})$ related by the adjunction

$$\mathrm{Hom}(\mathcal{P}, \underline{\mathrm{Hom}}(\mathcal{M},\mathcal{N})) \cong \mathrm{Hom}(\mathcal{P}\otimes\mathcal{M},\mathcal{N}).$$

The forgetful functor Cond(Ab) $\to$ Cond(Set) has a left adjoint $\mathcal{X} \mapsto \mathbb{Z}[\mathcal{X}]$, the "free condensed abelian group on $\mathcal{X}$". In the case $\mathcal{X} = \underline{S}$ for $S$ profinite, we abuse notation and write $\mathbb{Z}[S]$ for $\mathbb{Z}[\underline{S}]$.

The category Cond(Ab) has enough projectives [CS, Theorem 2.2], so we can form its derived category $D$(Cond(Ab)). Then $D$(Cond(Ab)) admits a derived tensor product $\otimes$ and a derived internal hom functor RHom satisfying the usual adjunction relation.

In the language of condensed mathematics, there is a notion of "solidity" which is akin to completeness. For a profinite set $S = \varprojlim S_i$ with each $S_i$ finite, the free solid abelian group on $S$ is defined as

$$\mathbb{Z}[S]^{\blacksquare} = \varprojlim \underline{\mathbb{Z}[S_i]},$$

where $\mathbb{Z}[S_i]$ is the free abelian group on $S_i$, considered as a discrete topological group. A *solid abelian group* is a condensed abelian group $\mathcal{M}$ such that for all profinite sets $S$, any morphism $\underline{S} \to \mathcal{M}$ extends uniquely to a morphism $\mathbb{Z}[S]^{\blacksquare} \to \mathcal{M}$. Let Solid denote the category of solid abelian groups; by [CS, Theorem 5.8], Solid is closed under all limits and colimits in Cond(Ab). Then Solid is an abelian subcategory of Cond(Ab).

The category Solid has enough projectives: products $\prod_I \mathbb{Z}$, where $I$ is any set, form a family of compact projective generators [CS, Theorem 5.8(i)]; alternatively one may take the free solid abelian groups $\mathbb{Z}[S]^{\blacksquare}$ as generators [CS, Corollary 5.5]. Therefore we can form the derived category $D$(Solid), from which the functor into $D$(Cond(Ab)) is fully faithful [CS, Theorem 5.8(ii)]. For an object $\mathcal{C}$ of $D$(Cond(Ab)), the following are equivalent (loc. cit.):

(1) $\mathcal{C}$ lies in the essential image of $D(\mathrm{Solid}) \to D(\mathrm{Cond(Ab)})$.
(2) $H^i(\mathcal{C})$ is a solid abelian group for all $i \in \mathbb{Z}$.

[4] They can be dealt with as in [CS, Lecture I].

(3) For all profinite sets $S$, the natural map

$$\mathrm{RHom}(\mathbb{Z}[S]^{\blacksquare}, \mathcal{C}) \to \mathrm{RHom}(\mathbb{Z}[S], \mathcal{C})$$

is an isomorphism.

The inclusion Solid $\subset$ Cond(Ab) admits a left adjoint $\mathcal{M} \mapsto \mathcal{M}^{\blacksquare}$, the "solidification functor"; this is the unique colimit-preserving extension of $\mathbb{Z}[S] \mapsto \mathbb{Z}[S]^{\blacksquare}$ [CS, Theorem 5.8]. There is a unique symmetric monoidal tensor product $\mathcal{M} \otimes^{\blacksquare} \mathcal{N}$ such that $\mathcal{M} \mapsto \mathcal{M}^{\blacksquare}$ is symmetric monoidal [CS, Theorem 6.2].

The following lemma shows how solidity and completeness are related.

**Lemma 3.2.1.** *Let $M$ be a topological abelian group which is separated and complete for a linear topology. Then $\underline{M}$ is a solid abelian group.*

*Proof.* The hypothesis means there is a directed system of open subgroups $M_n \subset M$ inducing the topology on $M$, such that the map $M \to \varprojlim M/M_n$ is an isomorphism. Then $\underline{M} \cong \varprojlim \underline{M/M_n}$. Each $\underline{M/M_n}$ is discrete, hence solid, hence the limit is also solid. $\square$

3.3. **The solid perspective on continuous cohomology.** Now suppose $\mathcal{G}$ is a condensed group. A $\mathcal{G}$-action on a condensed abelian group $\mathcal{M}$ is a morphism $\mathcal{G} \times \mathcal{M} \to \mathcal{M}$ satisfying the usual axioms.

**Lemma 3.3.1.** *Let $M$ be a topological abelian group which is separated and complete for a linear topology. Let $G$ be a profinite group which acts continuously on $M$. Then $\underline{M}$ is a solid abelian group with $\underline{G}$-action.*

*Proof.* Keeping in mind Lemma 3.2.1, this follows from the fact that the functor $X \mapsto \underline{X}$ preserves limits by [CS, Proposition 1.7]. $\square$

For a profinite group $G$, let $\mathrm{Solid}_G$ be the category of solid abelian groups admitting an action of $\underline{G}$. Then $\mathrm{Solid}_G$ is an abelian category. Importantly, it has enough projectives:

**Lemma 3.3.2.** *Let $G$ be a profinite group. The solid abelian groups $\mathbb{Z}[G \times S]^{\blacksquare}$, where $S$ is a profinite set and $G$ acts by left multiplication on the first factor of $G \times S$, form a family of compact projective generators of* $\mathrm{Solid}_G$*.*

*Proof.* We first claim that $\mathbb{Z}[G \times S]^{\blacksquare}$ is projective. Indeed, let $\pi\colon \mathcal{M} \to \mathbb{Z}[G \times S]^{\blacksquare}$ be a surjection in $\mathrm{Solid}_G$. Then the composition $\mathcal{M} \to \mathbb{Z}[G \times S]^{\blacksquare} \to \mathbb{Z}[S]^{\blacksquare}$ is surjective. Since $\mathbb{Z}[S]^{\blacksquare}$ is projective in Solid, there is a section $\mathbb{Z}[S]^{\blacksquare} \to \mathcal{M}$. Consider the composition

$$\underline{G} \times \mathbb{Z}[S]^{\blacksquare} \to \underline{G} \times \mathcal{M} \stackrel{\alpha}{\to} \mathcal{M}$$

in the category of condensed sets with $\underline{G}$-action ($\alpha$ is the action map). Since $\mathcal{M}$ is solid, this morphism extends uniquely to a morphism $\mathbb{Z}[G \times S]^{\blacksquare} \to \mathcal{M}$ in $\mathrm{Solid}_G$, which is the required section to $\pi$.

Now suppose $\mathcal{M}$ is an object in $\mathrm{Solid}_G$. Forgetting the $G$-structure, $\mathcal{M}$ receives a surjective morphism from a colimit over $i \in I$ of free solid abelian groups $\mathbb{Z}[S_i]^{\blacksquare}$, each of which admits a section $\mathbb{Z}[S_i]^{\blacksquare} \to \mathcal{M}$ in Solid. Repeating the argument of the last paragraph, this can be upgraded to a morphism $\mathbb{Z}[G \times S_i]^{\blacksquare} \to \mathcal{M}$ in $\mathrm{Solid}_G$, which after taking the colimit over $i$ becomes a surjective map from a colimit of projective objects of the form $\mathbb{Z}[G \times S]^{\blacksquare}$. $\square$

Lemma 3.3.2 allows us to define the derived category $D(\mathrm{Solid}_G)$: this is the $\infty$-category formed from the category of complexes of projective objects of $\mathrm{Solid}_G$ by inverting quasi-isomorphisms. For an object $\mathcal{M}$ of $D(\mathrm{Solid}_G)$, there is a functor $D(\mathrm{Solid}) \to D(\mathrm{Solid}_G)$ given by $\mathcal{P} \mapsto \mathcal{P} \otimes \mathcal{M}$; this functor admits a right adjoint $D(\mathrm{Solid}_G) \to D(\mathrm{Solid})$, which we notate as $\mathcal{N} \mapsto \underline{\mathrm{RHom}}_G(\mathcal{M}, \mathcal{N})$.

**Notation 3.3.3** (Continuous homotopy fixed points)**.** Let $G$ be a profinite group. We write

$$\begin{array}{rcl} D(\mathrm{Solid}_G) & \to & D(\mathrm{Solid}) \\ \mathcal{M} & \mapsto & \mathcal{M}^{hG} \end{array}$$

for the functor which sends $\mathcal{M}$ to $\underline{\mathrm{RHom}}_G(\mathbb{Z}^{\blacksquare}, \mathcal{M})$, where $\mathbb{Z}^{\blacksquare}$ represents the trivial $G$-action on the discrete group $\mathbb{Z}$. For $i \in \mathbb{Z}$ we write $H^i(G, \mathcal{M})$ for the $i$th cohomology group of $\mathcal{M}^{hG}$, an object of Solid.

As a special case, suppose $M$ is a topological abelian group which is separated and complete for the topology generated by subgroups $M_n$, and $G$ is a profinite group which acts continuously on $M$. On the one hand, viewing $\underline{M}$ as an object in the heart of $D(\mathrm{Solid}_G)$, we have the solid abelian groups $H^i(G, \underline{M})$. On the other hand, we have the complex of continuous cocycles $C_{\mathrm{cts}}(G, M)$ which computes $H^i_{\mathrm{cts}}(G, M)$. We can give each $C^i_{\mathrm{cts}}(G, M)$ the topology generated by the $C^i_{\mathrm{cts}}(G, M_n)$; then the $C^i_{\mathrm{cts}}(G, M)$ are separated and complete, and thereby we can view $C_{\mathrm{cts}}(G, M)$ as an object of $D(\mathrm{Solid})$, and each $H^i_{\mathrm{cts}}(G, M)$ as an object of Solid. The following lemma shows that the two points of view are compatible.

**Lemma 3.3.4.** *Let $M$ be an abelian group which is separated and complete for a linear topology, and let $G$ be a profinite group which acts continuously on $M$. There is a natural isomorphism in* Solid*:*

$$H^i(G, \underline{M}) \cong H^i_{\mathrm{cts}}(G, M).$$

*Proof.* The standard projective resolution of the trivial $G$-module $\mathbb{Z}$,

$$\cdots \to \mathbb{Z}[G \times G] \to \mathbb{Z}[G] \to \mathbb{Z}$$

has an extra degeneracy and thus is homotopic to zero. Applying the additive solidification functor therefore yields an exact complex

$$\cdots \to \mathbb{Z}[G \times G]^{\blacksquare} \to \mathbb{Z}[G]^{\blacksquare} \to \mathbb{Z}^{\blacksquare}. \tag{3.3.5}$$

For each $n \geq 1$, the free solid abelian group $\mathbb{Z}[G^n]^{\blacksquare}$ with $G$ acting diagonally on the left is a projective object of $D(\mathrm{Solid}_G)$ by Lemma 3.3.2, so (3.3.5) is in fact a projective resolution of the trivial module $\mathbb{Z}^{\blacksquare}$ in $\mathrm{Solid}_G$. Therefore $H^i(G, \underline{M})$ is computed by the complex in $D(\mathrm{Solid})$:

$$\underline{\mathrm{Hom}}(\mathbb{Z}[G]^{\blacksquare}, \underline{M}) \to \underline{\mathrm{Hom}}(\mathbb{Z}[G \times G]^{\blacksquare}, \underline{M}) \to \cdots$$

Since each $\mathbb{Z}[G^i]^{\blacksquare}$ is free, the $i$th term of the complex is isomorphic to $\underline{\mathrm{Hom}}(G^i, \underline{M})$, which is in turn isomorphic to $C^i_{\mathrm{cts}}(G, M)$, the solid group of continuous cocycles. ☐

Lemma 3.3.4 shows that as long as $M$ is separated and complete for a linear topology, the continuous cohomology $H^i_{\mathrm{cts}}(G, M)$ admits a good interpretation in terms of homological algebra. As an application, we establish a Lyndon–Hochschild–Serre spectral sequence in this context.

**Proposition 3.3.6** (Lyndon–Hochschild–Serre spectral sequence)**.** *Let $\mathcal{M}$ be an object of $D(\mathrm{Solid}_G)$ Let $H \subset G$ be a closed normal subgroup. Then for each $j \geq 0$ there is a natural interpretation of $H^j(H, \mathcal{M})$ as an object of $\mathrm{Solid}_{G/H}$, and there is a spectral sequence with $E_2$ page $H^i(G/H, H^j(H, \mathcal{M}))$ which converges to $H^{i+j}(G, \mathcal{M})$.*

*In particular, if $M$ is an abelian group which is separated and complete for a linear topology, and $G$ acts continuously on $M$, then there is a spectral sequence with $E_2$ page $H^i_{\mathrm{cts}}(G/H, H^j_{\mathrm{cts}}(H, M))$ which converges to $H^{i+j}_{\mathrm{cts}}(G, M)$.*

*Remark* 3.3.7. The category $\mathrm{Solid}_G$ does not necessarily have enough injectives, and so one cannot compute $\mathcal{M} \mapsto \mathcal{M}^{hG}$ in the usual way as the right derived functor of a left exact functor via an injective resolution; thus the Grothendieck spectral sequence is not directly available in this setting.

It is possible to extend the Lyndon–Hochschild–Serre spectral sequence to the case where $G$ is a general locally compact group, but one needs additional hypotheses on the representation, see for instance [CW74, Proposition 5].

*Proof.* In preparation for the proof, we record some functorial properties of the assignment $G \mapsto D(\mathrm{Solid}_G)$. Suppose $G \to G'$ is a continuous homomorphism between profinite groups. Let us write $\mathrm{Res}_G^{G'}$ for the restriction of scalars functor $D(\mathrm{Solid}_{G'}) \to D(\mathrm{Solid}_G)$. Then $\mathrm{Res}_G^{G'}$ admits both a left adjoint (extension of scalars, or induction)

$$\mathrm{Ind}_G^{G'} \mathcal{M} = \mathcal{M} \otimes_{\mathbb{Z}[G]^{\blacksquare}} \mathbb{Z}[G']^{\blacksquare}$$

and a right adjoint (coextension of scalars, or coinduction)

$$\mathrm{coInd}_G^{G'} \mathcal{M} = \underline{\mathrm{RHom}}_G(\mathbb{Z}[G']^{\blacksquare}, \mathcal{M}).$$

In the coextension of scalars, the $G'$-action appears as the right action of $\mathbb{Z}[G']^{\blacksquare}$ on itself. Note that the continuous homotopy fixed points functor $\mathcal{M} \mapsto \mathcal{M}^{hG}$ is a special case of coinduction, in the case $G' = 1$.

Now suppose $H \subset G$ is a closed normal subgroup as in the lemma, so that $G \to G/H$ is a continuous homomorphism of profinite groups. The factorization $\mathrm{Res}_1^G = \mathrm{Res}_1^{G/H} \circ \mathrm{Res}_{G/H}^G$ induces a factorization of right adjoints, so that for all objects $\mathcal{M}$ of $D(\mathrm{Solid}_G)$, we have $\mathcal{M}^{hG} \cong (\mathrm{coInd}_G^{G/H} \mathcal{M})^{h(G/H)}$.

Consider the coextension of scalars $\mathrm{coInd}_G^{G/H} \mathcal{M}$. When its $G/H$-module structure is forgotten, we claim it is isomorphic to $\mathcal{M}^{hH}$. The key observation is that the map

$$\mathrm{Ind}_H^G \mathbb{Z}^{\blacksquare} = \mathbb{Z}^{\blacksquare} \otimes_{\mathbb{Z}[H]^{\blacksquare}} \mathbb{Z}[G]^{\blacksquare} \to \mathbb{Z}[G/H]^{\blacksquare}$$

is an isomorphism in $D(\mathrm{Solid}_G)$. To see this, consider the resolution of $G/H$ by free $G$-sets:

$$\cdots G \times H^2 \mathrel{\substack{\to\\ \to\\ \to}} G \times H \rightrightarrows G.$$

The functor $S \mapsto \mathbb{Z}[S]^{\blacksquare}$ is left adjoint, and thus the colimit of the diagram

$$\cdots \to \mathbb{Z}[G \times H^2]^{\blacksquare} \to \mathbb{Z}[G \times H]^{\blacksquare} \to \mathbb{Z}[G]^{\blacksquare} \tag{3.3.8}$$

is $\mathbb{Z}[G/H]^{\blacksquare}$. On the other hand, (3.3.8) is the result of tensoring the free resolution of $\mathbb{Z}^{\blacksquare}$ by $\mathbb{Z}[H]^{\blacksquare}$-modules with $\mathbb{Z}[G]^{\blacksquare}$, so its colimit is isomorphic to $\mathbb{Z}^{\blacksquare} \otimes_{\mathbb{Z}[H]^{\blacksquare}} \mathbb{Z}[G]^{\blacksquare}$.

Turning to the coinduction $\mathrm{coInd}_G^{G/H} \mathcal{M}$, we find that its underling solid complex is

$$\begin{aligned} \mathrm{Res}_e^{G/H} \mathrm{coInd}_G^{G/H} \mathcal{M} &\cong \underline{\mathrm{RHom}}_G(\mathbb{Z}[G/H]^{\blacksquare}, \mathcal{M}) \\ &\cong \underline{\mathrm{RHom}}_G(\mathrm{Ind}_H^G \mathbb{Z}^{\blacksquare}, \mathcal{M}) \\ &\cong \underline{\mathrm{RHom}}_H(\mathbb{Z}^{\blacksquare}, \mathrm{Res}_H^G \mathcal{M}) \\ &\cong (\mathrm{Res}_H^G \mathcal{M})^{hH}. \end{aligned}$$

Our observation justifies a slight abuse of notation: if $\mathcal{M}$ is an object of $D(\mathrm{Solid}_G)$, we simply write $\mathcal{M}^{hH}$ for $\mathrm{coInd}_G^{G/H} \mathcal{M}$, understanding that there is a residual $G/H$-action. Then our factorization of $hG$ becomes:

$$\mathcal{M}^{hG} \simeq (\mathcal{M}^{hH})^{h(G/H)}.$$

We rely on the following general principle for constructing spectral sequences [Lur17, Proposition 1.2.2.14]: suppose $\mathcal{C}$ is a stable $\infty$-category equipped with a t-structure. Let

$$\cdots \to X(-1) \to X(0) \to X(1) \to \cdots$$

be a filtered object of $\mathcal{C}$ such that $X(n) \cong 0$ for $n \ll 0$. Assume $\mathcal{C}$ admits sequential colimits, and write $X = \operatorname{colim} X(n)$. Write $\mathrm{Gr}\, X(n)$ for the $n$th graded piece, meaning the cofiber of $X(n-1) \to X(n)$. Then there is a spectral sequence with $E_2$ page $\pi_p \operatorname{Gr} X(q)$ converging to $\pi_{p+q}X$. (Note that [Lur17, Proposition 1.2.2.14] produces an $E_1$-spectral sequence; we have changed the coordinates in order to produce the more familiar $E_2$-spectral sequence appearing in the proposition.)

In the situation of the proposition, we apply the above construction in the category $D(\mathrm{Solid})$, to the filtered object

$$X(n) := (\tau_{\geq -n}(\mathcal{M}^{hH}))^{h(G/H)}.$$

Since $h(G/H)$ is exact, the graded pieces of this object are $\operatorname{Gr} X(q) = H^q(H, M)^{h(G/H)}$, with homotopy groups $\pi_{-p} \operatorname{Gr} X(q) = H^p(G/H, H^q(H, M))$.

It remains to verify that $\operatorname{colim} X(n) \cong (\mathcal{M}^{hH})^{h(G/H)} \cong \mathcal{M}^{hG}$. To simplify notation, set $N := \mathcal{M}^{hH}$. Write $F(-) := (-)^{h(G/H)}$, so that $X(n) = F(\tau_{\geq -n}N)$ and $F(N) \simeq \mathcal{M}^{hG}$. Applying the exact functor $\operatorname{colim}_n F(-)$ to the natural truncation sequences for $N$, we obtain a cofiber sequence

$$\operatorname{colim}_n F(\tau_{\geq -n}N) \to \operatorname{colim}_n F(N) \to \operatorname{colim}_n F(\tau_{< -n}N).$$

We need to verify that $\operatorname{colim}_n F(\tau_{< -n}N) \simeq 0$. To this end, it suffices to show that

$$\pi_m \operatorname{colim}_n F(\tau_{< -n}N) \cong \operatorname{colim}_n \pi_m F(\tau_{< -n}N) = 0$$

for all $m \in \mathbb{Z}$. But $F$ preserves coconnectivity, so $\pi_m F(\tau_{< -n}N) = 0$ whenever $-n \leq m$, so certainly also in the colimit when $n \to \infty$. This implies the claim and finishes the proof. ☐

3.4. **$p$-adic solid complexes.** Here, we gather some results regarding the notion of derived $p$-completeness (for $p$ a prime number), in the contexts of both $D(\mathrm{Ab})$ and $D(\mathrm{Cond}(\mathrm{Ab}))$. Many of these results could be extended to the setting of derived $I$-completeness of $R$-modules, where $I$ is a finitely generated ideal in a ring $R$. The main result Lemma 3.4.7 below is that $D(\mathrm{Solid})$ contains the category $D(\mathrm{Ab})_p$ of derived $p$-complete complexes as a full tensor subcategory.

First we discuss the notion of derived $p$-completeness for general stable $\infty$-categories and then specialize to the case at hand; this is a straightforward extension of the notion of $p$-completeness for $D(\mathrm{Ab})$, for which some standard results are found in [Aut, Tag 091N].

**Definition 3.4.1.** Let $\mathcal{C}$ be a stable $\infty$-category. An object $A \in \mathcal{C}$ of $\mathcal{C}$ is *$p$-complete* (the word "derived" being understood) if $\underline{\operatorname{Hom}}(B, A) \cong 0$ for all $B \in \mathcal{C}$ such that $B/p \cong 0$. Let $\mathcal{C}_p$ be the full subcategory of $\mathcal{C}$ consisting of $p$-complete objects.

**Example 3.4.2.** This definition makes it clear that $D(\mathrm{Ab})_p$ is the Bousfield localization of $D(\mathrm{Ab})$ at the object $\mathbb{Z}/p\mathbb{Z}$, as described in Section 2. It suffices to check the condition in the definition with the object $B = \mathbb{Z}[1/p]$. Consequently, the inclusion $D(\mathrm{Ab})_p \subset D(\mathrm{Ab})$ has a left adjoint $A \mapsto \hat{A}$, namely:

$$\hat{A} = \varprojlim A/p^n,$$

where $A/p^n$ is the cofiber of $p^n \colon A \to A$. There is a unique symmetric monoidal operation $\hat{\otimes}$ on $D(\mathrm{Ab})_p$ which makes $A \mapsto \hat{A}$ a symmetric monoidal functor. (Here and elsewhere we write $\otimes$ instead of $\otimes^{\mathbb{L}}$ for the symmetric monoidal operation on $D(\mathrm{Ab})$.)

Specializing these notions to $D(\mathrm{Cond}(\mathrm{Ab}))$, we obtain:

**Lemma 3.4.3.** *For the category $D(\mathrm{Cond}(\mathrm{Ab}))_p$ of $p$-adically complete complexes of condensed abelian groups, we have:*

1. *$D(\mathrm{Cond}(\mathrm{Ab}))_p$ is closed under limits in $D(\mathrm{Cond}(\mathrm{Ab}))$.*
2. *For any $A \in D(\mathrm{Cond}(\mathrm{Ab}))_p$, if $A/p \cong 0$, then $A \cong 0$.*
3. *The inclusion $D(\mathrm{Cond}(\mathrm{Ab}))_p \to D(\mathrm{Cond}(\mathrm{Ab}))$ admits a left adjoint $A \mapsto A^{\wedge}_p$, namely:*

$$A^{\wedge}_p = \varprojlim A/p^n.$$

*Proof.* For part (1): if $A_i$ is a diagram in $D(\mathrm{Cond}(\mathrm{Ab}))_p$ and $B \in D(\mathrm{Cond}(\mathrm{Ab}))$ with $B/p \cong 0$, then $\underline{\mathrm{Hom}}(B, \lim A_i) \cong \lim \underline{\mathrm{Hom}}(B, A_i) \cong 0$, since the $A_i$'s are $p$-complete. In other words, $\lim A_i$ lies in $D(\mathrm{Cond}(\mathrm{Ab}))_p$.

For part (2): $A/p \cong 0$ means that multiplication by $p$ is an isomorphism on $A$, so that the canonical map $\underline{\mathrm{Hom}}(\underline{\mathbb{Z}[1/p]}, A) \to \underline{\mathrm{Hom}}(\underline{\mathbb{Z}}, A) \cong A$ is an isomorphism. But $\underline{\mathrm{Hom}}(\underline{\mathbb{Z}[1/p]}, A) \cong 0$ by definition of derived $p$-completeness, hence $A \cong 0$.

For part (3): it is enough to show that: (a) $A^\wedge_p$ is $p$-adically complete, and (b) if $A$ is $p$-adically complete, then $A \to A^\wedge_p$ is an isomorphism.

For (a): We first claim $A/p^n$ is $p$-complete for all $n$. Indeed, if $B/p \cong 0$, then multiplication by $p^n$ is an isomorphism on $B$, and thus on $\underline{\mathrm{Hom}}(B, A)$. As $D(\mathrm{Cond}(\mathrm{Ab}))_p$ is closed under limits, it must contain $A^\wedge_p$ as well.

For (b): Suppose $A$ is $p$-adically complete, and let $F$ be the fiber of $A \to A^\wedge_p$. Then $F$ is also $p$-adically complete (being a limit in $D(\mathrm{Cond}(\mathrm{Ab}))_p$), and so by part (2) it suffices to show that $F/p \cong 0$. But this is clear, as $F$ is the limit of multiplication by $p$ on $A$. $\square$

At this point we start relating the notions of $p$-completeness in $D(\mathrm{Ab})$ and $D(\mathrm{Cond}(\mathrm{Ab}))$. First we need a notion of discreteness for objects of the latter category.

**Definition 3.4.4.** Let $\Gamma_* \colon \mathrm{Cond}(\mathrm{Ab}) \to \mathrm{Ab}$ be the functor which evaluates an object on the point: $\Gamma_*(A) = A(*)$. Then $\Gamma_*$ is right adjoint to the pullback functor $\Gamma^*$ which sends an abelian group $A$ to the constant condensed abelian group $\underline{A}$. Both $\Gamma^*$ and $\Gamma_*$ are exact (the latter because all covers of $*$ have a section), so we continue to write $\Gamma^*$ and $\Gamma_*$ for their derived functors. Finally, note that, by construction, the unit $\mathrm{id} \to \Gamma_*\Gamma^*$ of the adjunction is an equivalence.

For $A \in D(\mathrm{Cond}(\mathrm{Ab}))$ we let $A^\delta = \Gamma^*\Gamma_*(A)$ be the "discretization" of $A$. The counit of the adjunction gives rise to a natural map $A^\delta \to A$. If this map is an isomorphism we say that $A$ is *discrete.*

**Lemma 3.4.5.** *Discretization enjoys the following properties:*

*(1) Discretization is $t$-exact with respect to the standard $t$-structures on the derived categories $D(\mathrm{Ab})$ and $D(\mathrm{Cond}(\mathrm{Ab}))$, respectively. That is, if $A$ lies in the heart $D(\mathrm{Cond}(\mathrm{Ab}))^\heartsuit \simeq \mathrm{Cond}(\mathrm{Ab})$, then $A^\delta \in D(\mathrm{Cond}(\mathrm{Ab}))^\heartsuit$ as well.*

*(2) The subcategory of discrete objects in $D(\mathrm{Cond}(\mathrm{Ab}))$ is closed under fibers, cofibers and retracts.*

*(3) If $A \in D(\mathrm{Ab})$, then $\Gamma^*A$ is discrete.*

*(4) An object $A \in D(\mathrm{Cond}(\mathrm{Ab}))$ has the property that, for all $i \in \mathbb{Z}$, $H^i(A)$ is discrete if and only if $A$ is discrete.*

*(5) The functor $\Gamma^*$ is symmetric monoidal.*

*Proof.* (1) Since $\Gamma^*$ is given by pulling back sheaves, it is $t$-exact. Since $\Gamma_*$ is corepresented by $\underline{\mathbb{Z}}$, which is projective as a condensed abelian group, $\Gamma_*$ is also $t$-exact. For (2), this follows from the fact that both $\Gamma^*$ and $\Gamma_*$ are exact, i.e., they both commute with fibers, cofibers and retracts. For (3), this is equivalent to the statement that $\Gamma^*$ is fully faithful. This in turn is the same as the statement that the unit of the adjunction is an equivalence, which we noted above. Statement (4) follows from the previous items. Finally, since $\Gamma^*$ is a pullback functor between derived categories of sheaves on sites, it is symmetric monoidal. $\square$

**Definition 3.4.6.** An object $A \in D(\mathrm{Cond}(\mathrm{Ab}))$ is called *$p$-adic* if it is $p$-complete and $A/p$ is discrete.

**Lemma 3.4.7.** *Consider the functor*

$$\begin{array}{rcl} D(\mathrm{Ab}) & \to & D(\mathrm{Cond}(\mathrm{Ab})) \\ A & \mapsto & (\Gamma^*A)^\wedge_p = \varprojlim \Gamma^*(A)/p^n. \end{array}$$

(1) *$A \mapsto (\Gamma^* A)^\wedge_p$ factors as $D(\mathrm{Ab}) \to D(\mathrm{Ab})_p \to D(\mathrm{Cond}(\mathrm{Ab}))_p \to D(\mathrm{Cond}(\mathrm{Ab}))$. In this composition, the outer functors are the natural ones, and the inner functor $D(\mathrm{Ab})_p \to D(\mathrm{Cond}(\mathrm{Ab})_p$ is fully faithful, with essential image consisting of the $p$-adic objects.*

(2) *The $p$-adic objects in $D(\mathrm{Cond}(\mathrm{Ab}))_p$ are solid; thus $A \mapsto (\Gamma^* A)^\wedge_p$ factors through a functor $D(\mathrm{Ab})_p \to D(\mathrm{Solid})$.*

(3) *The functor $D(\mathrm{Ab})_p \to D(\mathrm{Solid})$ is symmetric monoidal with respect to $\hat{\otimes}$ and $\otimes^{\blacksquare}$.*

*Proof.* (1) Since $(\Gamma^* A)^\wedge_p$ is $p$-complete by construction, the functor factors through $D(\mathrm{Cond}(\mathrm{Ab}))_p$. To see that it also factors through the functor $D(\mathrm{Ab}) \to D(\mathrm{Ab})_p$ sending $A$ to $A^\wedge_p$, it suffices to show that the canonical map

$$(\Gamma^* A)^\wedge_p \to (\Gamma^*(A^\wedge_p))^\wedge_p$$

is an isomorphism. Since the source and target are both $p$-complete, we can verify this after taking the quotient by $p$, see Lemma 3.4.3(2). We are reduced to showing that

$$(\Gamma^* A)/p \to (\Gamma^*(A^\wedge_p))/p$$

is an isomorphism. Since $\Gamma^*$ is exact, the map is isomorphic to $\Gamma^*(A/p) \to \Gamma^*(A^\wedge_p/p)$, which is indeed an isomorphism because $A/p \to A^\wedge_p/p$ is.

Next we observe that $(\Gamma^* A)^\wedge_p$ is $p$-adic. Indeed, it is $p$-complete by construction, while $(\Gamma^* A)^\wedge_p/p \cong \Gamma^*(A/p)$ is discrete.

Conversely, if $A \in D(\mathrm{Cond}(\mathrm{Ab}))$ is $p$-adic, then $\Gamma_* A$ is an object of $D(\mathrm{Ab})$. We claim that the counit of the adjunction between $\Gamma^*$ and $\Gamma_*$ induces an isomorphism $(\Gamma^*(\Gamma_* A))^\wedge_p = (A^\delta)^\wedge_p \to A^\wedge_p \xrightarrow{\sim} A$. Once again this may by checked after quotienting by $p$. And indeed $A^\delta/p \to A/p$ is an isomorphism, because $A/p$ is discrete.

(2) The object $(\Gamma^* A)^\wedge_p$ is solid because it is a limit of discrete objects $(\Gamma^* A)/p^n$, which are all solid.

(3) Since $\Gamma^*$ is symmetric monoidal by Lemma 3.4.5(5), and $p$-adic completion is lax symmetric monoidal, the functor $A \mapsto (\Gamma^* A)^\wedge_p$ from $D(\mathrm{Ab})_p$ into $D(\mathrm{Cond}(\mathrm{Ab}))$ is lax symmetric monoidal. That is, for all objects $A$ and $B$ of $D(\mathrm{Ab})_p$, there is a morphism in $D(\mathrm{Cond}(\mathrm{Ab}))$:

$$(\Gamma^* A)^\wedge_p \otimes (\Gamma^* B)^\wedge_p \to (\Gamma^*(A \hat{\otimes} B))^\wedge_p. \tag{3.4.8}$$

We now apply the solidification functor $\mathcal{M} \to \mathcal{M}^{\blacksquare}$ to (3.4.8), using the facts that (a) solidification is symmetric monoidal and (b) each of $(\Gamma^* A)^\wedge_p$, $(\Gamma^* B)^\wedge_p$, and $(\Gamma^*(A \hat{\otimes} B))^\wedge_p$ are already solid, so $\mathcal{M}^{\blacksquare} \cong \mathcal{M}$ for each of them. The result is a morphism

$$(\Gamma^* A)^\wedge_p \otimes^{\blacksquare} (\Gamma^* B)^\wedge_p \to (\Gamma^*(A \hat{\otimes} B))^\wedge_p \tag{3.4.9}$$

which expresses the lax symmetric monoidal structure on the functor $D(\mathrm{Ab})_p \to D(\mathrm{Solid})$.

The claim in (3) is that (3.4.9) is an isomorphism. Since the tensor products appearing on either side commute with colimits (they are both left adjoints) and shifts, it suffices to prove the claim for $A = B = \mathbb{Z}_p$, since direct sums of $\mathbb{Z}_p$ and its shifts generate $D(\mathrm{Ab})_p$. Since $(\Gamma^* \mathbb{Z}_p)^\wedge_p = \mathbb{Z}_p$ (considered as a condensed abelian group), and the natural map $\mathbb{Z}_p \hat{\otimes} \mathbb{Z}_p \to \mathbb{Z}_p$ is an isomorphism in $D(\mathrm{Ab})_p$, we are reduced to the claim that $\mathbb{Z}_p \otimes^{\blacksquare} \mathbb{Z}_p \to \mathbb{Z}_p$ is an isomorphism, which is explained in [CS, Example 6.4]. □

Given a derived $p$-complete abelian group $A$ with a continuous action of a profinite group $G$, Lemma 3.3.4(2) implies that we can consider continuous cohomology as the cohomology groups of the functor $(-)^{hG}$ applied to $(\Gamma^* A)^\wedge_p$. We get the same identification for complexes of derived $p$-complete continuous $G$-modules.

Let $D(\mathrm{Ab})_{p,G}$ be the full subcategory of $D(\mathrm{Solid}_G)$ consisting of objects whose underlying solid abelian group is $p$-adic in the sense of Definition 3.4.6, i.e., those in the image of $(\Gamma^*(-))^\wedge_p$. We will find the following "projection formula" useful.

**Lemma 3.4.10.** *Let $G$ be a profinite group and let $M \in D(\mathrm{Ab})_p$ and $N \in D(\mathrm{Ab})_{p,G}$ be bounded above objects. There is a canonical isomorphism in $D(\mathrm{Ab})_p$:*

$$M \otimes^{\blacksquare} N^{hG} \xrightarrow{\sim} (M \otimes^{\blacksquare} N)^{hG}.$$

*Proof.* Note that there is a canonical map from the left hand side to the right. Since both the source and the target are $p$-complete by Lemma 3.4.3(1) and Lemma 3.4.7, it is enough to establish the isomorphism after taking the cofiber of the multiplication by $p$ map. By the exactness of the functors involved, it is enough to prove the statement when $M$ is a bounded above complex of $\mathbb{F}_p$-modules. For $M = \mathbb{F}_p$, the result follows from the exactness of the functor $hG$. The same exactness implies the statement in all cases that $M$ is a perfect complex of $\mathbb{F}_p$-modules since these are precisely the objects in the thick subcategory generated by $\mathbb{F}_p$. To conclude, it is enough to show that both sides of the map commute with filtered colimits of uniformly bounded above $\mathbb{F}_p$-modules in the variable $M$. This is clear for the left hand side. For the right hand side, the compactness of $G$ implies that the homotopy fixed point functor commutes with filtered colimits of uniformly bounded above complexes, as can be checked in cohomology. □

3.5. **Sheaves of condensed abelian groups.** Let $X$ be a site, and let $\mathcal{A}$ be a sheaf of topological abelian groups on $X$. This means in particular that if $\{U_i\}$ is a cover of an object $U$ in $X$, then

$$\mathcal{A}(U) \to \prod_i \mathcal{A}(U_i) \rightrightarrows \prod_{i,j} \mathcal{A}(U_i \cap U_j)$$

is an equalizer diagram in the category of topological abelian groups.

In this situation one would like the cohomology $R\Gamma(X, \mathcal{A})$ to carry a topology, and indeed this is possible using condensed mathematics. Define a presheaf $\mathcal{A}_{\mathrm{cond}}$ of condensed abelian groups on $X$ by setting

$$\mathcal{A}_{\mathrm{cond}}(U) = \underline{\mathcal{A}(U)}$$

for every object $U$ of $X$. Then in fact $\mathcal{A}_{\mathrm{cond}}$ is a sheaf, which we call the *condensed enhancement* of $\mathcal{A}$. This follows from the fact that $A \mapsto \underline{A}$ is limit-preserving since it admits a left adjoint (see [CS, Proposition 1.7]); in particular it preserves equalizers.

From here we can form the *condensed cohomology* $R\Gamma(X, \mathcal{A}_{\mathrm{cond}}) \in D(\mathrm{Cond}(\mathrm{Ab}))$ as the derived global sections of $\mathcal{A}_{\mathrm{cond}}$. This is a condensed enhancement of the cohomology $R\Gamma(X, \mathcal{A})$ of the underlying sheaf of abelian groups, in the sense that $\Gamma_* R\Gamma(X, \mathcal{A}_{\mathrm{cond}}) \cong R\Gamma(X, \mathcal{A})$. Indeed, $\Gamma_*$ commutes with all limits and colimits and in particular simplicial limits and filtered colimits. Thus, it commutes with the computation of derived global sections $R\Gamma$.

An important example for us comes from formal schemes. Recall the notion of an *admissible* topological ring from [Aut, Tag0AHY]. For instance, if $A$ is separated and complete with respect to the topology induced by an ideal $I$, then $A$ with its $I$-adic topology is an admissible topological ring. An *ideal of definition* in an admissible topological ring is an open ideal $I$, such that every neighborhood of 0 contains $I^n$ for some $I$. An admissible topological ring always contains an ideal of definition.

For an admissible topological ring $A$ and any ideal of definition $I$, the map $\mathrm{Spec}\, A/I \to \mathrm{Spec}\, A$ factors through a homeomorphism from $\mathrm{Spec}\, A/I$ onto the subset $\mathrm{Spf}\, A \subset \mathrm{Spec}\, A$ consisting of open prime ideals. The structure sheaf $\mathcal{O}_{\mathrm{Spf}\, A}$ is defined as $\mathcal{O}_{\mathrm{Spf}\, A} = \lim_I \mathcal{O}_{\mathrm{Spec}\, A/I}$, where the limit runs over all ideals of definition for $A$. Here the limit is taken in the category of sheaves of topological rings, where each $\mathcal{O}_{\mathrm{Spec}\, A/I}$ is discrete. Finally, a *formal scheme* is a pair $(\mathfrak{X}, \mathcal{O})$ consisting of a topological space $\mathfrak{X}$ and a sheaf of topological rings $\mathcal{O}$ on $X$, such that $(\mathfrak{X}, \mathcal{O})$ is locally isomorphic to $(\mathrm{Spf}\, A, \mathcal{O}_{\mathrm{Spf}\, A})$ for an admissible topological ring.

For a formal scheme $\mathfrak{X}$, let $\mathcal{O}_{\mathrm{cond}}$ be the condensed structure sheaf; that is, the condensed enhancement of the structure sheaf $\mathcal{O}$ on $\mathfrak{X}$. Then $R\Gamma(\mathfrak{X}, \mathcal{O}_{\mathrm{cond}})$ is a ring object in $D(\mathrm{Cond}(\mathrm{Ab}))$. We will need the condensed version of a basic acyclicity theorem for affine formal schemes.

**Lemma 3.5.1.** *Let $A$ be an admissible topological ring, so that $\mathfrak{X} = \operatorname{Spf} A$ is an affine formal scheme. The condensed structure sheaf is acyclic on $\mathfrak{X}$, in the sense that the natural map $\underline{A} \to R\Gamma(\mathfrak{X}, \mathcal{O}_{\mathrm{cond}})$ is an isomorphism.*

*Proof.* Considering that $\mathcal{O}_{\operatorname{Spf} A} = \lim_I \mathcal{O}_{\operatorname{Spec} A/I}$ is a limit of sheaves of discrete rings, we have that $\mathcal{O}_{\mathrm{cond}} = \lim_I \mathcal{O}_{\operatorname{Spec} A/I}$ is a limit of sheaves of condensed abelian groups which are discrete. (No $R\lim$ terms appear: this is an instance of the fact that light condensed sets form a replete topos.) On derived global sections, this implies that

$$R\Gamma(\mathfrak{X}, \mathcal{O}_{\mathrm{cond}}) \cong \lim_I R\Gamma(\operatorname{Spec} A/I, \mathcal{O}_{\operatorname{Spec} A/I}).$$

We claim that $R\Gamma(\operatorname{Spec} A/I, \mathcal{O}_{\operatorname{Spec} A/I}) \cong \underline{A/I}$, which implies the lemma, since $\lim_I \underline{A/I} \cong \underline{A}$.

In other words, we have reduced the lemma to the case that $A$ is discrete. By a standard argument [Aut, Tag 01EW], it suffices to show that the Čech cohomology of $\mathcal{O}_{\mathrm{cond}}$ is acyclic with respect to the standard open cover of $\operatorname{Spec} A$ by opens $\{f_i \neq 0\}$ for elements $f_1, \dots, f_n \in A$ generating the unit ideal. Recall [Aut, Tag 01X9] the uncondensed version of this statement, which is the statement that the extended Čech complex

$$0 \to A \to \prod_i A[1/f_i] \to \prod_{i,j} A[1/f_i f_j] \to \cdots$$

is exact. To get the condensed version, we need the fact that the functor $X \mapsto \underline{X}$ from discrete abelian groups to condensed abelian groups preserves exact complexes. This follows from the formula $\underline{X}(S) = \varinjlim_i X^{S_i}$ whenever $X$ is discrete and $S = \varprojlim S_i$ is profinite. □

3.6. **Pro-étale cohomology of rigid-analytic spaces.** Let $X$ be a rigid-analytic space over $K$. We swiftly recall some material from [Sch13a, §3,4] regarding the pro-étale site $X_{\text{proét}}$.

Objects in $X_{\text{proét}}$ are formal limits $U = \varprojlim U_i$, where $i$ runs over a cofiltered index set, the $U_i$ are rigid-analytic spaces étale over $X$, and each transition map $U_i \to U_j$ commutes with the maps to $X$. It is required that $U_i \to U_j$ is finite étale and surjective for large $i > j$. For an object $U = \varprojlim U_i$, let $|U| = \varprojlim |U_i|$ be its underlying topological space. A covering in $X_{\text{proét}}$ is a family of pro-étale morphisms $f_j \colon V_j \to U$ such that the underlying topological space $|U|$ is covered by the $f_j(|V_j|)$.

There is a natural morphism of sites $\nu \colon X_{\text{proét}} \to X_{\text{ét}}$ from the pro-étale site to the étale site. The (uncompleted) integral structure sheaf $\mathcal{O}^+$ on $X_{\text{proét}}$ is $\mathcal{O}^+ = \nu^* \mathcal{O}_{X_{\text{ét}}}$, which is to say that if $U = \varprojlim U_i$ is an object in $X_{\text{proét}}$, then

$$\mathcal{O}^+(U) = \varinjlim \mathcal{O}^+(U_i).$$

The (completed) integral structure sheaf $\hat{\mathcal{O}}^+$ on $X_{\text{proét}}$ is defined as the $p$-adic completion:

$$\hat{\mathcal{O}}^+ = \varprojlim \mathcal{O}^+/p^n,$$

and the structure sheaf is $\hat{\mathcal{O}} := \hat{\mathcal{O}}^+[1/p]$.

The pro-étale cohomology $H^i(X_{\text{proét}}, \hat{\mathcal{O}}^+)$ of a rigid-analytic space $X$ will be of special interest to us. To investigate it, we will make crucial use of perfectoid spaces. Let us recall the relevant definitions from [SW20], adapted to the case where all structures live in characteristic 0. A topological $\mathbb{Q}_p$-algebra $R$ is *perfectoid* if the following conditions hold:

(1) $R$ is uniform, meaning its subring $R^\circ$ of power-bounded elements is bounded,
(2) $R^\circ$ is $p$-adically complete,

(3) There exists an element $\varpi \in R^\circ$ such that $\varpi^p | p$ holds in $R^\circ$, and such that the $p$th power Frobenius map

$$R^\circ/\varpi \to R^\circ/\varpi^p$$

is an isomorphism.

In particular there exists for all $n \geq 1$ an element $\varpi_n$ whose $p^n$th power is the product of $\varpi$ by a unit in $R^\circ$. A Huber pair $(R, R^+)$ is a *perfectoid affinoid algebra* over $\operatorname{Spa}(\mathbb{Q}_p, \mathbb{Z}_p)$ if $R$ is perfectoid. An adic space over $\operatorname{Spa}(\mathbb{Q}_p, \mathbb{Z}_p)$ is perfectoid if it admits a cover by adic spectra of perfectoid affinoid algebras.

Suppose that $X$ is a rigid-analytic space over a field containing $\mathbb{Q}_p$. We say that an object $U \in X_{\text{proét}}$ is *affinoid perfectoid* if $U = \varprojlim U_i$, where:

(1) $U_i = \operatorname{Spa}(R_i, R_i^+)$ is affinoid; and
(2) $R = R^+[1/p]$ is a perfectoid $K$-algebra, where $R^+$ is the $p$-adic completion of $\varinjlim R_i^+$.

More generally, an object $U \in X_{\text{proét}}$ is perfectoid if it admits an open cover by perfectoid affinoid subobjects.

Perfectoid objects are useful to the calculation of $H^i(X_{\text{proét}}, \hat{\mathcal{O}}^+)$ for the following two reasons:

(1) The affinoid perfectoid objects $U \in X_{\text{proét}}$ form a basis for the topology [Sch13a, Proposition 4.8].
(2) Suppose $U \in X_{\text{proét}}$ is affinoid perfectoid. Then for all $i \geq 1$, $H^i(U_{\text{proét}}, \hat{\mathcal{O}}^+)$ is almost zero, in the sense that it is annihilated by $\varpi_n$ for all $n$ [Sch13a, Lemma 4.10]. In particular $H^i(U_{\text{proét}}, \hat{\mathcal{O}}) = 0$ for $i \geq 1$.

The facts above suggest a strategy for computing $H^i(X_{\text{proét}}, \hat{\mathcal{O}}^+)$: By (1), there exists a pro-étale cover $f_i \colon U_i \to X$ where each $U_i$ is perfectoid affinoid, and by (2), $R\Gamma(X_{\text{proét}}, \hat{\mathcal{O}}^+)$ is "almost" computed by the Čech complex $\mathcal{O}^+(X) \to \prod_i \mathcal{O}^+(U_i) \to \prod_{i,j} \mathcal{O}^+(U_i \times_X U_j) \to \cdots$. The precise consequences for the pro-étale cohomology of rigid-analytic spaces will be reviewed in Section 5.

It will be important to consider the pro-étale cohomology $R\Gamma(X_{\text{proét}}, \hat{\mathcal{O}}^+)$ of a rigid-analytic space as a complex of condensed abelian groups. Once again, we proceed by defining a condensed enhancement $\hat{\mathcal{O}}^+_{\text{cond}}$ of $\hat{\mathcal{O}}^+$. The integral structure sheaf $\mathcal{O}^+$ on an adic space is already a sheaf of topological rings; therefore so are $\mathcal{O}^+_{X_{\text{ét}}}$, the uncompleted integral structure sheaf $\mathcal{O}^+ = \nu^* \mathcal{O}^+_{X_{\text{ét}}}$, and the completed structure sheaf $\hat{\mathcal{O}}^+ = \varprojlim \mathcal{O}^+/p^n$. Let $\hat{\mathcal{O}}^+_{\text{cond}}$ be the condensed enhancement of $\hat{\mathcal{O}}^+$. Then $R\Gamma(X_{\text{proét}}, \hat{\mathcal{O}}^+_{\text{cond}})$ is a condensed enhancement of $R\Gamma(X_{\text{proét}}, \hat{\mathcal{O}}^+)$.

There is another natural way to define a condensed enhancement of $R\Gamma(X_{\text{proét}}, \hat{\mathcal{O}}^+)$. Observe that we have a morphism of sites:

$$\lambda \colon X_{\text{proét}} \to *_{\text{proét}}$$

Indeed, if $S = \varprojlim S_i$ is a profinite set, then $X \times S = \varprojlim X \times S_i$ is an object of $X_{\text{proét}}$. Consequently, if $\mathcal{F}$ is a sheaf of abelian groups on $X_{\text{proét}}$, we may define $\Gamma_{\text{cond}}(X_{\text{proét}}, \mathcal{F}) = \lambda_* \mathcal{F}$, an object of $\operatorname{Cond}(\operatorname{Ab})$, and $R\Gamma_{\text{cond}}(X_{\text{proét}}, \mathcal{F}) = R\lambda_* \mathcal{F}$, an object of $D(\operatorname{Cond}(\operatorname{Ab}))$. Applied to $\mathcal{F} = \hat{\mathcal{O}}^+$, we find a complex $R\Gamma_{\text{cond}}(X_{\text{proét}}, \hat{\mathcal{O}}^+)$. In fact we have not found anything new:

**Lemma 3.6.1.** *We have an isomorphism in $D(\operatorname{Cond}(\operatorname{Ab}))$:*

$$R\Gamma_{\text{cond}}(X_{\text{proét}}, \hat{\mathcal{O}}^+) \cong R\Gamma(X_{\text{proét}}, \hat{\mathcal{O}}^+_{\text{cond}})$$

*Proof.* Since $X_{\text{proét}}$ admits a basis consisting of perfectoid affinoids, this reduces to the claim that for all perfectoid affinoids $U = \operatorname{Spa}(R, R^+) \in X_{\text{proét}}$, we have a natural isomorphism in $\operatorname{Cond}(\operatorname{Ab})$:

$$\Gamma_{\text{cond}}(U, \hat{\mathcal{O}}^+) \cong \hat{\mathcal{O}}^+_{\text{cond}}(U)$$

This in turn reduces to the following calculation, for any profinite set $S = \varprojlim S_i$:

$$\begin{aligned} H^0(U \times S, \hat{\mathcal{O}}^+) &\cong \left(\varinjlim H^0(U \times S_i, \hat{\mathcal{O}}^+)\right)^{\wedge}_{(p)} \\ &\cong \left(\varinjlim C_{\mathrm{cts}}(S_i, R^+)\right)^{\wedge}_{(p)} \\ &\cong C_{\mathrm{cts}}(S, R^+) \\ &\cong \hat{\mathcal{O}}^+_{\mathrm{cond}}(U)(S) \end{aligned}$$

In the penultimate step we used the fact that $R^+$ is $p$-adically complete. This follows from the fact that $R^\circ$ is $p$-adically complete (by definition), together with the fact that $R^\circ/R^+$ is $p$-torsion. For the claim that $R^\circ/R^+$ is $p$-torsion: For $f \in R^\circ$ the sequence $(pf)^N$ converges to 0 as $N \to \infty$. Since $R^+$ is open, we have $(pf)^N \in R^+$ for sufficiently large $N$. Since $R^+$ is integrally closed in $R$ we have $pf \in R^+$. (In fact $R^\circ/R^+$ is almost zero.) □

In the context of the proof of Lemma 3.6.1 we have

$$C_{\mathrm{cts}}(S, R^+) \cong \underline{\mathrm{Hom}}(\mathbb{Z}[S], \underline{R^+}).$$

Applying this over a covering of $X$ by perfectoid affinoids $U = \mathrm{Spa}(R, R^+)$ in $X_{\text{proét}}$, we obtain an isomorphism in $D(\mathrm{Solid})$:

$$R\Gamma((X \times S)_{\text{proét}}, \hat{\mathcal{O}}^+_{\mathrm{cond}}) \cong \underline{\mathrm{RHom}}(\mathbb{Z}[S]^{\blacksquare}, R\Gamma(X, \hat{\mathcal{O}}^+_{\mathrm{cond}})). \tag{3.6.2}$$

Let $G$ be a profinite group. A *perfectoid pro-étale $G$-torsor* over $X$ is a surjective pro-étale morphism $Y \to X$ from a perfectoid space $Y$ admitting an action $G \times Y \to Y$ lying over the trivial action of $X$, such that the map $G \times Y \to Y \times_X Y$ given by $(g, y) \mapsto (y, gy)$ is an isomorphism.

**Proposition 3.6.3.** *Let $X$ be a rigid-analytic space, and let $Y \to X$ be a perfectoid pro-étale $G$-torsor. There is an isomorphism in $D(\mathrm{Solid})$:*

$$R\Gamma(X_{\text{proét}}, \hat{\mathcal{O}}^+_{\mathrm{cond}}) \cong R\Gamma(Y_{\text{proét}}, \hat{\mathcal{O}}^+_{\mathrm{cond}})^{hG}.$$

*Proof.* Since $Y \to X$ is pro-étale, the pro-étale cohomology of $X$ can be computed by means of the simplicial cover:

$$\cdots \substack{\to\\\to\\\to} Y \times_X Y \rightrightarrows Y$$

Namely, $R\Gamma(X_{\text{proét}}, \hat{\mathcal{O}}^+_{\mathrm{cond}})$ is quasi-isomorphic to the totalization of the corresponding simplicial object

$$R\Gamma(Y, \hat{\mathcal{O}}^+_{\mathrm{cond}}) \to R\Gamma(Y \times_X Y, \hat{\mathcal{O}}^+_{\mathrm{cond}}) \rightrightarrows \cdots.$$

Since $Y \to X$ is a $G$-torsor, the $n$th term in the complex is quasi-isomorphic to

$$R\Gamma((G^n \times Y)_{\text{proét}}, \hat{\mathcal{O}}^+_{\mathrm{cond}}) \cong R\underline{\mathrm{Hom}}(\mathbb{Z}[G^n]^{\blacksquare}, R\Gamma(Y_{\text{proét}}, \hat{\mathcal{O}}^+_{\mathrm{cond}}))$$

by (3.6.2). Lemma 3.3.4 identifies the latter as the $n$th term in a cosimplicial resolution whose totalization computes $R\Gamma(Y_{\text{proét}}, \hat{\mathcal{O}}^+_{\mathrm{cond}})^{hG}$. □

**Example 3.6.4.** Let $K$ be a nonarchimedean field of characteristic $(0, p)$. Let $\overline{K}$ be an algebraic closure, and let $C$ be the metric completion of $\overline{K}$. Then $\mathrm{Spa}(C, \mathcal{O}_C) \to \mathrm{Spa}(K, \mathcal{O}_K)$ is a pro-étale torsor for the group $\mathrm{Gal}(\overline{K}/K)$. Furthermore, since $C$ is algebraically closed, every pro-étale cover of $\mathrm{Spa}(C, \mathcal{O}_C)$ is split, meaning that $H^i(\mathrm{Spa}(C, \mathcal{O}_C)_{\text{proét}}, \hat{\mathcal{O}}^+) = 0$ for $i > 0$. Therefore by Proposition 3.6.3 we have an isomorphism in $D(\mathrm{Solid})$:

$$R\Gamma(\mathrm{Spa}(K, \mathcal{O}_K)_{\text{proét}}, \hat{\mathcal{O}}^+_{\mathrm{cond}}) \cong \mathcal{O}_C^{h\,\mathrm{Gal}(\overline{K}/K)}.$$

3.7. $R\Gamma(X_{\text{proét}}, \hat{\mathcal{O}}^+_{\text{cond}})$ **for affinoid** $X$**.** The pro-étale cohomology $R\Gamma(X_{\text{proét}}, \hat{\mathcal{O}}^+)$ of a rigid-analytic space $X$ is a complex of abelian groups. Whereas in applications like our Theorem 3.9.4, it will be necessary to control $R\Gamma(X_{\text{proét}}, \hat{\mathcal{O}}^+_{\text{cond}})$ as a complex of condensed abelian groups. The main theorem of this subsection is that, when $X$ is affinoid, $R\Gamma(X_{\text{proét}}, \hat{\mathcal{O}}^+_{\text{cond}})$ is $p$-adic in the sense of Definition 3.4.6, and indeed is the $p$-adic completion of $R\Gamma(X_{\text{proét}}, \hat{\mathcal{O}}^+)$ in the sense of Lemma 3.4.7. This will allow us to prove condensed enhancements of comparison theorems relating to $R\Gamma(X_{\text{proét}}, \hat{\mathcal{O}}^+)$.

Recall the $p$-adic completion functor $A \mapsto (\Gamma^* A)^\wedge_p$ from $D(\text{Ab})$ to $D(\text{Solid})$, whose image consists of the $p$-adic solid complexes.

**Lemma 3.7.1.** *Let $X$ be an affinoid rigid-analytic space over a nonarchimedean field of characteristic $(0,p)$. There is a natural isomorphism in $D(\text{Solid})$:*

$$R\Gamma(X_{\text{proét}}, \hat{\mathcal{O}}^+_{\text{cond}}) \cong (\Gamma^* R\Gamma(X_{\text{proét}}, \hat{\mathcal{O}}^+))^\wedge_p$$

To prove the lemma, we need the notion of a *strictly totally disconnected* perfectoid space $Y$ [Sch22, Definition 1.14]: this means that $Y$ is quasi-compact, and every étale cover of $Y$ has a section.

**Lemma 3.7.2.** *Let $Y$ be a strictly totally disconnected perfectoid space. Then $H^i(Y_{\text{proét}}, \hat{\mathcal{O}}^+) = 0$ for all $i > 0$.*

*Proof.* First we claim that the global sections functor $\mathcal{F} \mapsto \mathcal{F}(Y)$ on sheaves of abelian groups on $Y_{\text{ét}}$ is exact. Indeed, if $f\colon \mathcal{F} \to \mathcal{G}$ is surjective, and $s \in \mathcal{G}(Y)$, there exists an étale cover $c\colon Y' \to Y$ and $t \in \mathcal{F}(Y')$ such that $f(t) = c^*(s)$. Since $Y$ is strictly totally disconnected, there exists a section $\sigma$ to $c$, i.e., $c\sigma = \text{id}_Y$. Then $s = \text{id}_Y^*(s) = \sigma^* c^*(s) = \sigma^* f(t) = f\sigma^*(t)$ lies in the image of $\mathcal{F}(Y)$ under $f$. As a result, $H^i(Y_{\text{ét}}, \mathcal{F}) = 0$ for all $i > 0$.

Let $\nu\colon Y_{\text{proét}} \to Y_{\text{ét}}$ be the projection. For a complex $A$ of sheaves of abelian groups on $Y_{\text{ét}}$, the unit of the adjunction $A \to R\nu_*\nu^* A$ is a quasi-isomorphism. A proof of this fact is found in [BS15, Corollary 5.1.6] (in the context of schemes, but the formalism is the same).

We apply the above observations to the sheaf $\mathcal{F} = \mathcal{O}^+/p^n$ on $Y_{\text{ét}}$, for each $n \geq 1$. The pullback of this sheaf to $Y_{\text{proét}}$, evaluated on a perfectoid affinoid $U = \varprojlim U_i$, is

$$\nu^* \mathcal{O}^+/p^n(U) = \varinjlim (\mathcal{O}^+/p^n)(U_i) = (\hat{\mathcal{O}}^+/p^n)(U),$$

which is to say that $\nu^*\mathcal{O}^+/p^n \cong \hat{\mathcal{O}}^+/p^n$ as sheaves on $Y_{\text{proét}}$. Therefore

$$R\Gamma(Y_{\text{proét}}, \hat{\mathcal{O}}^+/p^n) \cong R\Gamma(Y_{\text{ét}}, R\nu_*\nu^*\mathcal{O}^+/p^n) \cong R\Gamma(Y_{\text{ét}}, \mathcal{O}^+/p^n) \cong (\mathcal{O}^+/p^n)(Y)$$

has vanishing cohomology in degree $> 0$. Taking the limit in $n$ shows that $R\Gamma(Y_{\text{proét}}, \hat{\mathcal{O}}^+) \cong \mathcal{O}^+(Y)$ has no higher cohomology either. ☐

**Lemma 3.7.3.** *Let $X$ be an affinoid rigid-analytic space over a nonarchimedean field of characteristic $(0,p)$. There exists a pro-étale cover of $X$ by an affinoid perfectoid space $\tilde{X}$ which is strictly totally disconnected.*

*Proof.* There exists a pro-étale cover $X' \to X$ which is affinoid perfectoid by [Sch13a, Proposition 4.8], and then there exists a pro-étale cover $\tilde{X} \to X'$ with $\tilde{X}$ strictly totally disconnected by [Sch22, Lemma 7.18]. ☐

**Lemma 3.7.4.** *Let $Y$ be a strictly totally disconnected perfectoid space, and let $Y' \to Y$ be a pro-étale cover, with $Y'$ perfectoid affinoid. Then $Y'$ is also strictly totally disconnected.*

*Proof.* Since $Y'$ is affinoid, $Y' \to Y$ is quasicompact and separated, in which case [Sch22, Lemma 7.19] applies to conclude that $Y'$ is strictly totally disconnected. ☐

We can now give the proof of Lemma 3.7.1. By Lemma 3.7.3, there exists a pro-étale cover $\tilde{X} \to X$ with $\tilde{X}$ strictly totally disconnected. For $i \geq 1$, let $\tilde{X}^{(i)}$ be the $i$-fold fiber product of $\tilde{X}$ over $X$. Any of the projections $\tilde{X}^{(i)} \to \tilde{X}$, being a composite of base changes of pro-étale morphisms, is itself pro-étale. By Lemma 3.7.4, each $\tilde{X}^{(i)}$ is strictly totally disconnected and hence affinoid. By Lemma 3.7.2, $H^j(\tilde{X}^{(i)}_{\text{proét}}, \hat{\mathcal{O}}^+_{\text{cond}}) = 0$ for all $j > 0$. Therefore $R\Gamma(X_{\text{proét}}, \hat{\mathcal{O}}^+_{\text{cond}})$ is computed by the Čech complex associated to the simplicial complex associated to the pro-étale cover $\tilde{X} \to X$. Explicitly, $R\Gamma(X_{\text{proét}}, \hat{\mathcal{O}}^+_{\text{cond}})$ is quasi-isomorphic to the complex of condensed abelian groups with terms $\hat{\mathcal{O}}^+_{\text{cond}}(\tilde{X}^{(i)})$. By our definition of $\hat{\mathcal{O}}^+_{\text{cond}}$:

$$\hat{\mathcal{O}}^+_{\text{cond}}(\tilde{X}^{(i)}) = \underline{\hat{\mathcal{O}}^+(\tilde{X}^{(i)})} \cong \hat{\mathcal{O}}^+(\tilde{X}^{(i)})^\wedge_p,$$

where the isomorphism holds because $\hat{\mathcal{O}}^+(\tilde{X}^{(i)})$ carries the $p$-adic topology (which is in turn true because $\tilde{X}^{(i)}$ is affinoid perfectoid). Since the complex with terms $\hat{\mathcal{O}}^+(\tilde{X}^{(i)})$ computes $R\Gamma(X_{\text{proét}}, \hat{\mathcal{O}}^+)$ (again by Lemma 3.7.2), we find that the $p$-adic completion of $R\Gamma(X_{\text{proét}}, \hat{\mathcal{O}}^+)$ is quasi-isomorphic to $R\Gamma(X_{\text{proét}}, \hat{\mathcal{O}}^+_{\text{cond}})$.

3.8. **Continuous cohomology of $p$-adic Lie groups.** As Theorem B is a statement about the continuous cohomology of the Morava stabilizer group $\mathbb{G}_n \cong \mathcal{O}_D^\times \rtimes \hat{\mathbb{Z}}$, it will be useful to collect some basic results on the continuous cohomology of $p$-adic Lie groups such as $\mathcal{O}_D^\times$, with coefficients in $\mathbb{Q}_p$.

All the relevant techniques stem from Lazard [Laz65], who studied the continuous cohomology $H^i_{\text{cts}}(G, V)$ for a certain class of topological groups $G$ ("complete $p$-valued groups of finite rank") and a certain category of $\mathbb{Q}_p$-linear representations $V$ of $G$. Lazard defines a Lie $\mathbb{Q}_p$-algebra $\operatorname{Lie} G$ which acts on $V$, and proves a comparison theorem [Laz65, Théorème V.2.4.9] which identifies $H^i_{\text{cts}}(G, V)$ with the Lie algebra cohomology $H^i_{\text{cts}}(\operatorname{Lie} G, V)$.

As an example one may take $G$ to be a sufficiently small subgroup of the group of $\mathbb{Q}_p$-points of a reductive group $\mathbf{G}/\mathbb{Q}_p$, and $V$ to arise from an algebraic representation of $\mathbf{G}$ on a finite-dimensional $\mathbb{Q}_p$-vector space. In this case $\operatorname{Lie} G \cong \operatorname{Lie} \mathbf{G}$, and therefore:

$$H^i_{\text{cts}}(G, V) \cong H^i(\operatorname{Lie} \mathbf{G}, V).$$

See [CW74] for a discussion of this case of Lazard's theorem.

Morava [Mor85, Remark 2.2.5] noted the structure of $H^*_{\text{cts}}(\mathcal{O}_D^\times, \mathbb{Q}_p)$ and its application to homotopy theory. We summarize that structure here, together with its complement for the group $\operatorname{GL}_n(\mathbb{Z}_p)$.

**Proposition 3.8.1.** *Let $G$ be either of the groups $\operatorname{GL}_n(\mathbb{Z}_p)$ or $\mathcal{O}_D^\times$. Consider the trivial action of $G$ on $\mathbb{Q}_p$. There is an isomorphism of graded $\mathbb{Q}_p$-algebras:*

$$H^*_{\text{cts}}(G, \mathbb{Q}_p) \cong \Lambda_{\mathbb{Q}_p}(x_1, x_3, \ldots, x_{2n-1}).$$

*In the case of $G = \mathcal{O}_D^\times$, the conjugation action of $D^\times$ on $H^*_{\text{cts}}(\mathcal{O}_D^\times, \mathbb{Q}_p)$ is trivial.*

Before discussing the proof of Proposition 3.8.1, we recall the calculation of the relevant Lie algebra cohomology, following [CE48], [Kos50].

**Lemma 3.8.2.** *Let $k$ be a field of characteristic 0. Suppose $G/k$ is a form of $\operatorname{GL}_n$, and let $\mathfrak{g} = \operatorname{Lie} G$. There is an isomorphism of graded $k$-algebras:*

$$H^*(\mathfrak{g}, k) \cong \Lambda_k(x_1, x_3, \ldots, x_{2n-1})$$

*The adjoint action of $G$ on $\mathfrak{g}$ induces the trivial action of $G$ on $H^*(\mathfrak{g}, k)$.*

*Proof.* When $K$ is a compact Lie group, $H^*(\mathrm{Lie}\, K, \mathbb{R})$ is isomorphic to the de Rham cohomology ring $H^*_{\mathrm{dR}}(K)$. Generally, this is a graded-commutative $\mathbb{R}$-algebra whose primitive elements live in odd degree. In the special case $K = U(n)$, the sequence of fibrations $U(n-1) \to U(n) \to S^{2n-1}$ allows one to identify the cohomology $H^*_{\mathrm{dR}}(U(n))$ with the $\mathbb{R}$-cohomology of a product of spheres $S^1 \times S^3 \times \cdots \times S^{2n-1}$. The complexified Lie algebra of $U(n)$ is $\mathfrak{gl}_n(\mathbb{C})$, and so:

$$H^*(\mathfrak{gl}_n(\mathbb{C}), \mathbb{C}) \cong \Lambda_{\mathbb{C}}(x_1, x_3, \dots, x_{2n-1}) \tag{3.8.3}$$

(See [Kos50] for an explicit description, due to Dynkin, of the elements $x_1, x_3, \dots, x_{2n-1}$ in terms of cocycles.). The Lie algebra cohomology $H^*(\mathfrak{gl}_n(\mathbb{C}), \mathbb{C})$ is a representation of the algebraic group $\mathrm{GL}_n /\mathbb{C}$ via its adjoint action, which we claim is trivial. The action factors through the semisimple group $\mathrm{PGL}_n /\mathbb{C}$. An inductive argument on the degree $i \geq 0$ reduces the claim to showing that $\mathrm{PGL}_n /\mathbb{C}$ acts trivially on the quotient of $H^i(\mathfrak{gl}_n(\mathbb{C}), \mathbb{C})$ by its imprimitive subspace. This quotient is at most 1-dimensional, and there are no nontrivial 1-dimensional representations of $\mathrm{PGL}_n /\mathbb{C}$. This concludes the case of $k = \mathbb{C}$.

For general $k$, reduce to the case that $k$ is finitely generated, and embed $k$ into $\mathbb{C}$. Since $\mathfrak{g} \otimes \mathbb{C} \cong \mathfrak{gl}_n(\mathbb{C})$, we find that $H^*(\mathfrak{g}, k) \otimes \mathbb{C}$ is an exterior algebra. This implies the result for $H^*(\mathfrak{g}, k)$, as one sees by comparing dimensions in each degree. □

Returning to the situation of Proposition 3.8.1, let $G$ be either of the groups $\mathrm{GL}_n(\mathbb{Z}_p)$ or $\mathcal{O}_D^\times$. Then $G$ is a subgroup of the group of $\mathbb{Q}_p$-points of an algebraic group $\mathbf{G}/\mathbb{Q}_p$ which is a form of $\mathrm{GL}_n$.

Let $U \subset G$ be a finite-index normal sugroup which is sufficiently small so as to satisfy the hypothesis of Lazard's theorem, so that $H^*_{\mathrm{cts}}(U, \mathbb{Q}_p) \cong H^*(\mathrm{Lie}\, \mathbf{G}, \mathbb{Q}_p)$, and then $H^*(\mathrm{Lie}\, \mathbf{G}, \mathbb{Q}_p)$ is calculated in Lemma 3.8.2: it is an exterior algebra over $\mathbb{Q}_p$. The group $G$ acts on $H^*_{\mathrm{cts}}(U, \mathbb{Q}_p)$ by outer automorphisms; by naturality of Lazard's isomorphism, this corresponds to the action of $G$ on $H^*(\mathrm{Lie}\, \mathbf{G}, \mathbb{Q}_p)$ via the adjoint action; this action is trivial by Lemma 3.8.2.

To extend the result to $H^*_{\mathrm{cts}}(G, \mathbb{Q}_p)$, use the Lyndon–Hochschild–Serre spectral sequence Proposition 3.3.6:

$$H^i(G/U, H^j_{\mathrm{cts}}(U, \mathbb{Q}_p)) \implies H^{i+j}_{\mathrm{cts}}(G, \mathbb{Q}_p),$$

noting that all terms with $i > 0$ vanish (because $G/U$ is finite), and that $G/U$ acts trivially on $H^j_{\mathrm{cts}}(U, \mathbb{Q}_p)$ as noted above. In the case $G = \mathcal{O}_D^\times$, the triviality of the adjoint action of $G$ on $H^j_{\mathrm{cts}}(U, \mathbb{Q}_p)$ extends to $D^\times$.

*Remark* 3.8.4. In the context of Proposition 3.8.1, it is possible to obtain sharper results on the integral cohomology $H^*_{\mathrm{cts}}(G, \mathbb{Z}_p)$ using the integral version of Lazard's theorem developed by Huber–Kings–Naumann [HKN11].

We record one more lemma concerning the extension of these results from $\mathcal{O}_D^\times$ to the full Morava stabilizer group $\mathbb{G}_n$.

**Lemma 3.8.5.** *Let* $W = W(\overline{\mathbb{F}}_p)$ *and* $K = W[1/p]$.

1. *The continuous cohomology* $H^i_{\mathrm{cts}}(\mathrm{Gal}(\overline{\mathbb{F}}_p/\mathbb{F}_p), W)$ *is* $\mathbb{Z}_p$ *if* $i = 0$, *and is* 0 *otherwise.*
2. *Let* $\mathbb{G}_n$ *act on* $K$ *through its quotient* $\mathrm{Gal}(\overline{\mathbb{F}}_p/\mathbb{F}_p)$. *There is an isomorphism of graded* $\mathbb{Q}_p$*-algebras:*

$$H^*_{\mathrm{cts}}(\mathbb{G}_n, K) \cong \Lambda_{\mathbb{Q}_p}(x_1, x_3, \dots, x_{2n-1})$$

*Proof.* (1) Since $\mathrm{Gal}(\overline{\mathbb{F}}_p/\mathbb{F}_p) \cong \hat{\mathbb{Z}}$, it is enough to show vanishing of cohomology in degree 1. Since $W$ is $p$-adically complete, this is further reduced to showing that $H^1_{\mathrm{cts}}(\mathrm{Gal}(\overline{\mathbb{F}}_p/\mathbb{F}_p), \overline{\mathbb{F}}_p) = 0$; this is true because $x \mapsto x^p - x$ is surjective on $\overline{\mathbb{F}}_p$.

(2) Consider the Lyndon–Hochschild–Serre spectral sequence Proposition 3.3.6:

$$H^i_{\mathrm{cts}}(\mathrm{Gal}(\overline{\mathbb{F}}_p/\mathbb{F}_p), H^j_{\mathrm{cts}}(\mathcal{O}_D^\times, K)) \implies H^{i+j}_{\mathrm{cts}}(\mathbb{G}_n, K).$$

Consider the action of $\mathrm{Gal}(\overline{\mathbb{F}}_p/\mathbb{F}_p)$ on $H^j_{\mathrm{cts}}(\mathcal{O}_D^\times, K) = H^j_{\mathrm{cts}}(\mathcal{O}_D^\times, \mathbb{Q}_p) \otimes_{\mathbb{Z}_p} W$. The action on the first factor is trivial by Proposition 3.8.1, and the action on the second factor has no higher cohomology by (1). Therefore:

$$H^*_{\mathrm{cts}}(\mathbb{G}_n, K) \cong H^*_{\mathrm{cts}}(\mathcal{O}_D^\times, \mathbb{Q}_p),$$

at which point we apply Proposition 3.8.1. □

3.9. **Overview of the proof of Theorem B.** Let $p$ be a prime number, and let $\Gamma$ be a formal group of dimension 1 and height $n$ over $\overline{\mathbb{F}}_p$. Let LT be the functor of deformations of $\Gamma$. Then LT is representable by a formal scheme $\mathrm{Spf}\, A$ whose coordinate ring is isomorphic to a power series ring over $W = W(\overline{\mathbb{F}}_p)$:

$$A \cong W[\![u_1, \dots, u_{n-1}]\!].$$

Let $\mathcal{O}_D = \mathrm{End}\, \Gamma$; then $\mathcal{O}_D$ is the ring of integers in a division algebra $D$ over $\mathbb{Q}_p$ of invariant $1/n$. Explicitly, $\mathcal{O}_D$ is generated over $W(\mathbb{F}_{p^n})$ by an element $\Pi$ satisfying $\Pi^n = p$ and $\Pi\alpha = \sigma(\alpha)\Pi$, where $\alpha \in W(\mathbb{F}_{p^n})$ and $\sigma \in \mathrm{Aut}\, W(\mathbb{F}_{p^n})$ is the Frobenius automorphism. The Morava stabilizer group $\mathbb{G}_n$ from Section 2 is the profinite completion of $D^\times$; it fits into an exact sequence:

$$1 \to \mathcal{O}_D^\times \to \mathbb{G}_n \to \mathrm{Gal}(\overline{\mathbb{F}}_p/\mathbb{F}_p) \to 1.$$

There is a continuous action of $\mathbb{G}_n$ on $A$.

Our main theorem (Theorem B) concerns the continuous cohomology ring $H^*_{\mathrm{cts}}(\mathbb{G}_n, A)$. The action of $\mathbb{G}_n$ on $A$ is rather inexplicit, and it seems difficult to compute $H^*_{\mathrm{cts}}(\mathbb{G}_n, A)$ directly in terms of cocycles. To address this problem, we pass from the formal scheme LT to its rigid-analytic generic fiber $\mathrm{LT}_K$, which is isomorphic to the open unit ball over $K = W[1/p]$. The use of $\mathrm{LT}_K$ has a precedent in chromatic homotopy theory, namely in applications of the *Gross–Hopkins period morphism* [HG94]:

$$\pi \colon \mathrm{LT}_K \to \mathbb{P}(M(\Gamma)).$$

Here $M(\Gamma) \cong K^n$ is the rational Dieudonné module of $\Gamma$, and $\mathbb{P}(M(\Gamma)) \cong \mathbb{P}^{n-1}_K$ is the corresponding projective space, considered as a rigid-analytic space over $K$. The morphism $\pi$ is étale, surjective, equivariant for the action of $\mathcal{O}_D^\times$, and fairly explicit in terms of the variables $u_1, \dots, u_{n-1}$. Using the period morphism $\pi$, it is possible in principle to give formulas for the action of $\mathcal{O}_D^\times$ on $A$, see [DH95].

For our work we will not use the period morphism but rather a different structure possessed by $\mathrm{LT}_K$: *the isomorphism between the two towers*. This was discovered by Faltings [Fal02b], see also [FGL08] for more details. [SW13, Theorem D] poses the isomorphism as relating two perfectoid spaces.

We use the form of Faltings' isomorphism presented in [BSSW24, Theorem 3.7.3], which treats this phenomenon as an isomorphism between $v$-stacks. A $v$-stack is a sheaf of groupoids on the $v$-site of perfectoid spaces [Sch22, Definition 8.1(iii)]. For a rigid-analytic space $X$ over a $p$-adic field, let $X^\diamond$ be the sheaf on the $v$-site represented by $X$. Then $X^\diamond$ is a diamond, meaning the quotient of a perfectoid space by a pro-étale equivalence relation [Sch22, Theorem 1.5].

Let $\mathcal{H}$ be *Drinfeld's symmetric space* (also called Drinfeld's upper half-space) in dimension $n-1$. This is a rigid-analytic space over $\mathbb{Q}_p$, defined as

$$\mathcal{H} = \mathbb{P}^{n-1}_{\mathbb{Q}_p} \backslash \bigcup_H H$$

where $\mathbb{P}^{n-1}_{\mathbb{Q}_p}$ is rigid-analytic projective space, and $H$ runs over all $\mathbb{Q}_p$-rational hyperplanes in $\mathbb{P}^{n-1}_{\mathbb{Q}_p}$. Drinfeld shows that $\mathcal{H}$ is open in $\mathbb{P}^{n-1}_{\mathbb{Q}_p}$, and so is actually a rigid-analytic space (see [Dri74, Proposition 6.1], note that Drinfeld calls this object $\Omega^n$). Then $\mathcal{H}$ admits an action of the group $\mathrm{GL}_n(\mathbb{Q}_p)$ and in particular of its subgroup $\mathrm{GL}_n(\mathbb{Z}_p)$.

**Theorem 3.9.1.** *There is an isomorphism of v-stacks:*

$$[\mathrm{LT}_K^\diamond/\mathbb{G}_n] \cong [\mathcal{H}^\diamond/\mathrm{GL}_n(\mathbb{Z}_p)]$$

*Remark* 3.9.2. The stack $[\mathrm{LT}_K^\diamond/\mathbb{G}_n]$ has a clean interpretation, as the *v-stack of formal groups of height n*. For an affinoid perfectoid space $S = \mathrm{Spa}(R, R^+)$ over $\mathbb{Q}_p$, the $S$-points of $[\mathrm{LT}_K^\diamond/\mathbb{G}_n]$ classify, $v$- locally on $S$, 1-dimensional formal groups $H$ over $R^\circ$, such that for some pseudo-uniformizer $\varpi \in R^\circ$ dividing $p$, the special fiber $H \otimes R^\circ/\varpi$ has height $n$.

We refer the reader to [BSSW24, Theorem 3.7.3] for a detailed proof. For the purposes of being self-contained, we offer a discussion of the isomorphism in Theorem 3.9.1 on the level of $\mathrm{Spa}(C, \mathcal{O}_C)$-points, where $C$ is an algebraically closed nonarchimedean field containing $\mathbb{Q}_p$.

A $\mathrm{Spa}(C, \mathcal{O}_C)$-point of $[\mathrm{LT}_K^\diamond/\mathbb{G}_n]$ is a 1-dimensional formal group $H/\mathcal{O}_C$ whose special fiber has height $n$. Then the torsion $H[p^\infty]$ in $H$ is a connected $p$-divisible group over $\mathcal{O}_C$ of height $n$ and dimension 1, which in turn determines $H$. Now we apply the classification theorem [SW13, Theorem B], which places into equivalence $p$-divisible groups $G$ over $\mathcal{O}_C$ and pairs $(T, W)$, where $T$ is a finite free $\mathbb{Z}_p$-module and $W \subset T \otimes C(-1)$ is a vector subspace. The rank of $T$ is the height of $G$, and the dimension of $W$ is the dimension of $G$. The classification theorem shows that $\mathrm{Spa}(C, \mathcal{O}_C)$-points of $[\mathbb{P}^{n-1,\diamond}_{\mathbb{Q}_p}/\mathrm{GL}_n(\mathbb{Z}_p)]$ classify $p$-divisible groups $G$ over $\mathcal{O}_C$ of height $n$ and dimension 1. We claim the substack $[\mathcal{H}^\diamond_{\mathbb{Q}_p}/\mathrm{GL}_n(\mathbb{Z}_p)]$ classifies those $p$-divisible groups $G$ which are connected. Indeed, again using the classification theorem, there exists a nonzero $\mathbb{Q}_p$-linear form on $T \otimes C(-1)$ vanishing on $W$ if and only if there exists a nonzero map from $G$ towards a 0-dimensional (i.e., étale) $p$-divisible group, which means exactly that $G$ is not connected.

Theorem 3.9.1 suggests a strategy for accessing the cohomology ring $H^*_{\mathrm{cts}}(\mathbb{G}_n, A)$. Theorem 3.9.1 is equivalent to the existence of a diagram of adic spaces in which $\mathfrak{X}$ admits commuting actions of $\mathbb{G}_n$ and $\mathrm{GL}_n(\mathbb{Z}_p)$, and each arrow is a pro-étale torsor for the indicated group. This diagram witnesses an isomorphism in $D(\mathrm{Solid})$, using Proposition 3.6.3:

$$R\Gamma(\mathrm{LT}_{K,\mathrm{pro\acute{e}t}}, \hat{\mathcal{O}}^+_{\mathrm{cond}})^{h\mathbb{G}_n} \xrightarrow{\sim} R\Gamma(\mathcal{H}_{\mathrm{pro\acute{e}t}}, \hat{\mathcal{O}}^+_{\mathrm{cond}})^{h\,\mathrm{GL}_n(\mathbb{Z}_p)}, \tag{3.9.3}$$

since both objects are isomorphic to $R\Gamma(\mathfrak{X}_{\mathrm{pro\acute{e}t}}, \hat{\mathcal{O}}^+_{\mathrm{cond}})^{h(\mathbb{G}_n \times \mathrm{GL}_n(\mathbb{Z}_p))}$. The isomorphism in (3.9.3) is helpful because it translates the opaque action of $\mathcal{O}_D^\times$ on $\mathrm{LT}_K$ into the transparent action of $\mathrm{GL}_n(\mathbb{Z}_p)$ on $\mathcal{H}$. To completely leverage (3.9.3), we will have to say something about the $\hat{\mathcal{O}}^+_{\mathrm{cond}}$-cohomology of $X_{\mathrm{pro\acute{e}t}}$, where $X$ is $\mathrm{LT}_K$ or $\mathcal{H}$, respectively. The following comparison statements appear as Theorem 6.1.4 and Theorem 6.2.6.

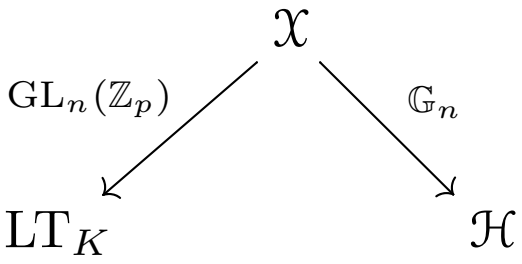

**Theorem 3.9.4.** *The pro-étale cohomology of* $\hat{\mathcal{O}}^+_{\mathrm{cond}}$ *on* $\mathrm{LT}_K$ *and* $\mathcal{H}$ *can be approximated as follows.*

*(1) There is a morphism of algebra objects in* $D(\mathrm{Solid}_{\mathbb{G}_n})$*:*

$$A[\varepsilon] \to R\Gamma(\mathrm{LT}_{K,\mathrm{pro\acute{e}t}}, \hat{\mathcal{O}}^+_{\mathrm{cond}})$$

*The action of* $\mathbb{G}_n$ *on* $\varepsilon$ *is trivial.*

*(2) There is a morphism of algebra objects in* $D(\mathrm{Solid}_{\mathrm{GL}_n(\mathbb{Z}_p)})$*:*

$$\mathbb{Z}_p[\varepsilon] \to R\Gamma(\mathcal{H}_{\mathrm{pro\acute{e}t}}, \hat{\mathcal{O}}^+_{\mathrm{cond}})$$

*The action of* $\mathrm{GL}_n(\mathbb{Z}_p)$ *on* $\mathbb{Z}_p[\varepsilon]$ *is trivial.*

*Here $\varepsilon$ is a square-zero element in cohomological degree 1. Let $\mathcal{C}$ be the cofiber of either of the above morphisms in $D(\mathrm{Solid})$. Then $H^i(\mathcal{C})=0$ for $i\leq 0$, and there exists a single power of $p$ which annihilates $H^i(\mathcal{C})$ for every $i\geq 1$.*

The proof of Theorem 3.9.4 requires the full force of the integral $p$-adic Hodge theory theorems of [BMS18] and [ČK19]. We apply these theorems in Section 5 to control $R\Gamma(X_{\text{proét}},\hat{\mathcal{O}}^+_{\text{cond}})$, where $X$ is a rigid-analytic space $X$ over a local field $K$ of characteristic $(0,p)$.

As a necessary first step, we dispense with the case $X=\mathrm{Spa}(K,\mathcal{O}_K)$. By Example 3.6.4, the pro-étale cohomology of $\mathrm{Spa}(K,\mathcal{O}_K)$ agrees with the Galois cohomology $\mathcal{O}_C^{h\,\mathrm{Gal}(\overline{K}/K)}$, where $C$ is the completion of an algebraic closure of $K$. This cohomology was controlled by Tate [Tat67]. Expressed in our language, Tate's result is that there is a morphism of solid $\mathcal{O}_K$-algebras

$$\mathcal{O}_K[\varepsilon]\to\mathcal{O}_C^{h\,\mathrm{Gal}(\overline{K}/K)}$$

which becomes an isomorphism after base change to $K$. This form of Tate's result will be explained in more detail in the next section.

Combining (3.9.3) with Theorem 3.9.4, we obtain a diagram in $D(\mathrm{Solid})$:

$$\begin{aligned}
A^{h\mathbb{G}_n}\otimes_{\mathbb{Z}_p}\mathbb{Z}_p[\varepsilon] &\cong A[\varepsilon]^{h\mathbb{G}_n}\\
&\to R\Gamma(\mathrm{LT}_{K,\text{proét}},\hat{\mathcal{O}}^+_{\text{cond}})^{h\mathbb{G}_n}\\
&\cong R\Gamma(\mathcal{H}_{\text{proét}},\hat{\mathcal{O}}^+_{\text{cond}})^{h\,\mathrm{GL}_n(\mathbb{Z}_p)}\\
&\leftarrow \mathbb{Z}_p[\varepsilon]^{h\,\mathrm{GL}_n(\mathbb{Z}_p)}\\
&\cong \mathbb{Z}_p^{h\,\mathrm{GL}_n(\mathbb{Z}_p)}\otimes_{\mathbb{Z}_p}\mathbb{Z}_p[\varepsilon]
\end{aligned}$$

Here, each of the two arrows not labeled as an isomorphism becomes an isomorphism after inverting $p$.

We will briefly indicate how this is used to prove Theorem B, leaving the details for Section 6. By Proposition 2.5.1, we have a $\mathbb{G}_n$-equivariant splitting $A=W\oplus A^c$. After inverting $p$ in the above diagram, we arrive at an isomorphism in cohomology:

$$(H^*_{\mathrm{cts}}(\mathbb{G}_n,W)\oplus H^*_{\mathrm{cts}}(\mathbb{G}_n,A^c))\otimes_{\mathbb{Z}_p}\mathbb{Q}_p[\varepsilon]\cong H^*_{\mathrm{cts}}(\mathrm{GL}_n(\mathbb{Z}_p),\mathbb{Q}_p)\otimes_{\mathbb{Q}_p}\mathbb{Q}_p[\varepsilon] \tag{3.9.5}$$

By Lemma 3.8.5, $H^*_{\mathrm{cts}}(\mathbb{G}_n,K)$ and $H^*_{\mathrm{cts}}(\mathrm{GL}_n(\mathbb{Z}_p),\mathbb{Q}_p)$ are isomorphic to the same exterior $\mathbb{Q}_p$-algebra. In particular $\dim_{\mathbb{Q}_p}H^i_{\mathrm{cts}}(\mathbb{G}_n,K)=\dim_{\mathbb{Q}_p}H^i_{\mathrm{cts}}(\mathrm{GL}_n(\mathbb{Z}_p),\mathbb{Q}_p)$ for all $i$.

Comparing dimensions of the $\mathbb{Q}_p$-vector spaces in (3.9.5) shows that $H^*_{\mathrm{cts}}(\mathbb{G}_n,A^c)\otimes_{\mathbb{Z}_p}\mathbb{Q}_p=0$, which is to say that $H^*_{\mathrm{cts}}(\mathbb{G}_n,A^c)$ is torsion: this is the assertion of Theorem B.

## 4. The Galois cohomology of $\mathcal{O}_C$

Let us fix some definitions. A *nonarchimedean field* is a field $K$ which is complete with respect to the topology induced from a nontrivial nonarchimedean valuation $|\,|\colon K\to\mathbb{R}_{\geq 0}$. (Some authors do not require $K$ to be complete, but for our purposes it will be useful to always assume this.) Let $\mathcal{O}_K$ be its valuation ring; i.e., the subring of elements with $|\alpha|\leq 1$. Let $\kappa$ be the residue field of $\mathcal{O}_K$. The characteristic of a nonarchimedean field $K$ refers to the pair $(\mathrm{char}\,K,\mathrm{char}\,\kappa)$.

A *local field* is a nonarchimedean field satisfying the additional conditions: (a) the valuation on $K$ is discrete, in other words $\mathcal{O}_K$ is a discrete valuation ring, and (b) the residue field $\kappa$ is perfect. (Some authors require that the residue field of a local field be finite, which is identical to $K$ being locally compact. However, we want to allow for fields such as $W(\overline{\mathbb{F}}_p)[1/p]$.)

Let $L$ be a (possibly infinite) Galois extension of a nonarchimedean field $K$. The valuation on $K$ extends uniquely to a valuation on $L$, and the completion $\hat{L}$ is a nonarchimedean field admitting a continuous action of $\mathrm{Gal}(L/K)$. (If $L/K$ is finite then it is not necessary to complete.) It is an interesting problem to compute the continuous cohomology groups $H^i_{\mathrm{cts}}(\mathrm{Gal}(L/K), \mathcal{O}_{\hat{L}})$, or at least to approximate these as $\mathcal{O}_K$-modules. For simplicity we will assume throughout that $\operatorname{char} K = 0$.

A basic result along these lines is due to Ax, which settles the problem in degree $i = 0$.

**Theorem 4.0.1** ([Ax70])**.** *Let $K$ be a nonarchimedean field with $\operatorname{char} K = 0$. Let $L/K$ be a Galois extension of nonarchimedean fields. Then the subfield of $\hat{L}$ fixed by $\mathrm{Gal}(L/K)$ is exactly $K$. Consequently $H^0(\mathrm{Gal}(L/K), \mathcal{O}_{\hat{L}}) = \mathcal{O}_K$.*

Results on higher cohomology tend to require that $K$ be a local field. A classical result [Frö83, Theorem 3] attributed to Noether states that if $L/K$ is a finite tamely ramified Galois extension of local fields, then $\mathcal{O}_L$ is a free $\mathcal{O}_K[\mathrm{Gal}(L/K)]$-module, and therefore $H^i(\mathrm{Gal}(L/K), \mathcal{O}_L) = 0$ for $i > 0$. For arbitrary finite extensions there is the following result of Sen:

**Theorem 4.0.2** ([Sen69])**.** *Let $K$ be a local field of characteristic $(0, p)$. Let $L/K$ be a finite Galois extension. Then $H^1(\mathrm{Gal}(L/K), \mathcal{O}_L)$ is $\alpha$-torsion for any $\alpha \in \mathcal{O}_K$ with $|\alpha| \leq |p|^{1/(p-1)}$.*

Sen's theorem does not imply that $H^1_{\mathrm{cts}}(\mathrm{Gal}(L/K), \mathcal{O}_{\hat{L}})$ is torsion for $L/K$ infinite. Indeed this $H^1_{\mathrm{cts}}$ is not torsion when $L = \overline{K}$, as the following classical result of Tate shows.

**Theorem 4.0.3** ([Tat67])**.** *Let $K$ be a local field of characteristic $(0, p)$, and let $C$ be the completion of an algebraic closure $\overline{K}/K$. Let $\Gamma_K = \mathrm{Gal}(\overline{K}/K)$. Then:*

*(1) $H^i_{\mathrm{cts}}(\Gamma_K, C) \cong K$ for $i = 0, 1$,*
*(2) $H^i_{\mathrm{cts}}(\Gamma_K, C) = 0$ for $i \geq 2$.*

*For $j \in \mathbb{Z}$, let $C(j)$ be the $j$th Tate twist: this is $C$ as a topological group, but the $\Gamma_K$-action has been twisted by the $j$th power of the cyclotomic character. For $j \neq 0$ and all $i \geq 0$ we have $H^i_{\mathrm{cts}}(\Gamma_K, C(j)) = 0$.*

Our contribution to this story is the following refinement of Tate's theorem, which controls torsion in the cohomology of $\mathcal{O}_C$. We call a local field $K$ of characteristic $(0, p)$ *tame* if $p$ does not divide its absolute ramification index $e_K$ (defined in the following subsection).

**Theorem 4.0.4.** *Let $K$ be a local field of characteristic $(0, p)$. Let $C$ be the completion of an algebraic closure $\overline{K}/K$. Let $\Gamma_K = \mathrm{Gal}(\overline{K}/K)$. Then:*

*(1) $H^0(\Gamma_K, \mathcal{O}_C) = \mathcal{O}_K$.*
*(2) There is an isomorphism $H^1_{\mathrm{cts}}(\Gamma_K, \mathcal{O}_C) \cong \mathcal{O}_K \oplus T$, where $T$ is $p^4$-torsion (resp., $p^5$-torsion) as $p$ is odd or even, respectively. If $K$ is tame, or if it is a cyclotomic extension of a tame field, this bound can be improved to $p^3$ (resp., $p^4$).*
*(3) For $i > 1$, $H^i_{\mathrm{cts}}(\Gamma_K, \mathcal{O}_C)$ is $p$-torsion.*

Theorem 4.0.4 appears later in this section as Theorem 4.3.9.

We also obtain an explicit bound on the continuous cohomology of the nontrivial Tate twists $\mathcal{O}_C(j)$.

**Theorem 4.0.5.** *With notation as in Theorem 4.0.4, let $j \neq 0$ be an integer, and let $\mathcal{O}_C(j)$ be the $j$th Tate twist of $\mathcal{O}_C$ as a $\mathrm{Gal}(\overline{K}/K)$-module. Then:*

*(1) $H^0(\Gamma_K, \mathcal{O}_C(j)) = 0$.*
*(2) $H^1_{\mathrm{cts}}(\Gamma_K, \mathcal{O}_C(j))$ is $p^{M+v(j)}$-torsion. Here $v(j)$ is the $p$-adic valuation of $j$, and $M = M_K$ is a constant which only depends on $K$ and which is insensitive to passage to a tamely ramified extension of $K$. If $p \nmid e_K$ we may take $M = 2$ if $p$ is odd and $M = 5$ if $p = 2$.*
*(3) For $i \geq 2$, $H^i_{\mathrm{cts}}(\Gamma_K, \mathcal{O}_C(j))$ is $p$-torsion if $p$ is odd and $p^3$-torsion if $p = 2$.*

The theorem appears later as Theorem 4.4.3.

4.1. **Ramification theory for extensions of local fields.** All material in this subsection is taken from Serre's book [Ser79], to which we refer the reader for more details.

Suppose $K$ is a local field of characteristic $(0,p)$. We fix here some notation regarding valuations. The ring of integers $\mathcal{O}_K$ is a discrete valuation ring with maximal ideal $\pi_K\mathcal{O}_K$, where $\pi_K$ is a uniformizing parameter for $K$. Any element $\alpha\in K^\times$ can be written as $\alpha=\pi_K^{v_K(\alpha)}u$, where $u$ is a unit in $\mathcal{O}_K$. The map

$$v_K\colon K^\times\to\mathbb{Z}$$

is then a surjective homomorphism not depending on the choice of $\pi$; this is the *relative valuation map*. For any real number $0<a<1$, the absolute value $|\alpha|=a^{v_K(\alpha)}$ defines the topology on $K$. Note that $p\in K^\times$ (since $K$ has characteristic 0) and $v_K(p)>0$ (since the residue field $\mathcal{O}_K/\pi_K\mathcal{O}_K$ has characteristic $p$); relatedly, $K$ contains $\mathbb{Q}_p$ as a subfield. The number $e_K=v_K(p)$ is the *absolute ramification index* of $K$. We define the *absolute valuation map*

$$v\colon K^\times\to\mathbb{Q}$$

by $v(\alpha)=v_K(\alpha)/e_K$, so that $v(p)=1$; thus $v$ extends the valuation $v_{\mathbb{Q}_p}$ on $\mathbb{Q}_p$.

If $L/K$ is a field extension of finite degree, let $\mathcal{O}_L$ be the integral closure of $\mathcal{O}_K$ in $L$. Then $\mathcal{O}_L$ is also a complete discrete valuation ring [Ser79, Chapter 2, §2, Proposition 3], so that $L$ is a local field of characteristic $(0,p)$. The number $e_{L/K}=v_L(\pi_K)$ is the *relative ramification index* of $L/K$; the extension $L/K$ is *totally ramified* if $e_{L/K}=[L:K]$, *ramified* if $e_{L/K}>1$ and *unramified* if $e_{L/K}=1$. Note that for $\alpha\in K^\times$ we have $v_L(\alpha)=e_{L/K}v_K(\alpha)$. Applied to $\alpha=p$, this last relation gives $e_{L/K}=e_L/e_K$. Thus for $\alpha\in K^\times$ we have $v_K(\alpha)/e_K=v_L(\alpha)/e_L$. As a result, the absolute valuation map $v$ defined above extends unambiguously to $L^\times$. Indeed, if $\overline{K}$ is an algebraic closure of $K$, we continue to write $v\colon\overline{K}^\times\to\mathbb{Q}$ for the absolute valuation map.

Besides the ramification index $e_{L/K}$, there are two other quantities which measure ramification in the extension $L/K$: the *discriminant* $\mathfrak{d}_{L/K}$ and the *different* $\mathcal{D}_{L/K}$. To define these, we need the trace map $\mathrm{tr}_{L/K}\colon L\to K$. The discriminant $\mathfrak{d}_{L/K}$ is the ideal of $\mathcal{O}_K$ generated by the determinant of the matrix $(\mathrm{tr}_{L/K}(e_ie_j))$, where $e_i$ runs through a basis for $\mathcal{O}_L$ as an $\mathcal{O}_K$-module. The different $\mathcal{D}_{L/K}$ is the ideal of $\mathcal{O}_L$ whose inverse is the dual of $\mathcal{O}_L$ under the trace pairing $L\times L\to K$ defined by $(x,y)\mapsto\mathrm{tr}_{L/K}(xy)$. The discriminant and different are related by $N_{L/K}(\mathcal{D}_{L/K})=\mathfrak{d}_{L/K}$, where $N_{L/K}$ is the norm map [Ser79, Chapter 3, §3, Proposition 6]. An extension $L/K$ is unramified if and only if $\mathcal{D}_{L/K}=\mathcal{O}_L$ [Ser79, Chapter 3, §5, Theorem 1], which in turn is true if and only if $\mathfrak{d}_{L/K}=\mathcal{O}_K$.

Now suppose $L/K$ is a finite Galois extension. The Galois group $\mathrm{Gal}(L/K)$ admits an exhaustive decreasing filtration $\mathrm{Gal}(L/K)_u$ by real numbers $u\geq-1$, defined by

$$\mathrm{Gal}(L/K)_u=\left\{\sigma\in\mathrm{Gal}(L/K)\;\middle|\;v_L(\sigma(\alpha)-\alpha)\geq u+1\text{ for all }\alpha\in\mathcal{O}_L\right\}.$$

Since $v_L$ takes integer values, we have $\mathrm{Gal}(L/K)_u=\mathrm{Gal}(L/K)_{\lceil u\rceil}$, where $\lceil u\rceil$ is the ceiling function of $u$.

This filtration of $\mathrm{Gal}(L/K)$ is known as the *lower numbering filtration*. The first three subgroups appearing in this filtration are:

$$\begin{aligned}\mathrm{Gal}(L/K)_{-1}&=\mathrm{Gal}(L/K),\\ \mathrm{Gal}(L/K)_0&=\text{the inertia group of }L/K,\\ \mathrm{Gal}(L/K)_1&=\text{the maximal }p\text{-subgroup of }\mathrm{Gal}(L/K)_0,\end{aligned}$$

see [Ser79, Chapter IV, §2, Corollary 1 and Corollary 3]. The cardinality of the inertia group $\mathrm{Gal}(L/K)_0$ equals the relative ramification index $e_{L/K}$. The extension $L/K$ is *tame* if $p\nmid e_{L/K}$; this condition is equivalent to $\mathrm{Gal}(L/K)_1=1$. Otherwise, $L/K$ is *wild*.

The lower numbering filtration of $\mathrm{Gal}(L/K)$ recovers the different $\mathcal{D}_{L/K}$ by the relation [Ser79, Chapter IV, §1, Proposition 4]:

$$v_L(\mathcal{D}_{L/K}) = \sum_{i=0}^{\infty} (\#\,\mathrm{Gal}(L/K)_i - 1) \tag{4.1.1}$$

The lower number filtration is adapted to subgroups of $\mathrm{Gal}(L/K)$, in the sense that if $K'/K$ is a subextension of $L/K$, then $\mathrm{Gal}(L/K')_u = \mathrm{Gal}(L/K)_u \cap \mathrm{Gal}(L/K')$. There is an alternative descending filtration $\mathrm{Gal}(L/K)^u$ by real numbers $u \geq -1$, the *upper numbering filtration.* This filtration is adapted to quotients of $\mathrm{Gal}(L/K)$, in the sense that if $K'/K$ is a Galois subextension of $L/K$, then $\mathrm{Gal}(K'/K)^u$ is the image of $\mathrm{Gal}(L/K)^u$ under the quotient map $\mathrm{Gal}(L/K) \to \mathrm{Gal}(K'/K)$. To define the upper numbering filtration, one needs the *Herbrand function* $\varphi_{L/K} \colon [-1,\infty) \to [-1,\infty)$, a continuous, increasing, and bijective function defined by the rules:

$$\varphi_{L/K}(u) = \begin{cases} u, & -1 \leq u \leq 0, \\ \int_0^u [\mathrm{Gal}(L/K)_0 : \mathrm{Gal}(L/K)_t]^{-1} dt, & u > 0. \end{cases}$$

The inverse function to $\varphi_{L/K}$ is denoted $\psi_{L/K}$. The upper numbering filtration on $\mathrm{Gal}(L/K)$ is now defined by:

$$\mathrm{Gal}(L/K)^u = \mathrm{Gal}(L/K)_{\psi_{L/K}(u)}$$

for all $u \geq -1$. Note that $\mathrm{Gal}(L/K)^0 = \mathrm{Gal}(L/K)_0$ is the inertia group of $L/K$.

In terms of the upper numbering filtration, the formula for the different in (4.1.1) becomes:

$$v_L(\mathcal{D}_{L/K}) = e_{L/K} \int_{-1}^{\infty} \left(1 - \frac{1}{\#\,\mathrm{Gal}(L/K)^u}\right) du \tag{4.1.2}$$

The property claimed above about the upper numbering filtration being adapted to quotients is [Ser79, Chapter IV, §3, Proposition 14]. This property allows for the definition of $\mathrm{Gal}(L/K)^u$ when $L/K$ is an infinite Galois extension. Indeed, the $\mathrm{Gal}(L_i/K)^u$ form a projective system as $L_i$ ranges over Galois subextensions of $L/K$, and one defines $\mathrm{Gal}(L/K)^u = \varprojlim \mathrm{Gal}(L_i/K)^u$, a closed subgroup of $\mathrm{Gal}(L/K)$.

Suppose $L/K$ is a possibly infinite Galois extension. A real number $u$ is a *jump* for $L/K$ if $\mathrm{Gal}(L/K)^{u+\varepsilon} \subsetneq \mathrm{Gal}(L/K)^u$ for every $\varepsilon > 0$. The *Hasse–Arf theorem* states that if $\mathrm{Gal}(L/K)$ is abelian, the jumps for $L/K$ are all integers. For a proof, see [Ser79, Chapter V, §7].

Finally we recall some of the main theorems of local class field theory, which furnishes a description of the abelianization $\mathrm{Gal}(\overline{K}/K)^{\mathrm{ab}}$. Assume first that the residue field of $K$ is finite. For an integer $i \geq 0$, let $U_K^i$ equal $\mathcal{O}_K^\times$ if $i = 0$, and $1 + \pi_K^i \mathcal{O}_K$ if $i > 0$. The $U_K^i$ are profinite topological groups.

**Theorem 4.1.3.** *Let $K^{\mathrm{ab}}/K$ be the maximal abelian extension, so that $\mathrm{Gal}(\overline{K}/K)^{\mathrm{ab}} \cong \mathrm{Gal}(K^{\mathrm{ab}}/K)$. Let $\hat{K}^\times$ be the profinite completion of the group $K^\times$, so that $\hat{K}^\times = \pi_K^{\hat{\mathbb{Z}}} \times \mathcal{O}_K^\times$. There exists an isomorphism (the* reciprocity map*)*

$$\mathrm{rec}_K \colon \hat{K}^\times \to \mathrm{Gal}(K^{\mathrm{ab}}/K)$$

*which carries $U_K^i$ onto $\mathrm{Gal}(K^{\mathrm{ab}}/K)^i$ for all integers $i \geq 0$.*

*For a finite extension $L/K$, the inclusion $K^{\mathrm{ab}} \subset L^{\mathrm{ab}}$ induces a homomorphism $\mathrm{Gal}(L^{\mathrm{ab}}/L) \to \mathrm{Gal}(K^{\mathrm{ab}}/K)$. The following diagram commutes:*

$$\begin{array}{ccc} L^\times & \xrightarrow{\mathrm{rec}_L} & \mathrm{Gal}(L^{\mathrm{ab}}/L) \\ \downarrow {\scriptstyle N_{L/K}} & & \downarrow \\ K^\times & \xrightarrow[\mathrm{rec}_K]{} & \mathrm{Gal}(K^{\mathrm{ab}}/K) \end{array}$$

For a proof, see [Ser67].

For an abelian group $A$ written multiplicatively, write $A^{(p)}$ for the image of $x \mapsto x^p$ on $A$.

**Corollary 4.1.4.** *Let $L/K$ be a (possibly infinite) abelian extension. For $u > e_K/(p-1)$, we have* $(\mathrm{Gal}(L/K)^u)^{(p)} = \mathrm{Gal}(L/K)^{u+e_K}$.

*Proof.* By the Hasse–Arf theorem, the jumps of $L/K$ are integers, and we are reduced to the case that $u = i$ is an integer. By Theorem 4.1.3, specifically the equality $\mathrm{rec}_K(U_K^i) = \mathrm{Gal}(K^{\mathrm{ab}}/K)^i$, it is enough to prove the relation $(U_K^i)^{(p)} = U_K^{i+e_K}$ whenever $i > e_K/(p-1)$. For $(U_K^i)^{(p)} \subset U_K^{i+e_K}$: Suppose $\alpha \in \mathcal{O}_K$. Expanding $(1+\pi_K^i\alpha)^p$ gives $1+\pi_K^{ip}\alpha^p$ plus terms divisible by $p\pi_K^i$. We have $v_K(p\pi_K^i) = i + e_K$, and $v_K(\pi_K^{ip}) = ip > i + e_K$ by hypothesis.

For $U_K^{i+e_K} \subset (U_K^i)^{(p)}$: Again suppose $\alpha \in \mathcal{O}_K$. We claim that the binomial expansion of $(1+p\pi_K^i\alpha)^{1/p}$ converges to an element of $U_K^i$. For $n \geq 1$, the $n$th term in the expansion is $p^n \binom{1/p}{n} \pi_K^{ni}\alpha^n$, which equals an element of $\mathcal{O}_K$ times $\pi_K^{ni}/n!$. Using the inequality $v(n!) \leq (n-1)/(p-1)$ (a standard inequality), we have:

$$v_K\left(\frac{\pi_K^{ni}}{n!}\right) - i = ni - e_K v(n!) - i \geq (n-1)(i - \frac{e_k}{p-1}).$$

This shows that the terms of the expansion of $(1+p\pi_K^i\alpha)^{1/p}$ beyond the initial 1 lie in $\pi^i\mathcal{O}_K$ and converge to 0. ☐

We remark that Corollary 4.1.4 holds even when the residue field of $K$ is infinite. One reduces to the case that the residue field is algebraically closed (since the higher ramification groups are insensitive to unramified base change). In that case Serre [Ser61] proves an analogue of Theorem 4.1.3 in which the $U_K^i$ are replaced by pro-algebraic groups $\mathcal{U}_K^i$. It is still true that the $p$th power map carries $\mathcal{U}_K^i$ surjectively onto $\mathcal{U}_K^{i+e_K}$ for $i$ sufficiently large, which is all that is necessary for Corollary 4.1.4.

Some of the results in this section pertain to the stability of some property under base change by a tame extension $L/K$. We list some properties of tame extensions in the following proposition.

**Proposition 4.1.5.** *Suppose $L/K$ is a finite Galois extension which is tame (that is, $p \nmid e_{L/K}$). Let $e = e_{L/K}$.*

1. *For $u \geq 0$ we have $\varphi_{L/K}(u) = u/e$ and $\psi_{L/K}(u) = eu$.*
2. *For $u > 0$ we have $N_{L/K}(U_L^{eu}) = U_K^u$.*
3. *Let $M/K$ be a (possibly infinite) abelian extension. For $u > 0$, the image of* $\mathrm{Gal}(LM/L)^{eu}$ *under* $\mathrm{Gal}(LM/L) \to \mathrm{Gal}(M/K)$ *is* $\mathrm{Gal}(M/K)^u$.

*Proof.* Part (1) follows from the definition of $\varphi_{L/K}$, using the facts that $\#\,\mathrm{Gal}(L/K)_0 = e$ and $\mathrm{Gal}(L/K)_u$ is trivial for all $u > 0$.

For part (2): The claim is transitive in the sense that if it is true for $L/K$ and $M/L$, then it is true for $M/K$. Using the fact that inertia groups of local fields are solvable [Ser79, Chapter IV, §2, Corollary 5], we can reduce to the case that $L/K$ is either unramified or totally ramified and cyclic of prime degree $\ell \neq p$. The proof of the claim in those cases is found in [Ser79, Chapter V, §2, Proposition 1] and [Ser79, Chapter V, §3, Corollary 4], respectively.

Part (3) is proved by combining compatibility of $\mathrm{rec}_K$ with the norm map (Theorem 4.1.3) with part (2). ☐

4.2. **Ramified $\mathbb{Z}_p$-extensions.** Keep the assumption that $K$ is a local field of characteristic $(0, p)$.

Let $K_\infty/K$ be a Galois extension with $\mathrm{Gal}(K_\infty/K) \cong \mathbb{Z}_p$. We assume that $K_\infty/K$ is ramified, in the sense that the inertia group $\mathrm{Gal}(K_\infty/K)^0$ is nontrivial. The inertia group is always closed, so this condition is equivalent to $\mathrm{Gal}(K_\infty/K)^0$ being of finite index in $\mathrm{Gal}(K_\infty/K)$, or equivalently, that $K_\infty$ contains elements of arbitrarily small positive valuation.

Such an extension always exists. For instance, let $L/K$ be the extension obtained by adjoining a primitive $p^n$th root of unity $\zeta_{p^n}$ for all $n \geq 1$. Since $v(\zeta_{p^n} - 1) = 1/p^n$ gets arbitrarily small, $L/K$ is ramified. Its Galois group is an open subgroup of $\mathbb{Z}_p^\times$, and so admits a maximal quotient isomorphic to $\mathbb{Z}_p$. The extension $K_\infty/K$ corresponding to this quotient is a ramified $\mathbb{Z}_p$-extension of $K$, which we call the *cyclotomic $\mathbb{Z}_p$-extension*.

The Galois group $\mathrm{Gal}(K_\infty/K) \cong \mathbb{Z}_p$ has an obvious filtration by subgroups $p^n\mathbb{Z}_p$ for $n = 0, 1, \dots$. We write $K_n/K$ for the finite extension corresponding to the subgroup $p^n\mathbb{Z}_p \subseteq \mathbb{Z}_p$.

On the other hand we have the upper numbering ramification filtration $\mathrm{Gal}(K_\infty/K)^u$ for real numbers $u \geq -1$. By the Hasse–Arf theorem, the jumps of this filtration occur at rational integers. That is, there is a sequence of integers $-1 = u_0 \leq u_1 \leq u_2 \leq \dots$ tending to $\infty$ such that (as a subgroup of $\mathrm{Gal}(K_\infty/K) \cong \mathbb{Z}_p$):

$$\mathrm{Gal}(K_\infty/K)^u \cong p^n\mathbb{Z}_p \text{ whenever } u_n < u \leq u_{n+1}$$

for all $n = 0, 1, \dots$. Let us call the sequence $u_0, u_1, u_2, \dots$ the *jumps* of $K_\infty/K$.

The two filtrations on $\mathrm{Gal}(K_\infty/K)$ are "eventually compatible" in the following sense.

**Lemma 4.2.1.** *Let $u_0, u_1, u_2, \dots$ be the sequence of jumps for $K_\infty/K$ defined above.*

*(1) There exists $N$ such that $u_N \geq 0$ and $u_{n+1} = u_n + e_K$ for all $n \geq N$.*

*(2) If $e_K = 1$ and $K_\infty/K$ is the cyclotomic $\mathbb{Z}_p$-extension, then $u_0 = -1$ and $u_i = i$ for all $i \geq 1$.*

*(3) Suppose $L/K$ is a finite Galois extension which is tame, and let $e = e_{L/K}$. Then $L_\infty = LK_\infty$ is a ramified $\mathbb{Z}_p$-extension of $L$. The sequence of jumps for $L_\infty/L$ is $-1, eu_1, eu_2, eu_3, \dots$.*

*(4) Suppose $e_K = 1$ and $L/K$ is a finite Galois extension which is tame. If $L_\infty/L$ is the cyclotomic $\mathbb{Z}_p$-extension, then the sequence of jumps for $L_\infty/L$ is $-1, e_L, 2e_L, 3e_L, \dots$.*

*Proof.* For part (1): Let $N$ be large enough so that $u_N > e_K/(p-1)$. By Corollary 4.1.4, for $u \geq u_N$ we have $(\mathrm{Gal}(K_\infty/K)^u)^{(p)} = \mathrm{Gal}(K_\infty/K)^{u+e_K}$. Now suppose $n \geq N$. If $u_n < u \leq u_{n+1}$, this relation means $\mathrm{Gal}(K_\infty/K)^{u+e_K} \cong p^{n+1}\mathbb{Z}_p$. This means precisely that $u_{n+1} = u_n + e_K$.

For part (2): The computation of the higher ramification groups of the cyclotomic $\mathbb{Z}_p$-extension over $\mathbb{Q}_p$ (or any unramified extension thereof) is standard, see for instance [Ser79, Chapter IV, §4].

For part (3): By Proposition 4.1.5, for all $u > 0$ the image of $\mathrm{Gal}(L_\infty/L)^{eu}$ under $\mathrm{Gal}(L_\infty/L) \to \mathrm{Gal}(K_\infty/K)$ is $\mathrm{Gal}(K_\infty/K)^u$, so that the jumps of $L_\infty/L$ are exactly $e$ times the jumps of $K_\infty/K$.

For part (4): Combine parts (2) and (3). □

**Lemma 4.2.2.** *Let $K_\infty/K$ be a ramified $\mathbb{Z}_p$-extension with jumps $u_0, u_1, u_2, \dots$. Then for each $n \geq 1$ large enough so that $K_n/K_{n-1}$ is ramified, the valuation of the different of $K_n/K_{n-1}$ is:*

$$v_{K_n}(\mathcal{D}_{K_n/K_{n-1}}) = \frac{e_{K_n/K}}{p^n} \sum_{k=0}^{n-1} (u_{k+1} - u_k)(p^{k+1} - p^k).$$

*Proof.* Applying the formula for the different (4.1.2) to $K_n/K$, we find:

$$v_{K_n}(\mathcal{D}_{K_n/K}) = e_{K_n/K} \sum_{k=0}^{n-1} (u_{k+1} - u_k)\left(1 - \frac{1}{p^{n-k}}\right)$$

By the transitivity of the different [Ser79, Chapter III, §4, Proposition 8], we have $\mathcal{D}_{K_n/K} = \mathcal{D}_{K_n/K_{n-1}}\mathcal{D}_{K_{n-1}/K}$. Applying $v_{K_n}$ (and using the relation $v_{K_{n-1}} = p^{-1}v_{K_n}$, true by our assumption that $K_n/K_{n-1}$ is ramified) gives the result. □

The next lemma adapts the statement and proof of [Tat67, Proposition 5] to give a precise formula for the limiting behavior of the ramification growth as one ascends a ramified $\mathbb{Z}_p$-extension.

**Lemma 4.2.3.** *Let $K_\infty/K$ be a ramified $\mathbb{Z}_p$-extension. There exists an integer $N \geq 0$ such that for all $n \geq N+1$ we have:*

$$v_{K_n}(\mathfrak{D}_{K_n/K_{n-1}}) \geq p - 1 + e_{K_n}(1 - p^{N-n}).$$

*Proof.* Take for $N$ the integer appearing in Lemma 4.2.1(1), so that $u_N \geq 0$ (this means that $K_\infty/K_N$ is totally ramified) and $u_n = u_{n-1} + e_K$ for all $n \geq N+1$. For such an $n$, $K_n/K_{n-1}$ is totally ramified, and so Lemma 4.2.2 applies:

$$\begin{aligned} v_{K_n}(\mathfrak{D}_{K_n/K_{n-1}}) &= \frac{e_{K_n/K}}{p^n}\sum_{k=0}^{n-1}(u_{k+1}-u_k)(p^{k+1}-p^k) \\ &= \frac{e_{K_n/K}}{p^n}\sum_{k=0}^{N-1}(u_{k+1}-u_k)(p^{k+1}-p^k) + \frac{e_{K_n/K}}{p^n}\sum_{k=N}^{n-1}e_K(p^{k+1}-p^k) \\ &= \frac{e_{K_n/K}}{p^n}\sum_{k=0}^{N-1}(u_{k+1}-u_k)(p^{k+1}-p^k) + e_{K_n}(1-p^{N-n}) \end{aligned}$$

Suppose that the inertia group $\mathrm{Gal}(K_\infty/K)^0 \cong p^r\mathbb{Z}_p$ as a subgroup of $\mathrm{Gal}(K_\infty/K) \cong \mathbb{Z}_p$; this means exactly that $u_0 = u_1 = \cdots = u_r = -1$ and $u_{r+1} \geq 0$, and also that $p^r e_{K_n/K} = p^n$. Since $K_n/K_{n-1}$ is ramified, we must have $r < N$. Therefore the first sum above is bounded below by $(e_{K_n/K}/p^n)(p^{r+1} - p^r) = p - 1$. $\square$

In the context of Lemma 4.2.3, $K_\infty/K_N$ is a totally ramified $\mathbb{Z}_p$-extension. If it is relabeled so that $K = K_N$ is the base of the extension, then each $K_n/K$ is totally ramified, and the inequality in that lemma becomes

$$v_{K_n}(\mathfrak{D}_{K_n/K_{n-1}}) \geq e_K(p^n - 1) + (p-1). \tag{4.2.4}$$

**Definition 4.2.5.** Let $K_\infty/K$ be a ramified $\mathbb{Z}_p$-extension. We call $K_\infty/K$ *sufficiently ramified* if the inequality in (4.2.4) holds for each $n \geq 1$.

In light of the above considerations and Lemma 4.2.1, there are plenty of examples of sufficiently ramified $\mathbb{Z}_p$-extensions.

**Lemma 4.2.6.** *Sufficiently ramified $\mathbb{Z}_p$-extensions exist in the following situations:*

*(1) Let $K$ be a local field of characteristic $(0,p)$ such that $p \nmid e_K$. Then $K$ admits a sufficiently ramified $\mathbb{Z}_p$-extension.*
*(2) Let $K$ be as in part (1), and let $L/K$ be a finite subextension of the cyclotomic $\mathbb{Z}_p$-extension over $L$. Then $L$ admits a sufficiently ramified $\mathbb{Z}_p$-extension.*
*(3) Let $K$ be any local field of characteristic $(0,p)$, and let $K_\infty/K$ be a ramified $\mathbb{Z}_p$-extension. Then for sufficiently large $n$, $K_\infty/K_n$ is a sufficiently ramified $\mathbb{Z}_p$-extension.*
*(4) Let $K_\infty/K$ be a sufficiently ramified $\mathbb{Z}_p$-extension, and let $L/K$ be a finite Galois extension which is tame. Let $L_\infty = K_\infty L$. Then $L_\infty/L$ is a sufficiently ramified $\mathbb{Z}_p$-extension.*

We continue to follow [Tat67] very closely, taking care to give precise bounds. The next phase is to study the behavior of traces.

**Lemma 4.2.7** ([Ser79, Chapter V Lemma 4])**.** *Let $L/K$ be a cyclic extension of order $p$. Then for any $\alpha \in L$ we have:*

$$v_K(\mathrm{tr}_{L/K}(\alpha)) \geq \left\lfloor \frac{v_L(\alpha) + v_L(\mathfrak{D}_{L/K})}{p} \right\rfloor,$$

*where $\lfloor x \rfloor$ is the floor function.*

The next lemma is an adaptation of [Tat67, Corollaries 2 and 3]. Let us write $|x| = a^{v(x)}$ for $x \in K_\infty$ and an arbitrary real number $0 < a < 1$. Then $|x|$ is an absolute value inducing the topology on $K_\infty$.

**Lemma 4.2.8.** *Assume $K_\infty/K$ is sufficiently ramified. We have the following inequalities.*

*(1) For $x \in K_n$ we have:*

$$\left|\mathrm{tr}_{K_n/K_{n-1}}(x)\right| \leq |p|^{1-p^{-n}} |x|$$

*(2) For $x \in K_n$ we have:*

$$\left|\mathrm{tr}_{K_n/K}(x)\right| \leq |p|^{n-(p-1)^{-1}} |x|$$

*Proof.* For part (1): Let $d = v_{K_n}(\mathfrak{D}_{K_n/K_{n-1}})$. By definition of sufficiently ramified we have $d \geq e_K(p^n-1)+p-1$. Applying Lemma 4.2.7 to $K_{n+1}/K_n$, we find:

$$\begin{aligned} v_{K_{n-1}}(\mathrm{tr}_{K_n/K_{n-1}}(x)) &\geq \left\lfloor \frac{v_{K_n}(x)+d}{p} \right\rfloor \\ &\geq \frac{v_{K_n}(x)+d-(p-1)}{p} \\ &= e_K p^{n-1} - \frac{e_K}{p} + \frac{v_{K_n}(x)}{p}. \end{aligned}$$

Translating this inequality in terms of the absolute valuation gives:

$$v(\mathrm{tr}_{K_n/K_{n-1}}(x)) \geq 1 - \frac{1}{p^n} + v(x),$$

which (recalling that $v(p)=1$) means exactly the inequality in part (1).

For part (2), we apply part (1) inductively on $x \in K_n$ to obtain:

$$\begin{aligned} \left|\mathrm{tr}_{K_n/K}(x)\right| &= \left|\mathrm{tr}_{K_1/K}\,\mathrm{tr}_{K_2/K_1}\cdots\mathrm{tr}_{K_n/K_{n-1}}(x)\right| \\ &\leq |p|^{n-p^{-1}-p^{-2}-\cdots-p^{-n}} |x| \\ &\leq |p|^{n-(p-1)^{-1}} |x| \end{aligned}$$

□

Keep the assumption that $K_\infty/K$ is sufficiently ramified. Define the normalized trace $t\colon K_\infty \to K$ by $t(x) = p^{-n}\,\mathrm{tr}_{K_n/K}(x)$ whenever $x \in K_n$ (note that this is unambiguous). For $n \geq 1$ we define $t_n\colon K_n \to K_{n-1}$ by $t_n(x) = p^{-1}\,\mathrm{tr}_{K_n/K_{n-1}}(x)$. Then Lemma 4.2.8 translates into:

$$|t(x)| \leq |p|^{-(p-1)^{-1}} |x|\,,\ x \in K_\infty \tag{4.2.9}$$

and

$$|t_n(x)| \leq |p|^{-p^{-n}} |x|\,. \tag{4.2.10}$$

Finally note that $t(t_n(x)) = t(x)$ for all $x \in K_n$.

Let $\sigma$ be a topological generator of $\mathrm{Gal}(K_\infty/K) \cong \mathbb{Z}_p$. Elements of the polynomial ring $K[\sigma]$ act on $K_\infty$ in the evident manner. For $x \in K_n$ we record the relation

$$t_n(x) = p^{-1}(1+\sigma^{p^{n-1}}+\sigma^{2p^{n-1}}+\cdots+\sigma^{(p-1)p^{n-1}})(x). \tag{4.2.11}$$

The next lemma adapts [Tat67, Proposition 6].

**Lemma 4.2.12.** *Assume that $K_\infty/K$ is sufficiently ramified. For $x \in K_\infty$ we have:*

$$|x - t(x)| \leq |p|^{-1-\frac{1}{p(p-1)}} |\sigma(x) - x|$$

*Proof.* For each $x \in K_n$, (4.2.11) implies:

$$\begin{aligned} x - t_n(x) &= p^{-1}\sum_{i=1}^{p-1}(1-\sigma^{ip^{n-1}})(x) \\ &= p^{-1}\sum_{i=1}^{p-1}(1+\sigma^{p^{n-1}}+\cdots+\sigma^{(i-1)p^{n-1}})(1-\sigma^{p^{n-1}})(x). \end{aligned}$$

Noting that $|\alpha(x)| \leq |x|$ whenever $\alpha \in \mathcal{O}_K[\sigma]$, we find:

$$|x - t_n(x)| \leq |p|^{-1}\left|(1-\sigma^{p^{n-1}})x\right| \leq |p|^{-1}\left|\sigma(x) - x\right|. \tag{4.2.13}$$

We will prove by induction on $n \geq 1$ the following statement which implies the lemma: for $x \in K_n$, we have

$$|x - t(x)| \leq |p|^{-1-\frac{1}{p^2}-\cdots-\frac{1}{p^n}}\left|\sigma(x) - x\right|. \tag{4.2.14}$$

(The exponent of $|p|$ is meant to be $-1$ in the case $n = 1$.) The base case $n = 1$ is an instance of (4.2.13). Now suppose $n \geq 2$. Assume (4.2.14) for $n-1$, and then for $x \in K_n$ we have:

$$|x - t(x)| \leq \max\left\{|x - t_n(x)|, |t_n(x) - t(x)|\right\}$$

Treating each quantity on the right side in turn, we have

$$|x - t_n(x)| \leq |p|^{-1}\left|\sigma(x) - x\right|$$

by (4.2.13); this is $\leq$ the right side of (4.2.14). The other quantity is

$$\begin{aligned} |t_n(x) - t(x)| &= |t_n(x) - t(t_n(x))| \\ &\leq |p|^{-1-\frac{1}{p^2}-\cdots-\frac{1}{p^{n-1}}}\left|(\sigma-1)t_n(x)\right| \end{aligned}$$

by applying the inductive hypothesis to $t_n(x) \in K_{n-1}$. Noting that $\sigma$ commutes with $t_n$, we have

$$|(\sigma-1)t_n(x)| = |t_n(\sigma-1)(x)| \leq |p|^{-p^{-n}}\left|\sigma(x) - x\right|$$

by (4.2.10) and so

$$|t_n(x) - t(x)| \leq |p|^{-1-\frac{1}{p^2}-\cdots-\frac{1}{p^n}}\left|\sigma(x) - x\right|$$

as required. □

The bounds in (4.2.9) and Lemma 4.2.12 show that $t\colon K_\infty \to K$ is continuous and extends uniquely to a $K$-linear map $t\colon \hat{K}_\infty \to K$ satisfying the same bounds for all $x \in \hat{K}_\infty$:

$$\begin{aligned} |t(x)| &\leq |p|^{-\frac{1}{p-1}}\,|x| &(4.2.15)\\ |x - t(x)| &\leq |p|^{-1-\frac{1}{p(p-1)}}\,|\sigma(x) - x| &(4.2.16) \end{aligned}$$

The following lemma is an integral version of [Tat67, Proposition 7a].

**Lemma 4.2.17.** *Assume that $K_\infty/K$ is sufficiently ramified. Let $\mathcal{O}^{t=0}_{\hat{K}_\infty}$ denote the kernel of the normalized trace map on $\mathcal{O}_{\hat{K}_\infty}$. The cokernel of the inclusion map*

$$\mathcal{O}_K \oplus \mathcal{O}^{t=0}_{\hat{K}_\infty} \hookrightarrow \mathcal{O}_{\hat{K}_\infty}$$

*is killed by any $\alpha \in \mathcal{O}_K$ with $|\alpha| \leq |p|^{1/(p-1)}$.*

*Proof.* Given $x \in \mathcal{O}^{t=0}_{\hat{K}_\infty}$, we may write $x = (x - t(x)) + t(x)$; after multiplying by an element $\alpha$ as in the Lemma, $\alpha t(x)$ lies in $\mathcal{O}_K$ by (4.2.15), and then $\alpha(x - t(x))$ lies in $\mathcal{O}_{\hat{K}_\infty}$. □

We turn our attention now to the continuous cohomology $H^i_{\mathrm{cts}}(\mathrm{Gal}(K_\infty/K), \mathcal{O}_{\hat{K}_\infty})$ when $K_\infty/K$ is a sufficiently ramified $\mathbb{Z}_p$-extension. Since $\mathrm{Gal}(K_\infty/K) \cong \mathbb{Z}_p$, the continuous cohomology of $\mathrm{Gal}(K_\infty/K)$ acting continuously on a $p$-adically complete abelian group $M$ is computed by the complex

$$M^{h\,\mathrm{Gal}(K_\infty/K)}: \quad M \xrightarrow{\sigma-1} M\ ,$$

which is to say that $H^0_{\mathrm{cts}}(M) = M^\sigma$, $H^1_{\mathrm{cts}}(M) = M/(\sigma-1)M$, and $H^i_{\mathrm{cts}}(M) = 0$ for $i \geq 2$. In the case of $M = \mathcal{O}_{\hat{K}_\infty}$, the $H^0$ is $\mathcal{O}_K$ by Ax's theorem (Theorem 4.0.1), so we are left with describing the $H^1$.

The following proposition adapts [Tat67, Proposition 7b] to the setting of integral cohomology.

**Proposition 4.2.18.** *Let $K_\infty/K$ be a ramified $\mathbb{Z}_p$-extension. Let $N \geq 0$ be large enough so that $K_\infty/K_N$ is sufficiently ramified. Define an abelian group $X$ with* $\mathrm{Gal}(K_\infty/K)$*-action by the exact sequence:*

$$0 \to \mathcal{O}_K \to \mathcal{O}_{\hat{K}_\infty} \to X \to 0.$$

*Then $X$ is $p$-adically complete, $H^i_{\mathrm{cts}}(\mathrm{Gal}(K_\infty/K), X) = 0$ for all $i \neq 1$, and $H^1_{\mathrm{cts}}(\mathrm{Gal}(K_\infty/K), X)$ is $p^{N+2}$-torsion for $p \neq 2$ and $p^{N+3}$-torsion for $p = 2$.*

*Remark* 4.2.19. Proposition 4.2.18 may be stated in terms of complexes this way: Let $\mathcal{X} = X^{h\,\mathrm{Gal}(K_\infty/K)}$, so that $\mathcal{X}$ sits in an exact triangle

$$\mathcal{O}_K^{h\,\mathrm{Gal}(K_\infty/K)} \to \mathcal{O}_{\hat{K}_\infty}^{h\,\mathrm{Gal}(K_\infty/K)} \to \mathcal{X}.$$

Then we have $H^i(\mathcal{X}) = 0$ for $i \neq 1$, and $H^1(\mathcal{X}) = H^1_{\mathrm{cts}}(\mathrm{Gal}(K_\infty/K), X)$ is a torsion group as claimed in the proposition. Note that since $\mathrm{Gal}(K_\infty/K)$ acts trivially on $\mathcal{O}_K$, a choice of isomorphism $\mathrm{Gal}(K_\infty/K) \cong \mathbb{Z}_p$ induces a quasi-isomorphism $\mathcal{O}_K^{h\,\mathrm{Gal}(K_\infty/K)} \cong \mathcal{O}_K[\varepsilon]$.

*Proof.* By Lemma 4.2.17, $X$ is an extension of a $p$-torsion group by the $p$-complete group $\mathcal{O}^{t=0}_{\hat{K}_\infty}$, and therefore it is $p$-complete.

Regarding the claims about $H^i_{\mathrm{cts}}(\mathrm{Gal}(K_\infty/K, X)$: Since $\mathrm{Gal}(K_\infty/K) \cong \mathbb{Z}_p$ has cohomological dimension 1, it is enough to consider $H^0$ and $H^1_{\mathrm{cts}}$.

*We first consider the case* $N = 0$, meaning that $K_\infty/K$ is sufficiently ramified. Let $\sigma$ be a topological generator for $\mathrm{Gal}(K_\infty/K)$. The claim that $H^0(\mathrm{Gal}(K_\infty/K), X) = 0$ is equivalent to the claim that $\mathcal{O}_K \to H^1_{\mathrm{cts}}(\mathrm{Gal}(K_\infty/K), \mathcal{O}_{\hat{K}_\infty})$ is injective. If $\alpha \in \mathcal{O}_K$ lies in the kernel, then $\alpha = (\sigma-1)\alpha'$ for $\alpha' \in \mathcal{O}_{\hat{K}_\infty}$, but then $\alpha = t(\alpha) = t(\sigma-1)\alpha'$. We claim that $t(\sigma-1)$ acts as 0 on $\hat{K}_\infty$. Indeed, $t(\sigma-1)$ is zero on each $K_n$, hence on $K_\infty$, and so by continuity it is zero on all of $\hat{K}_\infty$.

We now turn to $H^1_{\mathrm{cts}}(\mathrm{Gal}(K_\infty/K), \mathcal{O}_{\hat{K}_\infty})$. Let us write the superscript $t = 0$ to mean the kernel of the normalized trace $t$ wherever this is defined. For each $n \geq 1$, the $K$-linear operator $\sigma - 1$ is injective on the finite-dimensional $K$-vector space $K_n^{t=0}$, since the kernel is $K \cap K_n^{t=0} = 0$. (The operator $t$ is the identity on $K$.) Therefore the restriction of $\sigma-1$ to $K_n^{t=0}$ is an isomorphism; write $(\sigma-1)^{-1}$ for its inverse. The inequality (4.2.16) shows that $(\sigma-1)^{-1}$ is bounded on $K_\infty^{t=0}$ with operator norm $\leq |p|^{-1-1/p(p-1)}$, so it extends to an operator $(\sigma-1)^{-1}$ on $\hat{K}_\infty^{t=0}$ with the same operator norm. Therefore

$$H^1_{\mathrm{cts}}(\mathrm{Gal}(K_\infty/K), \mathcal{O}^{t=0}_{\hat{K}_\infty}) \cong \frac{\mathcal{O}^{t=0}_{\hat{K}_\infty}}{(\sigma-1)\mathcal{O}^{t=0}_{\hat{K}_\infty}}$$

is $\alpha$-torsion for any $\alpha \in \mathcal{O}_K$ with $v(\alpha) \geq 1 + 1/p(p-1)$.

Consider the following diagram of $p$-adically complete groups, in which both rows are exact:

$$\begin{array}{ccccccccc} 0 & \longrightarrow & \mathcal{O}_K & \longrightarrow & \mathcal{O}_{\hat{K}_\infty} & \longrightarrow & X & \longrightarrow & 0 \\ & & \downarrow & & {=}\downarrow & & \downarrow & & \\ 0 & \longrightarrow & \mathcal{O}_K \oplus \mathcal{O}^{t=0}_{\hat{K}_\infty} & \longrightarrow & \mathcal{O}_{\hat{K}_\infty} & \longrightarrow & W & \longrightarrow & 0 \end{array}$$

Lemma 4.2.17 says that $W$ is $\beta$-torsion for any $\beta \in \mathcal{O}_K$ with $v(\beta) \geq 1 + 1/p(p-1)$.

The long exact sequences in cohomology associated to the above exact sequences read in part:

$$\begin{array}{ccccccccc} 0 & \longrightarrow & H^1_{\mathrm{cts}}(\mathcal{O}_K) & \longrightarrow & H^1_{\mathrm{cts}}(\mathcal{O}_{\hat{K}_\infty}) & \longrightarrow & H^1_{\mathrm{cts}}(X) & \longrightarrow & 0 \\ & & \downarrow & & {=}\downarrow & & \downarrow & & \\ 0 & \longrightarrow & \left(H^1_{\mathrm{cts}}(\mathcal{O}_K) \oplus H^1_{\mathrm{cts}}(\mathcal{O}^{t=0}_{\hat{K}_\infty})\right)/H^0(W) & \longrightarrow & H^1_{\mathrm{cts}}(\mathcal{O}_{\hat{K}_\infty}) & \longrightarrow & H^1_{\mathrm{cts}}(W) & \longrightarrow & 0 \end{array}$$

By the snake lemma, there is an exact sequence

$$0 \to H^1_{\mathrm{cts}}(\mathcal{O}^{t=0}_{\hat{K}_\infty})/H^0(W) \to H^1_{\mathrm{cts}}(X) \to H^1_{\mathrm{cts}}(W) \to 0$$

As $H^1_{\mathrm{cts}}(\mathcal{O}^{t=0}_{\hat{K}_\infty})$ is killed by any element $\alpha \in \mathcal{O}_K$ with $v(\alpha) \geq 1 + 1/p(p-1)$, and $H^1_{\mathrm{cts}}(W)$ is killed by any element $\beta$ with $v(\beta) \geq 1/(p-1)$, we find that $H^1_{\mathrm{cts}}(X)$ is killed by any product $\alpha\beta$ of such elements.

If $e_K$ is sufficiently large (to wit $e_K \geq p(p-1)$), then any $\gamma \in \mathcal{O}_K$ with $v(\gamma) \geq 1 + 1/p(p-1) + 1/(p-1) = 1 + (p+1)/p(p-1)$ can be factored into such a product $\alpha\beta$, and the lemma would follow for $K$. But if for instance $e_K = 1$, we could only conclude from the above that $H^1_{\mathrm{cts}}(X)$ is $p^3$-torsion.

To improve the result as stated in the proposition, let $L/K$ be a sufficiently ramified tame Galois extension with group $G$. (For example, $L$ could be the splitting field of $x^n - \pi_K$, where $\pi_K$ is a uniformizer of $K$ and $n$ is sufficiently large and prime to $p$.) Since each $K_n/K$ is totally wildly ramified, $K_n$ and $L$ are linearly disjoint over $K$. Let $L_\infty = K_\infty L$, so that $L_\infty/L$ is a sufficiently ramified $\mathbb{Z}_p$-extension by Lemma 4.2.1(3). We abuse notation by referring to both $\mathrm{Gal}(K_\infty/K)$ and $\mathrm{Gal}(L_\infty/L)$ as $\mathbb{Z}_p$, and we identify $\mathrm{Gal}(L_\infty/K) = G \times \mathbb{Z}_p$. Define a complex $\mathcal{X}_L$ with $G$-action by the exact triangle

$$\mathcal{O}_L^{h\mathbb{Z}_p} \to \mathcal{O}_{\hat{L}_\infty}^{h\mathbb{Z}_p} \to \mathcal{X}_L.$$

Applying derived $G$-invariants, we obtain an exact triangle

$$\mathcal{O}_L^{h(G\times\mathbb{Z}_p)} \to \mathcal{O}_{\hat{L}_\infty}^{h(G\times\mathbb{Z}_p)} \to \mathcal{X}_L^{hG}$$

Since $L_n/K_n$ is a tame extension with group $G$, Noether's theorem states that $\mathcal{O}_{L_n}$ is a free $\mathcal{O}_{K_n}[G]$-module. For every $m$, $\mathcal{O}_{L_n}/p^m$ is a free $(\mathcal{O}_{K_n}/p^m)[G]$-module, and so $H^i(G, \mathcal{O}_{L_n}/p^m) = 0$ for all $i > 0$. Therefore

$$H^i(G, \mathcal{O}_{\hat{L}_\infty}) = \varprojlim_m \varinjlim_n H^i(G, \mathcal{O}_{L_n}/p^m) = 0,$$

meaning that $\mathcal{O}^G_{\hat{L}_\infty} = H^0(G, \mathcal{O}_{\hat{L}_\infty}) = \mathcal{O}_{\hat{K}_\infty}$. Our exact triangle now reads

$$\mathcal{O}_K^{h\mathbb{Z}_p} \to \mathcal{O}_{\hat{K}_\infty}^{h\mathbb{Z}_p} \to \mathcal{X}_L^{hG},$$

so that $\mathcal{X} = \mathcal{X}_L^{hG}$. By the argument above, $H^i(\mathcal{X}_L) = 0$ for all $i \neq 1$, and $H^1(\mathcal{X}_L)$ is $p^2$- or $p^3$-torsion as $p$ is odd or even. By the Lyndon–Hochschild–Serre spectral sequence Proposition 3.3.6, we have $H^1(\mathcal{X}) = H^0(G, H^1(\mathcal{X}_L))$, which is subject to the same bounds.

*Now consider the general case:* Let $N$ be large enough so that $K_\infty/K_N$ is sufficiently ramified. Let $X = \mathcal{O}_{\hat{K}_\infty}/\mathcal{O}_K$ and $X_N = \mathcal{O}_{\hat{K}_\infty}/\mathcal{O}_{K_N}$, so that we have an exact sequence of $p$-adically complete abelian groups with $\mathrm{Gal}(K_\infty/K)$-action:

$$0 \to \mathcal{O}_{K_N}/\mathcal{O}_K \to X \to X_N \to 0. \tag{4.2.20}$$

By the above argument we have that $H^0(\mathrm{Gal}(K_\infty/K_N), X_N) = 0$ and $H^1_{\mathrm{cts}}(\mathrm{Gal}(K_\infty/K_N), X_N)$ is $p^2$-torsion or $p^3$-torsion, as $p$ is odd or even.

Consider the $H^0$ part of the long exact sequence associated to (4.2.20):

$$H^0(\mathrm{Gal}(K_\infty/K), \mathcal{O}_{K_N}/\mathcal{O}_K) \to H^0(\mathrm{Gal}(K_\infty/K), X) \to H^0(\mathrm{Gal}(K_\infty/K), X_N).$$

We have $H^0(\mathrm{Gal}(K_\infty/K), \mathcal{O}_{K_N}/\mathcal{O}_K) = H^0(\mathrm{Gal}(K_N/K), \mathcal{O}_{K_N}/\mathcal{O}_K) = 0$ (for instance, because it injects into $H^1(\mathrm{Gal}(K_N/K), \mathcal{O}_K) = 0$). Also $H^0(\mathrm{Gal}(K_\infty/K), X_N) \subset H^0(\mathrm{Gal}(K_\infty/K_N), X_N) = 0$. Therefore $H^0(\mathrm{Gal}(K_\infty/K), X) = 0$.

Now consider the $H^1$ part of the long exact sequence associated to (4.2.20), relative to the group $\mathrm{Gal}(K_\infty/K_N)$:

$$H^1_{\mathrm{cts}}(\mathrm{Gal}(K_\infty/K_N), \mathcal{O}_{K_N}/\mathcal{O}_K) \to H^1_{\mathrm{cts}}(\mathrm{Gal}(K_\infty/K_N), X) \to H^1_{\mathrm{cts}}(\mathrm{Gal}(K_\infty/K_N), X_N) \tag{4.2.21}$$

The term $H^1_{\mathrm{cts}}(\mathrm{Gal}(K_\infty/K_N), \mathcal{O}_{K_N}/\mathcal{O}_K)$ is isomorphic to $\mathcal{O}_{K_N}/\mathcal{O}_K$ itself. The terms in (4.2.21) carry a residual action of $\mathrm{Gal}(K_N/K)$; taking invariants reveals that $H^1_{\mathrm{cts}}(\mathrm{Gal}(K_\infty/K_N), X)^{\mathrm{Gal}(K_N/K)}$ injects into $H^1_{\mathrm{cts}}(\mathrm{Gal}(K_\infty/K_N), X_N)$ and is therefore $p^2$- or $p^3$-torsion.

Finally, consider the inflation-restriction sequence which calculates $H^1_{\mathrm{cts}}(\mathrm{Gal}(K_\infty/K), X)$:

$$H^1(\mathrm{Gal}(K_N/K), X^{\mathrm{Gal}(K_\infty/K_N)}) \to H^1_{\mathrm{cts}}(\mathrm{Gal}(K_\infty/K), X) \to H^1_{\mathrm{cts}}(\mathrm{Gal}(K_\infty/K_N), X)^{\mathrm{Gal}(K_N/K)}$$

The term on the left is $p^N$-torson (this being the order of $\mathrm{Gal}(K_N/K)$). Therefore the middle term is $p^{N+2}$- or $p^{N+3}$-torsion as $p$ is odd or even. $\square$

4.3. **The Galois cohomology of $\mathcal{O}_C$.** Keep the assumption that $K$ is a local field of characteristic $(0, p)$. Let $\overline{K}$ be an algebraic closure, and let $C$ be the metric completion of $\overline{K}$. The Galois group $\Gamma_K = \mathrm{Gal}(\overline{K}/K)$ acts continuously on $\overline{K}$ and so extends to a continuous action on $C$ and on $\mathcal{O}_C$. We are interested in $\mathcal{O}_C^{h\Gamma_K}$, the complex which computes the continuous cohomology $H^*_{\mathrm{cts}}(\Gamma_K, \mathcal{O}_C)$. Tate's idea is to use a ramified $\mathbb{Z}_p$-extension $K_\infty/K$ as an intermediary:

$$\mathcal{O}_C^{h\Gamma_K} \cong \left(\mathcal{O}_C^{h\,\mathrm{Gal}(\overline{K}/K_\infty)}\right)^{h\,\mathrm{Gal}(K_\infty/K)}. \tag{4.3.1}$$

First we deal with the inner term on the right side of (4.3.1), relating to the continuous cohomology $H^*_{\mathrm{cts}}(\mathrm{Gal}(\overline{K}/K_\infty), \mathcal{O}_C)$.

The next lemma is the integral version of [Tat67, Theorem 1].

**Lemma 4.3.2.** *Define a complex $\mathcal{Y}_0$ of solid $\mathcal{O}_{\hat{K}_\infty}$-modules by the exact triangle*

$$\mathcal{O}_{\hat{K}_\infty} \to \mathcal{O}_C^{h\,\mathrm{Gal}(\overline{K}/K_\infty)} \to \mathcal{Y}_0, \tag{4.3.3}$$

*where the first morphism is induced from $\mathcal{O}_{\hat{K}_\infty} \hookrightarrow \mathcal{O}_C$. Then $H^0(\mathcal{Y}_0) = 0$, and the complex $\mathcal{Y}_0$ is almost zero.*

*Proof.* The proof of [Tat67, Theorem 1] in fact proves this stronger result. From the modern point of view, the result is a special case $U = \mathrm{Spa}(\hat{K}_\infty, \mathcal{O}_{\hat{K}_\infty})$ of the aforementioned result [Sch13a, Lemma 4.10], which states that whenever $U$ is an affinoid perfectoid space, $H^i(U_{\text{proét}}, \hat{\mathcal{O}}^+)$ is almost zero for $i \geq 1$.

We need to justify why $\hat{K}_\infty$ is perfectoid. The extension $K_\infty/K$ has open higher ramification groups, which means it is arithmetically profinite in the sense of [Win83]. For any such extension, $\hat{K}_\infty$ is perfectoid. Indeed, the *field of norms* construction of Fontaine-Wintenberger implicitly shows that the $p$th power map on $\mathcal{O}_{K_\infty}$ is surjective modulo some pseudo-uniformizer of $\hat{K}_\infty$.

The claim about $H^0(\mathcal{Y}_0)$ is Ax's result applied to $\hat{K}_\infty$. □

The complex $\mathcal{Y}_0$ in Lemma 4.3.2 admits a residual $\mathrm{Gal}(K_\infty/K)$-action, so we may define the derived invariants $\mathcal{Y} := \mathcal{Y}_0^{h\,\mathrm{Gal}(K_\infty/K)}$, which fit into an exact triangle:

$$\mathcal{O}_{\hat{K}_\infty}^{h\,\mathrm{Gal}(K_\infty/K)} \to \mathcal{O}_C^{h\Gamma_K} \to \mathcal{Y}. \tag{4.3.4}$$

**Lemma 4.3.5.** *We have $H^0(\mathcal{Y}) = 0$, and for all $i \geq 1$, $H^i(\mathcal{Y})$ is a vector space over the residue field of $K$. In particular, $H^i(\mathcal{Y})$ is $p$-torsion.*

*Proof.* The claim about $H^0(\mathcal{Y})$ follows by taking $\mathrm{Gal}(K_\infty/K)$-invariants in $H^0(\mathcal{Y}_0)$.

The claim about $H^i(\mathcal{Y})$ for $i \geq 1$ comes from a general phenomenon in almost mathematics. Let $R = \mathcal{O}_{\hat{K}_\infty}$, and let $I \subset \mathcal{O}_{\hat{K}_\infty}$ be the maximal ideal. Let $R$ be a ring, and let $I \subset R$ be an ideal. The key assumptions we need are that $I$ is a flat $R$-module and $I^2 = I$. Those conditions ensure that the commutative algebra object $R/I$ of the derived category $D(R)$ is idempotent [Lur17, Definition 4.8.2.8], meaning that the multiplication map $R/I \otimes_R^{\mathbb{L}} R/I \to R/I$ is an isomorphism.

By a general formalism [Lur17, Proposition 4.8.2.4], the base change morphism $D(R) \to D(R/I)$ is localizing, meaning that its right adjoint, namely the pullback functor $D(R/I) \to D(R)$, is fully faithful. Call an object $A$ of $D(R)$ almost zero if its cohomology is $I$-torsion. Evidently if $A$ is an object of $D(R/I)$ then its image in $D(R)$ is almost zero; we claim the converse as well. Suppose $A$ is almost zero. We claim $A \to A \otimes_R^{\mathbb{L}} R/I$ is a quasi-isomorphism; by the Tor spectral sequence, this reduces to the claim that $\mathrm{Tor}_R^i(H^j(A), R/I)$ is $H^j(A)$ for $i = 0$ and is 0 otherwise. This in turn is easily deduced from the flat resolution $I \to R$ of $R/I$, and from the fact that $H^j(A)$ is $I$-torsion.

We have shown that $D(R/I) \to D(R)$ is a fully faithful functor whose essential image is the almost zero objects. In our situation, we need an enriched variation: the pullback functor induces an equivalence between the derived category of solid $k$-vector spaces (here $k$ is the residue field of $K_\infty$) with semilinear $\mathrm{Gal}(K_\infty/K)$-action, and the category of almost zero objects in the derived category of solid $\mathcal{O}_{\hat{K}_\infty}$-modules with semilinear $\mathrm{Gal}(K_\infty/K)$-action. To prove the variation, one only needs to verify that $I$ is flat in the solid category; indeed it is the colimit over the nonunit principal ideals of $\mathcal{O}_{\hat{K}_\infty}$, each of which is flat.

In the situation of the lemma, $\mathcal{Y}$ is almost zero, so it is quasi-isomorphic to a complex of solid $k$-vector spaces with $\mathrm{Gal}(K_\infty/K)$-action. Therefore $\mathcal{Y} = \mathcal{Y}_0^{h\,\mathrm{Gal}(K_\infty/K)}$ is a complex of $\mathbb{F}_p$-vector spaces, and in particular $H^i(\mathcal{Y})$ is $p$-torsion for all $i$.

□

*Remark* 4.3.6. See [HS24] for a discussion of yet more general hypotheses under which one has a good theory of homological algebra in the almost setting.

We have a composition of morphisms:

$$\mathcal{O}_K^{h\,\mathrm{Gal}(K_\infty/K)} \to \mathcal{O}_{\hat{K}_\infty}^{h\,\mathrm{Gal}(K_\infty/K)} \to \mathcal{O}_C^{h\Gamma_K} \tag{4.3.7}$$

where the first is induced by the inclusion $\mathcal{O}_K \to \mathcal{O}_{\hat{K}_\infty}$ and the second appears in (4.3.4).

**Theorem 4.3.8.** *Let $K_\infty/K$ be a ramified $\mathbb{Z}_p$-extension, and let $N = N_{K_\infty/K}$ be large enough so that $K_\infty/K_N$ is sufficiently ramified. Define $\mathcal{Z}$ by the exact triangle*

$$\mathcal{O}_K^{h\,\mathrm{Gal}(K_\infty/K)} \to \mathcal{O}_C^{h\Gamma_K} \to \mathcal{Z},$$

*where the first morphism is the composition in* (4.3.7)*. Then:*

$$H^i(\mathcal{Z}) = \begin{cases} 0, & i = 0, \\ p^{N+3}\text{-torsion (resp., } p^{N+4}\text{-torsion)}, & i = 1, \\ p\text{-torsion}, & i \geq 2 \end{cases}$$

*as $p$ is odd or even, respectively.*

*If $L/K$ is a tame Galois extension and $L_\infty = K_\infty L$, then $L_\infty/L$ is a ramified $\mathbb{Z}_p$-extension, and we can take $N_{L_\infty/L} = N_{K_\infty/K}$.*

*Proof.* Consider the three exact triangles:

$$\begin{aligned} \mathcal{O}_K^{h\,\mathrm{Gal}(K_\infty/K)} &\to \mathcal{O}_{\hat{K}_\infty}^{h\,\mathrm{Gal}(K_\infty/K)} &&\to \mathcal{X} \\ \mathcal{O}_{\hat{K}_\infty}^{h\,\mathrm{Gal}(K_\infty/K)} &\to \mathcal{O}_C^{h\Gamma_K} &&\to \mathcal{Y} \\ \mathcal{O}_K^{h\,\mathrm{Gal}(K_\infty/K)} &\to \mathcal{O}_C^{h\Gamma_K} &&\to \mathcal{Z} \end{aligned}$$

Here the first triangle is from Proposition 4.2.18 (more precisely, the remark following it), the second from (4.3.4), and the third is as in the theorem. By the octahedral axiom, we have an exact triangle $\mathcal{X} \to \mathcal{Z} \to \mathcal{Y}$. The bounds in the theorem now follow by combining Proposition 4.2.18 with Lemma 4.3.5.

The statement about tame extensions follows from Lemma 4.2.6. ☐

Finally, we can prove Theorem 4.0.4, repeated here for convenience.

**Theorem 4.3.9.** *Let $K$ be a local field of characteristic $(0, p)$. Let $C$ be the completion of an algebraic closure $\overline{K}/K$. Let $\Gamma_K = \mathrm{Gal}(\overline{K}/K)$. Then:*

1. *$H^0(\Gamma_K, \mathcal{O}_C) = \mathcal{O}_K$.*
2. *There exists an isomorphism of $\mathcal{O}_K$-modules*
$$H^1_{\mathrm{cts}}(\Gamma_K, \mathcal{O}_C) \cong \mathcal{O}_K \oplus T$$
*where $T$ is $p^4$-torsion (resp., $p^5$-torsion) as $p$ is odd or even, respectively. If $K$ is tame, or if it is a cyclotomic extension of a tame field, this bound can be improved to $p^3$ (resp., $p^4$).*
3. *For $i > 1$, $H^i_{\mathrm{cts}}(\Gamma_K, \mathcal{O}_C)$ is $p$-torsion (resp., $p^2$-torsion) as $p$ is odd or even, respectively.*

*Proof.* Let $K_\infty/K$ be a ramified $\mathbb{Z}_p$-extension. Taking the long exact sequence associated to the exact triangle appearing in Theorem 4.3.8, we obtain $H^0(\Gamma_K, \mathcal{O}_C) = \mathcal{O}_K$, an exact sequence

$$0 \to \mathcal{O}_K \to H^1_{\mathrm{cts}}(\Gamma_K, \mathcal{O}_C) \to H^1(\mathcal{Z}) \to 0,$$

and isomorphisms $H^i_{\mathrm{cts}}(\Gamma_K, \mathcal{O}_C) \cong H^i(\mathcal{Z})$ for $i \geq 2$. The cases $i = 0$ and $i \geq 2$ of the theorem follow directly from Theorem 4.3.8. When $K_\infty/K$ can be taken to be sufficiently ramified (which includes the cases listed in part (2) of the theorem), then $H^1(\mathcal{Z})$ is $p^3$- or $p^4$-torsion as $p$ is odd or even, and then Lemma 4.3.10 below applies to show that $H^1_{\mathrm{cts}}(\Gamma_K, \mathcal{O}_C)$ has the structure claimed by the theorem.

In the general case, let $N$ be large enough so that $K_\infty/K_N$ is sufficiently ramified. Let $G = \mathrm{Gal}(K_N/K)$. The inflation-restriction sequence for the tower $K_\infty/K_n/K$ reads in part:

$$H^1(G, \mathcal{O}_{K_N}) \to H^1_{\mathrm{cts}}(\Gamma_K, \mathcal{O}_C) \to H^1_{\mathrm{cts}}(\Gamma_{K_N}, \mathcal{O}_C)^G$$

(For the term on the left, we have used Ax's theorem to identify $H^0(\Gamma_{K_N}, \mathcal{O}_C)$ with $\mathcal{O}_{K_N}$.) By Sen's theorem (Theorem 4.0.2), the term on the left is $p$-torsion. Since $K_\infty/K_N$ is sufficiently ramified, all torsion in $H^1_{\mathrm{cts}}(\Gamma_{K_N}, \mathcal{O}_C)$ is killed by $p^3$ (resp., $p^4$). Therefore $p^4$ (resp., $p^5$) kills all torsion in $H^1_{\mathrm{cts}}(\Gamma_K, \mathcal{O}_C)$. ☐

**Lemma 4.3.10.** *Given an exact sequence of $\mathcal{O}_K$-modules*

$$0 \to M_1 \to M \to M_2 \to 0$$

*with $M_1$ finite free of rank $r$ and $M_2$ killed by $p^n$, there exists an isomorphism $M \cong \mathcal{O}_K^{\oplus r} \oplus T$, where $T$ is $p^n$-torsion.*

*Proof.* Let $T \subset M$ be the submodule consisting of torsion elements, and let $N = M/T$. We may identify $M_1$ and $N$ with $\mathcal{O}_K$-submodules of $M \otimes_{\mathcal{O}_K} K \cong K^r$ in such a way that $M_1 \subset N \subset p^{-n} M_1$. Since $\mathcal{O}_K$ is Noetherian and $M_1$ is finitely generated, $N$ is finitely generated as well, and therefore (since it is torsion-free and spans $K^r$) it must be free of rank $r$. Let $M_0 \subset M$ be the free submodule spanned by lifts of generators of $N$. Then $M = M_0 \oplus T$. □

4.4. **Galois cohomology of characters.** Once again suppose $K$ is a local field of characteristic $(0, p)$. Let

$$\chi \colon \operatorname{Gal}(\overline{K}/K) \to \mathbb{Z}_p^\times$$

be a character; i.e., a continuous homomorphism. We assume that $\chi$ is infinitely ramified, in the sense that the image of the inertia group under $\chi$ is infinite (equivalently, this image has finite index). Let $K_\infty$ be the fixed field of the kernel of $\chi$. Let $\mathbb{Z}_p(\chi)$ be $\mathbb{Z}_p$ with an action of $\operatorname{Gal}(\overline{K}/K)$ through $\chi$. Then if $M$ is any $p$-adically complete $\operatorname{Gal}(\overline{K}/K)$-module, we may define a new such module by $M(\chi) = M \otimes_{\mathbb{Z}_p} \mathbb{Z}_p(\chi)$. The present goal is to bound the continuous cohomology of $\operatorname{Gal}(K_\infty/K)$ acting on $\mathcal{O}_{\hat{K}_\infty}(\chi)$.

Let $U \subset \mathbb{Z}_p^\times$ be the largest subgroup which is isomorphic to $\mathbb{Z}_p$. Thus $U = 1 + p\mathbb{Z}_p$ for $p$ odd and $U = 1 + p^2\mathbb{Z}_p$ for $p = 2$. Let $K_0$ be the fixed field of $\chi^{-1}(U)$; then $K_\infty/K_0$ is a ramified $\mathbb{Z}_p$-extension. As usual, we let $K_n/K_0$ be the fixed field of $p^n \operatorname{Gal}(K_\infty/K_0)$. Finally, we define an integer $r \geq 1$ by

$$\chi(\operatorname{Gal}(K_\infty/K_0)) = 1 + p^r \mathbb{Z}_p.$$

**Lemma 4.4.1.** *Let $\chi$ be a non-trivial character and $N \geq 0$ be large enough so that $K_\infty/K_N$ is sufficiently ramified (see Lemma 4.2.6). Then*

$$H^i_{\mathrm{cts}}(\operatorname{Gal}(K_\infty/K_0), \mathcal{O}_{\hat{K}_\infty}(\chi)) = \begin{cases} 0, & i = 0 \\ (p^k\text{-torsion}), & i = 1. \end{cases}$$

*Here we may take $k = N + r + 1$ if $p$ is odd and $k = N + r + 2$ if $p = 2$.*

*Proof.* We first prove the result in the case that $K_\infty/K_0$ is sufficiently ramified. Write $t \colon \hat{K}_\infty \to K_0$ for the normalized trace. By Lemma 4.2.17, the inclusion $\mathcal{O}_{K_0} \oplus \mathcal{O}^{t=0}_{\hat{K}_\infty} \hookrightarrow \mathcal{O}_{\hat{K}_\infty}$ has $\alpha$-torsion cokernel for any $\alpha \in \mathcal{O}_{K_N}$ with $v(\alpha) \geq 1/(p-1)$.

We are therefore reduced to studying the continuous cohomology of $\operatorname{Gal}(K_\infty/K_0)$ acting on $\mathcal{O}_{K_0}(\chi)$ and $\mathcal{O}^{t=0}_{\hat{K}_\infty}(\chi)$, respectively. Let $\sigma \in \operatorname{Gal}(K_\infty/K_0)$ be a topological generator, and let $\lambda = \chi(\sigma)^{-1}$. Then for any $p$-adically complete $\operatorname{Gal}(K_\infty/K_0)$-module $M$, the continuous cohomology of $M(\chi)$ is computed by the complex $\sigma - \lambda \colon M \to M$. Note that $\lambda = 1 + p^r u$ for a $p$-adic unit $u \in \mathbb{Z}_p^\times$.

In the case $M = \mathcal{O}_{K_0}$, we have $H^0(\operatorname{Gal}(K_\infty/K_0), \mathcal{O}_{K_0}(\chi)) = 0$, since $\sigma$ acts on $\mathcal{O}_{K_0}(\chi)$ as the scalar $\lambda \neq 1$. On the other hand $H^1_{\mathrm{cts}}(\operatorname{Gal}(K_\infty/K_0), \mathcal{O}_{K_0}(\chi))$ is $p^r$-torsion, being the cokernel of multiplication by $\lambda - 1 = p^r u$.

To treat $M = \mathcal{O}^{t=0}_{\hat{K}_\infty}$, we apply Lemma 4.2.12 to the sufficiently ramified extension $K_\infty/K_0$. We find that $(\sigma - 1)^{-1}$ is defined on $\hat{K}^{t=0}_\infty$ and has operator norm $\leq |p|^{-1-1/p(p-1)}$.

Let

$$\mu = (\sigma - 1)^{-1}(\sigma - \lambda) = 1 - (\lambda - 1)(\sigma - 1)^{-1}.$$

In the case that $r \geq 2$, we have $|\lambda - 1|\,|\sigma - 1|^{-1} < 1$, and so $\mu$ has continuous inverse satisfying $\left|\mu^{-1}\right| \leq 1$. Therefore $(\sigma - \lambda)^{-1} = (\sigma - 1)^{-1}\mu^{-1}$ exists on $\hat{K}_\infty^{t=0}$ and has operator norm bounded by $|p|^{-1-1/p(p-1)}$. Thus $H^0(\mathrm{Gal}(K_\infty/K_0), \mathcal{O}^{t=0}_{\hat{K}_\infty}(\chi)) = 0$ and $H^1_{\mathrm{cts}}(\mathrm{Gal}(K_\infty/K_0), \mathcal{O}^{t=0}_{\hat{K}_\infty}(\chi))$ is annihilated by any element $\alpha \in \mathcal{O}_{K_0}$ with $v(\alpha) \geq 1 + 1/p(p-1)$.

If $r = 1$, the idea is to apply the same argument to the sufficiently ramified extension $K_\infty/K_1$, noting that $\lambda^p = \chi(\sigma^p)$ now satisfies $|\lambda^p - 1|\,|\sigma^p - 1|^{-1} < 1$. Thus $H^i_{\mathrm{cts}}(\mathrm{Gal}(K_\infty/K_1), \mathcal{O}^{t=0}_{\hat{K}_\infty}(\chi))$ is 0 for $i = 0$ and is $p^{1+1/p(p-1)}$-torsion for $i = 1$. The inflation-restriction sequence allows us to deduce the same results for $K_\infty/K_0$.

We have found that $H^0(\mathrm{Gal}(K_\infty/K_0), \mathcal{O}_{\hat{K}_\infty}(\chi)) = 0$ and $H^1_{\mathrm{cts}}(\mathrm{Gal}(K_\infty/K_0), \mathcal{O}_{\hat{K}_\infty}(\chi))$ is $\alpha$-torsion for any $\alpha \in \mathcal{O}_{K_0}$ with $v(\alpha) \geq \max\{r, 1 + 1/p(p-1)\} + 1/(p-1)$. The latter is $\geq r+1$ or $\geq r+2$ as $p$ is odd or even, respectively, which implies the lemma in the case that $K_\infty/K_0$ is sufficiently ramified. (To handle the problem that $\mathcal{O}_{K_0}$ may not contain elements of sufficiently small valuation, it may be necessary to pass to a tamely ramified extension $L/K_0$, establish the bound there, and then descend along $L/K_0$, using Noether's theorem to show that this descent does not introduce additional cohomology.)

In general, suppose $N \geq 0$ is large enough so that $K_\infty/K_N$ is sufficientlty ramified. Note that $\mathrm{Gal}(K_\infty/K_N)$ is generated by $\sigma^{p^N}$ and that $\chi(\sigma^{p^N}) \equiv 1 \pmod{p^{N+r}}$. Thus the preceding argument shows that $H^0_{\mathrm{cts}}(\mathrm{Gal}(K_\infty/K_N), \mathcal{O}_{\hat{K}_\infty}(\chi)) = 0$ and $H^1_{\mathrm{cts}}(\mathrm{Gal}(K_\infty/K_N), \mathcal{O}_{\hat{K}_\infty}(\chi))$ is $p^{N+r+1}$- or $p^{N+r+2}$-torsion as $p$ is odd or even, respectively. The inflation-restriction sequence for the tower $K_\infty/K_N/K_0$ reads:

$$\begin{aligned} H^1(\mathrm{Gal}(K_N/K_0), \mathcal{O}_{\hat{K}_\infty}(\chi)^{\mathrm{Gal}(K_\infty/K_N)}) &\to H^1_{\mathrm{cts}}(\mathrm{Gal}(K_\infty/K_0), \mathcal{O}_{\hat{K}_\infty}(\chi)^{\mathrm{Gal}(K_\infty/K_N)}) \\ &\to H^1_{\mathrm{cts}}(\mathrm{Gal}(K_\infty/K_N), \mathcal{O}_{\hat{K}_\infty})^{\mathrm{Gal}(K_N/K_0)} \end{aligned}$$

By the above argument, the term on the left is zero, and we get the result. □

**Lemma 4.4.2.** *Assume that $K_\infty/K_N$ is sufficiently ramified. Let $a_p = 0$ if $p$ is odd, and $a_p = 1$ if $p = 2$. We have:*

$$H^i_{\mathrm{cts}}(\mathrm{Gal}(K_\infty/K), \mathcal{O}_{\hat{K}_\infty}(\chi)) = \begin{cases} 0, & i = 0, \\ (p^{N+r+1+a_p}\text{-torsion}), & i = 1, \\ (p^{a_p}\text{-torsion}), & i > 1. \end{cases}$$

*Proof.* Combine the spectral sequence

$$H^i(\mathrm{Gal}(K_0/K), H^j_{\mathrm{cts}}(\mathrm{Gal}(K_\infty/K_0), \mathcal{O}_{\hat{K}_\infty}(\chi))) \implies H^{i+j}_{\mathrm{cts}}(\mathrm{Gal}(K_\infty/K), \mathcal{O}_{\hat{K}_\infty}(\chi))$$

with Lemma 4.4.1. The only terms on the left that contribute occur when $j = 1$. If $p$ is odd, the $i > 0$ terms on the left side vanish because $\#\,\mathrm{Gal}(K_0/K)$ is invertible in $\mathcal{O}_K$. If $p = 2$ then $\#\,\mathrm{Gal}(K_0/K) = 2$, and the best we can say is that the $i > 0$ terms on the left side are 2-torsion. □

At this point we specialize to the case of the cyclotomic character

$$\chi_{\mathrm{cycl}} \colon \Gamma_K \to \mathbb{Z}_p^\times,$$

defined by the relation $\tau(\zeta) = \zeta^{\chi_{\mathrm{cycl}}(\tau)}$ for any $\tau \in \Gamma_K$ and any $p$th power root of unity $\zeta$. For any $j \in \mathbb{Z}$ we write $\mathbb{Z}_p(j) = \mathbb{Z}_p(\chi^j_{\mathrm{cycl}})$ for the $j$th Tate twist, and similarly $\mathcal{O}_C(j) = \mathcal{O}_C \otimes \mathbb{Z}_p(j)$.

We conclude the section with a bound on the cohomology of $\mathcal{O}_C(j)$.

**Theorem 4.4.3.** *With notation as in Theorem 4.0.4, let $j \neq 0$ be an integer, and let $\mathcal{O}_C(j)$ be the $j$th Tate twist of $\mathcal{O}_C$ as a $\mathrm{Gal}(\overline{K}/K)$-module. Then:*

*(1)* $H^0(\Gamma_K, \mathcal{O}_C(j)) = 0$.

*(2)* $H^1_{\mathrm{cts}}(\Gamma_K, \mathcal{O}_C(j))$ *is* $p^{M+v(j)}$*-torsion. Here* $M = M_K$ *is a constant which only depends on* $K$ *and which is insensitive to passage to a tamely ramified extension of* $K$*. If* $p \nmid e_K$ *we may take* $M = 2$ *if* $p$ *is odd and* $M = 4$ *if* $p = 2$*.*

*(3) For* $i \geq 2$*,* $H^i_{\mathrm{cts}}(\Gamma_K, \mathcal{O}_C(j))$ *is* $p$*-torsion if* $p$ *is odd and* $p^2$*-torsion if* $p = 2$*.*

*(We recall that* $v(j)$ *the* $p$*-adic valuation of* $j$*, normalized so that* $v(p) = 1$*.)*

*Proof.* Let $K_\infty/K$ be the extension obtained by adjoining all $p$th power roots of unity. As in the beginning of this subsection, we let $K_0/K$ be the minimal subextension such that $K_\infty/K_0$ is a $\mathbb{Z}_p$-extension, and we define an integer $r \geq 1$ by $\chi_{\mathrm{cycl}}(\mathrm{Gal}(K_\infty/K_0)) = 1 + p^r\mathbb{Z}_p$. Then

$$\chi^j_{\mathrm{cycl}}(\mathrm{Gal}(K_\infty/K_0)) = (1 + p^r\mathbb{Z}_p)^j = 1 + p^{r+v(j)}\mathbb{Z}_p.$$

Lemma 4.4.2 bounds the cohomology $H^i_{\mathrm{cts}}(\mathrm{Gal}(K_\infty/K), \mathcal{O}_{\hat{K}_\infty}(j))$ in terms of the integer $N = N_K$ for which $K_\infty/K_N$ is infinitely ramified: we have $H^0 = 0$, $H^1$ is $p^{N+r+v(j)+a_p}$-torsion, and $H^2$ is $p^{a_p}$-torsion, where $a_p$ is 0 or 1 as $p$ is odd or even.

As in Lemma 4.3.5, define a complex $\mathcal{Y}_0$ by the exact triangle

$$\mathcal{O}_{\hat{K}_\infty} \to \mathcal{O}_C^{h\,\mathrm{Gal}(\overline{K}/K_\infty)} \to \mathcal{Y}_0.$$

Then $H^0(\mathcal{Y}_0) = 0$ and $H^i(\mathcal{Y}_0)$ is almost zero for $i > 0$. Twisting by $\chi^j_{\mathrm{cycl}}$ and taking derived $\mathrm{Gal}(K_\infty/K)$-invariants, we obtain an exact triangle

$$\mathcal{O}_{\hat{K}_\infty}(j)^{h\,\mathrm{Gal}(K_\infty/K)} \to \mathcal{O}_C(j)^{h\Gamma_K} \to \mathcal{Y},$$

where $\mathcal{Y} = \mathcal{Y}_0(j)^{h\,\mathrm{Gal}(K_\infty/K)}$. The same method of proof for Lemma 4.3.5 shows that $H^0(\mathcal{Y}) = 0$ and $H^i(\mathcal{Y})$ is $p$-torsion for $i > 0$. Combining these bounds with those obtained in Lemma 4.4.2, we obtain the bounds appearing in the theorem. $\square$

## 5. Pro-étale cohomology of rigid-analytic spaces

Let $K$ be a local field of characteristic $(0, p)$, and let $X$ be a smooth rigid-analytic space over $K$. We consider three topologies on $X$ in order of coarsest to finest:

(1) The analytic topology $X_{\mathrm{an}}$, whose opens are simply the open subsets of the underlying topological space of $X$,
(2) The étale topology $X_{\mathrm{\acute{e}t}}$, whose opens are étale morphisms $U \to X$,
(3) The pro-étale topology $X_{\mathrm{pro\acute{e}t}}$, whose opens are formal limits $U = \varprojlim U_i$, where each $U_i \to X$ is étale.

In this section, we present some comparison results relating the cohomology of the structure sheaf among these sites. For the purposes of exposition, we start with the case of the rational structure sheaf $\mathcal{O}_X$. Recall from Section 3 the notation $\hat{\mathcal{O}}_X$ for the completed structure sheaf on $X_{\mathrm{pro\acute{e}t}}$.

The case that $X = \mathrm{Spa}\,K$ is a single point is instructive. In this case $H^*(X_{\mathrm{an}}, \mathcal{O}_X) = K$ rather trivially, and $H^*(X_{\mathrm{\acute{e}t}}, \mathcal{O}_X)$ is the Galois cohomology $H^*_{\mathrm{cts}}(\mathrm{Gal}(\overline{K}/K), \overline{K})$, where one must put the discrete topology on $\overline{K}$; by the normal basis theorem this is $K$ again. By contrast, $H^*(X_{\mathrm{pro\acute{e}t}}, \hat{\mathcal{O}}_X)$ may be identified with $H^*_{\mathrm{cts}}(\mathrm{Gal}(\overline{K}/K), C)$, where $C$ is the $p$-adic completion of $\overline{K}$ (with its $p$-adic topology). Tate's theorem (Theorem 4.0.3) states that this is isomorphic to $K[\varepsilon]$.

In the general case where $X$ is a smooth rigid-analytic variety over $K$, it is still the case that $H^*(X_{\mathrm{an}}, \mathcal{O}_X)$ is isomorphic to $H^*(X_{\mathrm{\acute{e}t}}, \mathcal{O}_X)$ by étale descent, but of course $H^*(X_{\mathrm{pro\acute{e}t}}, \hat{\mathcal{O}}_X)$ may be different. The general pattern we observed is that the difference is accounted for entirely by the pro-étale cohomology of the base $\mathrm{Spa}\,K$. In the next subsection we make this explicit for the rational structure sheaf before passing to the study of the integral structure sheaf in subsequent subsections.

5.1. **The rational comparison isomorphism.** The following is a rapid discussion of the sort of comparison result we need, in the much simpler situation where $p$ has been inverted. We leave precise details to the subsections that follow.

Let $C$ be an algebraically closed nonarchimedean field of characteristic $(0, p)$, and let $X$ be a smooth rigid-analytic space over $C$. We recall here a theorem of Scholze, which relates the difference between $X_{\text{ét}}$ and $X_{\text{proét}}$ to differential forms on $X$. Let us write $\Omega^j_{X/C}$ for the sheaf of differential $j$-forms, considered as a sheaf of $\mathcal{O}_X$-modules on $X_{\text{ét}}$. Also we write

$$\nu\colon X_{\text{proét}} \to X_{\text{ét}}$$

for the projection map between sites.

**Theorem 5.1.1** ([Sch13b, Proposition 3.23])**.** *Let $X$ be a smooth rigid-analytic space over $C$. Let $\nu\colon X_{\text{proét}} \to X_{\text{ét}}$ be the projection. Then for each $j \geq 0$ there is an isomorphism of $\mathcal{O}_{X_{\text{ét}}}$-modules:*

$$\Omega^j_{X_{\text{ét}}/C}(-j) \cong R^j\nu_*\hat{\mathcal{O}}_X.$$

The essential calculation behind Theorem 5.1.1 goes back to Faltings [Fal88]. The assumption that $X$ is smooth means that (locally on $X_{\text{ét}}$) it admits an étale morphism to the $d$-dimensional torus

$$\begin{aligned} \mathbb{T}^d &= \operatorname{Spa}(R_d, R_d^+) \\ R_d^+ &= \mathcal{O}_C\left\langle T_1^{\pm}, \dots, T_d^{\pm 1}\right\rangle \\ R_d &= R_d^+[1/p]. \end{aligned}$$

Therefore let us explain the proof of Theorem 5.1.1 in the case of $X = \mathbb{T}^d$. There is an affinoid perfectoid pro-étale torsor $\tilde{\mathbb{T}}^d \to \mathbb{T}^d$ for the group $\mathbb{Z}_p(1)^d$, namely

$$\begin{aligned} \tilde{\mathbb{T}}^d &= \operatorname{Spa}(\tilde{R}_d, \tilde{R}_d^+) \\ \tilde{R}_d^+ &= \mathcal{O}_C\left\langle T_1^{\pm 1/p^\infty}, \dots, T_d^{\pm 1/p^\infty}\right\rangle \\ \tilde{R}_d &= \tilde{R}_d^+[1/p]. \end{aligned}$$

(Here $\mathbb{Z}_p(1) = \varprojlim \mu_{p^n}$, where $\mu_{p^n}$ is the group of $p^n$th roots of 1.) Affinoid perfectoid covers in $X_{\text{proét}}$ are convenient because if $U$ is affinoid perfectoid, then $H^i(U_{\text{proét}}, \hat{\mathcal{O}}_U) = 0$ is for $i > 0$ [Sch13a, Lemma 4.10]. As in the proof of Proposition 3.6.3, this allows us to compute the pro-étale cohomology of $X$ in terms of group cohomology:

$$H^*(\mathbb{T}^d_{\text{proét}}, \hat{\mathcal{O}}_{\mathbb{T}^d}) \cong H^*_{\text{cts}}(\mathbb{Z}_p(1)^d, \tilde{R}_d)$$

But now an explicit calculation in group cohomology shows that the natural map

$$H^*_{\text{cts}}(\mathbb{Z}_p(1)^d, R_d) \to H^*_{\text{cts}}(\mathbb{Z}_p(1)^d, \tilde{R}_d)$$

is an isomorphism. Since $\mathbb{Z}_p(1)^d$ is profinite free abelian of rank $d$ and acts trivially on $R_d$, we have $H^i_{\text{cts}}(\mathbb{Z}_p(1)^d, R_d) \cong \bigwedge^i_{R_d}(R_d)^d$ (see [Sch13a, Lemma 5.5]).

The upshot is that $R^i\nu_*\hat{\mathcal{O}}_X \cong \bigwedge^i R^1\nu_*\hat{\mathcal{O}}_X$, and $R^1\nu_*\hat{\mathcal{O}}_X$ is a locally free $\mathcal{O}_{X_{\text{ét}}}$-module of rank $d$. This already suggests that $R^j\nu_*\hat{\mathcal{O}}_X$ should be related to differentials. We refer the reader to [Sch13b, Lemma 3.24] for a functorial construction of the isomorphism in Theorem 5.1.1.

Now suppose once again that $X$ is a rigid-analytic space over a local field $K$. Write $C$ for the completion of an algebraic closure $\overline{K}$, and write $\Gamma_K = \operatorname{Gal}(\overline{K}/K)$. Let $\overline{X}$ be the base change of $X$ to $C$. Once again we write $\nu\colon X_{\text{proét}} \to X_{\text{ét}}$ for the projection; let $\overline{\nu}\colon \overline{X}_{\text{proét}} \to \overline{X}_{\text{ét}}$ be the corresponding projection for $\overline{X}$. We have

$$R\nu_*\hat{\mathcal{O}}_X = R\nu_*(\hat{\mathcal{O}}_{\overline{X}}^{h\Gamma_K}) = (R(\overline{\nu})_*\hat{\mathcal{O}}_{\overline{X}})^{h\Gamma_K}.$$

Applying Theorem 5.1.1, we find that $(R(\overline{\nu})_*\hat{\mathcal{O}}_{\overline{X}})^{h\Gamma_K}$ admits a filtration with graded pieces $\Omega^j_{\overline{X}/C}(-j)^{h\Gamma_K}$. For an affinoid object $U \in X_{\text{ét}}$, the global sections of the $j$th piece on $U$ are

$$\begin{aligned}\Gamma(U, \Omega^j_{\overline{X}/C}(-j)^{h\Gamma_K}) &\cong (\Gamma(U, \Omega^j_{U/K})\hat{\otimes}_K C(-j))^{h\Gamma_K} \\ &\cong \Gamma(U, \Omega^j_{U/K})\hat{\otimes}_K C(-j)^{h\Gamma_K}.\end{aligned}$$

(Here the $\hat{\otimes}$ means completed tensor product of Banach $K$-vector spaces; to commute the $h\Gamma_K$ into this operation, we have used the projection formula, Lemma 3.4.10.) By Theorem 4.0.3, the terms with nonzero $j$ vanish, and the $j = 0$ term is isomorphic to $\Gamma(U, \mathcal{O}_{X_{\text{ét}}})[\varepsilon]$. (We remind the reader that for an abelian group $Y$, $Y[\varepsilon]$ is the complex $Y \oplus Y[1]$.) Therefore

$$R\nu_*\hat{\mathcal{O}}_X \cong \mathcal{O}_{X_{\text{ét}}}[\varepsilon].$$

We have proved:

**Theorem 5.1.2.** *Let $K$ be a discretely valued nonarchimedean field of characteristic 0 with perfect residue field. Let $X/K$ be a smooth rigid-analytic space. There is an isomorphism*

$$R\Gamma(X_{\text{ét}}, \mathcal{O}_X)[\varepsilon] \cong R\Gamma(X_{\text{proét}}, \hat{\mathcal{O}}_X).$$

5.2. **Integral $p$-adic Hodge theory.** Suppose again that $C$ is an algebraically closed nonarchimedean field of characteristic $(0,p)$, and that $X$ is a smooth rigid-analytic space over $C$. When $X$ has a sufficiently nice formal model over $\mathcal{O}_C$, the theorems of [BMS18] and [ČK19] can be used to gain control over the integral pro-étale cohomology $R\Gamma(X_{\text{proét}}, \hat{\mathcal{O}}^+_X)$.

**Definition 5.2.1** ([ČK19])**.** Let $\mathfrak{X}$ be a formal scheme over $\operatorname{Spf}\mathcal{O}_C$. We say that $\mathfrak{X}$ is *semistable of dimension $d$* if $\mathfrak{X}$ can be covered by affine opens $\mathfrak{U}$ which admit an étale $\mathcal{O}_C$-morphism to a formal scheme of the form

$$\operatorname{Spf}\mathcal{O}_C\left\langle T_0, \dots, T_r, T^{\pm}_{r+1}, \dots, T^{\pm}_d\right\rangle/(T_0\cdots T_r - \pi) \qquad (5.2.2)$$

where $0 \le r \le d$, and $\pi \in \mathcal{O}_C$ a nonzero element, such that $|\pi|$ is a rational power of $|p|$. (The values of $r$ and $\pi$ may vary with $\mathfrak{U}$.)

By the discussion in [ČK19], $\mathfrak{X}$ is a log-smooth formal scheme over $\operatorname{Spf}\mathcal{O}_C$, and we have the sheaf of continuous log-differentials $\Omega^1_{\mathfrak{X},\log}$, a locally free $\mathcal{O}_{\mathfrak{X}}$-module of rank $d$. Finally, $\Omega^j_{\mathfrak{X},\log} = \bigwedge^j \Omega^1_{\mathfrak{X},\log}$.

There is an evident extension of Definition 5.2.1 to the setting of formal schemes over $\operatorname{Spf}\mathcal{O}_K$, where $K$ is a local field. If $\mathfrak{X}$ is a semistable formal scheme over $\operatorname{Spf}\mathcal{O}_K$, then its adic generic fiber $X = \mathfrak{X}_K$ is a smooth rigid-analytic variety.

**Lemma 5.2.3.** *Let $\mathfrak{X}$ be an affine semistable formal scheme over* $\operatorname{Spf}\mathcal{O}_K$*, with generic fiber $X$. The natural map $\mathcal{O}(\mathfrak{X}) \to \mathcal{O}^+(X)$ is an isomorphism.*

*Proof.* Let $\pi$ be a uniformizer for $K$. Following the recipe for calculating the generic fiber of $\mathfrak{X} \to \operatorname{Spa}(\mathcal{O}_K, \mathcal{O}_K)$ (see [SW20, Proposition 5.1.5], noting that this map is "adic"), we find that $\mathcal{O}^+(X)$ is the integral closure of $\mathcal{O}(\mathfrak{X})$ in $\mathcal{O}(\mathfrak{X})[1/\pi]$.

We claim that $\mathcal{O}(\mathfrak{X})$ is already integrally closed in $\mathcal{O}(\mathfrak{X})[1/\pi]$. All the same, we claim $\mathcal{O}(\mathfrak{X})/\pi$ is reduced. The claim is étale local on $\operatorname{Spec}\mathcal{O}(\mathfrak{X})/\pi$. Since $\mathfrak{X}$ is semistable, $\operatorname{Spec}\mathcal{O}(\mathfrak{X})/\pi$ is étale locally isomorphic to

$$\operatorname{Spec} k[T_0, \dots, T_r, T^{\pm 1}_{r+1}, \dots, T_r],$$

which is reduced. □

If $\mathfrak{U} \in \mathfrak{X}_{\text{ét}}$, then the adic generic fiber of $\mathfrak{U}$ is an object of $X_{\text{proét}}$, which is to say there is a morphism of sites:

$$\nu\colon X_{\text{proét}} \to \mathfrak{X}_{\text{ét}}.$$

The integral version of Theorem 5.1.1 relates the difference between $X_{\text{proét}}$ and $\mathfrak{X}_{\text{ét}}$ in terms of the sheaf of differentials $\Omega^j_{\mathfrak{X},\log}$. To make this work there are two additional ingredients: the Breuil–Kisin twist and the décalage functor.

The Breuil–Kisin twist $\mathcal{O}_C\{1\}$ is a free $\mathcal{O}_C$-module of rank 1 carrying a $\operatorname{Gal}(\overline{K}/K)$-action. See [BMS18, Definition 8.2] for its precise definition. There is a canonical Galois-equivariant injection $\mathcal{O}_C(1) \hookrightarrow \mathcal{O}_C\{1\}$ whose cokernel is killed by $(\zeta_p - 1)$, where $\zeta_p$ is a primitive $p$th root of 1. For $j \in \mathbb{Z}$ we let $\mathcal{O}_C\{j\}$ be the $j$th tensor power of $\mathcal{O}_C\{1\}$.

The décalage functor will be reviewed in Section 5.5. In our context it appears as an endofunctor $L\eta_{(\zeta_p-1)}$ on the derived category of $\mathcal{O}_{\mathfrak{X}_{\mathcal{O}_C}}$-modules. For the moment we need two facts concerning $L\eta_{(\zeta_p-1)}$: it is a lax symmetric monoidal functor, and also there is a natural transformation

$$a\colon L\eta_{(\zeta_p-1)}\mathcal{C} \to \mathcal{C} \tag{5.2.4}$$

whenever $\mathcal{C}$ is bounded below by 0 and $H^0(\mathcal{C})$ is torsion-free. Let

$$\widetilde{\Omega}_{\mathfrak{X}} = L\eta_{(\zeta_p-1)} R\nu_* \hat{\mathcal{O}}^+_X,$$

so that $\widetilde{\Omega}_{\mathfrak{X}}$ is a complex of $\mathcal{O}_{\mathfrak{X}}$-modules on $\mathfrak{X}_{\text{ét}}$.

The following theorem is the basis for all the comparison isomorphisms that follow.

**Theorem 5.2.5** ([BMS18, Theorem 8.3],[ČK19, Theorem 4.11])**.** *Let $\mathfrak{X}$ be a semistable affine formal scheme over $\mathcal{O}_C$. For each $j \geq 0$ there is a canonical isomorphism of sheaves of $\mathcal{O}_{\mathfrak{X}}$-modules on $\mathfrak{X}_{\text{ét}}$:*

$$\Omega^j_{\mathfrak{X},\log}\{-j\} \cong H^j(\widetilde{\Omega}_{\mathfrak{X}}).$$

Now let $K$ be a local subfield of $C$. Assume that $\mathfrak{X}$ is a semistable affine formal scheme over $\mathcal{O}_K$, in the sense that it has an étale chart made of formal schemes of the form (5.2.2), except that $\mathcal{O}_K$ replaces $\mathcal{O}_C$. Denote by $\overline{\mathfrak{X}}$ the base change to $\mathcal{O}_C$. We have

$$\Omega^j_{\overline{\mathfrak{X}},\log} \cong \Omega^j_{\mathfrak{X},\log} \hat{\otimes}_{\mathcal{O}_K} \mathcal{O}_C$$

(where $\hat{\otimes}$ means $p$-adically completed tensor product). The theorem above can be applied to $\mathfrak{X}_{\mathcal{O}_C}$. The group $\Gamma_K$ acts on $\Omega^j_{\mathfrak{X}_{\mathcal{O}_C}}$ and $\tilde{\Omega}$. (Note that $\Gamma_K$ preserves the ideal $(\zeta_p - 1)$ and so its action commutes with $L\eta_{(\zeta_p-1)}$.)

We get the following corollary of the theorem:

**Corollary 5.2.6.** *Let $\mathfrak{X}$ is a semistable affine formal scheme over $\mathcal{O}_K$. For each $j \geq 0$, there is a canonical $\Gamma_K$-equivariant isomorphism*

$$\Omega^j_{\mathfrak{X},\log} \hat{\otimes}_{\mathcal{O}_K} \mathcal{O}_C\{-j\} \cong H^j(\widetilde{\Omega}_{\overline{\mathfrak{X}}}).$$

5.3. **The integral comparison isomorphism for affine semistable formal schemes.** The idea now is to present integral versions of the comparison isomorphism in Theorem 5.1.2. The first version applies to the setting where $K$ is a local field of characteristic $(0,p)$, and $X$ is a rigid-analytic space over $K$ admitting an *affine* semistable model $\mathfrak{X}$ over $\mathcal{O}_K$. We will write $\mathcal{O}$ for the structure sheaf of $\mathfrak{X}_{\text{ét}}$, and $\hat{\mathcal{O}}^+$ for the completed integral structure sheaf on $X_{\text{proét}}$.

Consider the projection $\nu\colon X_{\text{proét}} \to \mathfrak{X}$. It pulls back functions to integral functions, which is to say we have an $\mathcal{O}_K$-algebra homomorphism $\mathcal{O} \to \nu_* \hat{\mathcal{O}}^+$. Taking derived global sections and noting that $\mathcal{O}$ is acyclic (since $\mathfrak{X}$ is affine), we obtain a homomorphism of $\mathcal{O}_K$-algebra objects in $D(\mathrm{Ab})$:

$$\nu^*\colon \mathcal{O}(\mathfrak{X}) \to R\Gamma(X_{\text{proét}}, \hat{\mathcal{O}}^+)$$

Similarly, the structure morphism $s\colon X \to \operatorname{Spa} K$ induces a pullback map on the level of cohomology (here we confuse $K$ with $\operatorname{Spa} K$):

$$s^*\colon R\Gamma(K_{\text{proét}}, \hat{\mathcal{O}}^+) \to R\Gamma(X_{\text{proét}}, \hat{\mathcal{O}}^+)$$

In Example 3.6.4 we have identified $R\Gamma(K_{\text{proét}}, \hat{\mathcal{O}}^+)$ with the Galois cohomology $R\Gamma_{\text{cts}}(\Gamma_K, \mathcal{O}_C) = R\Gamma_* \mathcal{O}_C^{h\Gamma_K}$, where $\Gamma_K = \operatorname{Gal}(\overline{K}/K)$.

Let $K_\infty/K$ be a ramified $\mathbb{Z}_p$-extension. There is no harm in assuming that $K_\infty/K$ is the cyclotomic $\mathbb{Z}_p$-extension (see Section 4.2). Theorem 4.3.8 states that there is a homomorphism of $\mathcal{O}_K$-algebra objects in $D(\mathrm{Ab})$:

$$\alpha_0\colon \mathcal{O}_K[\varepsilon] \to R\Gamma(K_{\text{proét}}, \hat{\mathcal{O}}^+) \tag{5.3.1}$$

whose cofiber has the bounds:

$$H^i(\operatorname{cof}(\alpha_0)) = \begin{cases} 0, & i = 0, \\ p^{M_K}\text{-torsion}, & i = 1, \\ p\text{-torsion}, & i \geq 2. \end{cases} \tag{5.3.2}$$

(Throughout this discussion, we identify $\operatorname{Gal}(K_\infty/K)$ with $\mathbb{Z}_p$, so that $\mathcal{O}_K^{h\mathbb{Z}_p}$ is identified with $\mathcal{O}_K[\varepsilon]$.)

Altogether we obtain a homomorphism of $\mathcal{O}_K$-algebra objects in $D(\mathrm{Ab})$:

$$\alpha_{\mathfrak{X}}\colon \mathcal{O}(\mathfrak{X}) \otimes_{\mathcal{O}_K} \mathcal{O}_K[\varepsilon] \to R\Gamma(X_{\text{proét}}, \hat{\mathcal{O}}^+) \tag{5.3.3}$$

obtained by sending $x \otimes y$ to $\nu^*(x) \otimes s^*\alpha_0(y)$.

**Theorem 5.3.4.** *Let $K$ be a local field of characteristic $(0, p)$. Let $\mathfrak{X}$ be an affine formal scheme which is semistable over $\mathcal{O}_K$ of dimension $d$. Let $\alpha_{\mathfrak{X}}$ be the map defined in (5.3.3). Then for all $i \in \mathbb{Z}$, $H^i(\operatorname{cof}(\alpha_{\mathfrak{X}}))$ is killed by $p^{M_{K^t}(d)}$, where $M_{K^t}(d)$ is a linear function of $d$ which only depends on the maximal tame extension $K^t$ of $K$. (This means the function is insensitive to replacing $K$ with any finite tame extension.)*

Note that the theorem recovers Theorem 4.3.8 in the case $X = \operatorname{Spa} K$.

*Proof.* We will use repeatedly this fact: if $a_1$ and $a_2$ are composable morphisms in $D(\mathrm{Ab})$ and integers $n_{ij}$ for each $i \in \mathbb{Z}$, $j \in \{1, 2\}$ such that $p^{n_{ij}}$ kills $H^i(\operatorname{cof}(a_j))$ for each $i, j$, then $p^{n_{i1}+n_{i2}}$ kills $H^i(\operatorname{cof}(a_2 \circ a_1))$ for all $i$. As a result, if each of $\operatorname{cof}(a_1), \operatorname{cof}(a_2)$ is bounded from below and subject to the same sort of bounds as appears in the theorem, then so is $\operatorname{cof}(a_2 \circ a_1)$.

The map $\alpha_{\mathfrak{X}}$ factors as

$$\mathcal{O}(\mathfrak{X}) \otimes_{\mathcal{O}_K} \mathcal{O}_K[\varepsilon] \overset{\mathrm{id}\otimes\alpha_0}{\to} \mathcal{O}(\mathfrak{X}) \hat{\otimes}_{\mathcal{O}_K} \mathcal{O}_C^{h\Gamma_K} \overset{\nu^*\otimes s^*}{\to} R\Gamma(X_{\text{proét}}, \hat{\mathcal{O}}^+).$$

Since $\mathcal{O}(\mathfrak{X})$ is a flat $\mathcal{O}_K$-module, the cofiber of $\mathrm{id} \otimes \alpha_0$ is subject to the same bounds as in (5.3.2). In light of the observation we made in the first paragraph, it is enough to show that $\operatorname{cof}(\nu^* \otimes s^*)$ is subject to the same bounds as in the theorem.

Let $C$ be the completion of an algebraic closure of $K$; as above we write $\overline{\mathfrak{X}}$ for the base change of $\mathfrak{X}$ to $\operatorname{Spf} \mathcal{O}_C$. Then $\overline{\mathfrak{X}}$ is an affine semistable formal scheme over $\mathcal{O}_C$ admitting a semilinear action of $\Gamma_K = \operatorname{Gal}(\overline{K}/K)$, whose rigid-analytic generic fiber is $\overline{X} = X \times_{\operatorname{Spa} K} \operatorname{Spa} C$. Let $\overline{\nu}\colon X_{C,\text{proét}} \to \overline{\mathfrak{X}}_{\text{ét}}$ be the projection.

We have a commutative diagram:

$$\begin{array}{ccc}
\mathcal{O}(\mathfrak{X})\hat{\otimes}_{\mathcal{O}_K}\mathcal{O}_C^{h\Gamma_K} & \xrightarrow{\nu^*\otimes s^*} & R\Gamma(X_{\text{proét}},\hat{\mathcal{O}}^+) \\
\cong\big\downarrow & & \big\downarrow\cong \\
\mathcal{O}(\overline{\mathfrak{X}})^{h\Gamma_K} & \longrightarrow & R\Gamma(X_{C,\text{proét}},\hat{\mathcal{O}}^+)^{h\Gamma_K} \\
\cong\big\downarrow & & \big\downarrow\cong \\
R\Gamma(\overline{\mathfrak{X}}_{\text{ét}},\mathcal{O})^{h\Gamma_K} & \longrightarrow & R\Gamma(\overline{\mathfrak{X}}_{\text{ét}},R\overline{\nu}_*\hat{\mathcal{O}}^+)^{h\Gamma_K}
\end{array}$$

The upper left vertical arrow is the projection formula from Lemma 3.4.10, and the upper right vertical arrow is descent along the $\Gamma_K$ pro-étale torsor $\overline{X}\to X$, see Proposition 3.6.3.

In light of this diagram, we wish to bound the cofiber of the lower horizontal morphism, namely $R\Gamma(\overline{\mathfrak{X}}_{\text{ét}},\operatorname{cof}\overline{\nu}^*)^{h\Gamma_K}$, where $\overline{\nu}^*$ refers to the pullack morphism $\mathcal{O}\to R\overline{\nu}_*\hat{\mathcal{O}}^+$.

Throughout this discussion we let $L\eta$ mean $L\eta_{\zeta_p-1}$. Let $\tilde{\Omega}=L\eta R\overline{\nu}_*\hat{\mathcal{O}}^+$ be the sheaf introduced in Section 5.2. The next key step is to observe that $\overline{\nu}^*\colon\mathcal{O}\to R\overline{\nu}_*\hat{\mathcal{O}}^+$ factors as a composition of $\Gamma_K$-equivariant maps:

$$\mathcal{O}\xrightarrow{i}\tilde{\Omega}\xrightarrow{a}R\overline{\nu}_*\hat{\mathcal{O}}^+$$

To see this, apply the naturality of $a$ to the morphism $\overline{\nu}^*\colon\mathcal{O}\to R\nu_*\hat{\mathcal{O}}^+$ to find a commutative diagram:

$$\begin{array}{ccc}
L\eta\mathcal{O} & \xrightarrow{L\eta(\overline{\nu}^*)} & L\eta R\overline{\nu}_*\hat{\mathcal{O}}^+ \\
a\big\downarrow & & \big\downarrow a \\
\mathcal{O} & \xrightarrow[\overline{\nu}^*]{} & R\nu_*\hat{\mathcal{O}}^+
\end{array}$$

Now note that since $\mathcal{O}$ is concentrated in degree 0 and is torsion-free, the map $a\colon L\eta\mathcal{O}\to\mathcal{O}$ is an isomorphism, and we get the desired factorization of $\overline{\nu}^*$.

The complex $\tilde{\Omega}$ is quasi-coherent by Theorem 5.2.5; that is, it is the complex of sheaves on $X_{\text{ét}}$ associated to the complex of $\mathcal{O}(\overline{\mathfrak{X}})$-modules $\tilde{\Omega}(\overline{\mathfrak{X}}):=R\Gamma(\overline{\mathfrak{X}}_{\text{ét}},\tilde{\Omega})$.

The theorem will follow from combining the following two statements:

(1) Let $\tilde{\Omega}_{\geq 1}$ be the cofiber of $i\colon\mathcal{O}(\overline{\mathfrak{X}})\to\tilde{\Omega}(\overline{\mathfrak{X}})$. Then $H^i(\tilde{\Omega}_{\geq 1}^{h\Gamma_K})$ obeys bounds of the same shape as in the theorem.

(2) Regarding $a\colon\tilde{\Omega}\to R\overline{\nu}_*\hat{\mathcal{O}}^+$, we have that $H^i(\operatorname{cof}(a)^{h\Gamma_K})$ is 0 for $i=0$, and is killed by a uniform power of $p$ for $i\geq 1$, namely $p^{2d+1}$.

These statements appear as Lemma 5.4.2 and Lemma 5.5.5 below. □

5.4. **Controlling $\tilde{\Omega}_{\geq 1}^{h\Gamma_K}$.** By Corollary 5.2.6, $\tilde{\Omega}_{\geq 1}$ admits a filtration whose associated graded pieces are

$$\Omega^j_{\log}(\mathfrak{X})\hat{\otimes}_{\mathcal{O}_K}\mathcal{O}_C\left\{-j\right\}.$$

Using the projection formula Lemma 3.4.10, we find that $\tilde{\Omega}_{\geq 1}^{h\Gamma_K}$ has a finite filtration with associated graded pieces

$$\Omega^j_{\log}(\mathfrak{X})\hat{\otimes}_{\mathcal{O}_K}\mathcal{O}_C\left\{-j\right\}^{h\Gamma_K}$$

for $j=1,\ldots,d$.

We can use Theorem 4.4.3 to put bounds on the torsion in the cohomology of $\Gamma_K$ with coefficients in the Breuil–Kisin twists.

**Lemma 5.4.1.** *For all $j\geq 0$ we have:*

(1) $H^0(\Gamma_K, \mathcal{O}_C\{-j\}) = 0$.
(2) $H^1_{\mathrm{cts}}(\Gamma_K, \mathcal{O}_C\{-j\})$ *is* $p^{B_{K^t}+v(j)}$*-torsion, where* $B_{K^t}$ *only depends on* $K^t$*. (In fact we can take* $B_{K^t} = M_{K^t} + 1$ *for the constant* $M_{K^t}$ *appearing in Theorem 4.0.5.)*
(3) *For* $i \geq 2$*,* $H^i_{\mathrm{cts}}(\Gamma_K, \mathcal{O}_C\{-j\})$ *is* $p^{B'}$*-torsion, where* $B'$ *is an absolute constant. (In fact we can take* $B' = 2$ *or* $B' = 4$ *as* $p$ *is odd or even.)*

*Proof.* The key point is that for each $j > 0$ there is an injective $\Gamma_K$-equivariant map $\mathcal{O}_C\{-j\} \to \mathcal{O}_C(-j)$ with cokernel killed by $p$.

Recall the definition of the Breuil–Kisin twist [BMS18, Definition 8.2]: Start with the module of Kähler differentials $\Omega^1_{\mathcal{O}_C/\mathbb{Z}_p}$, and let $\mathcal{O}_C\{1\}$ be its Tate module; i.e., $\mathcal{O}_C\{1\} = \varprojlim \Omega^1_{\mathcal{O}_C/\mathbb{Z}_p}$, where the transition maps in the limit are multiplication by $p$. Then $\mathcal{O}_C\{1\}$ is a free $\mathcal{O}_C$-module of rank 1 carrying a $\Gamma_K$-action.

Let $(\zeta_{p^r})_{r\geq 1}$ be a compatible system of primitive $p$th power roots of 1; then $\omega = (d\zeta_{p^r}/\zeta_{p^r})_{r\geq 1}$ represents an element of $\mathcal{O}_C\{1\}$. For all $\sigma \in \Gamma_K$ we have $\sigma(\omega) = \chi_{\mathrm{cycl}}(\sigma)\omega$, which is to say that $\omega$ induces a $\Gamma_K$-equivariant map $\mathcal{O}_C(1) \to \mathcal{O}_C\{1\}$. The cokernel of this map is killed by $(\zeta_p - 1)$.

After dualizing and taking $j$-fold tensor products, we obtain a $\Gamma_K$-equivariant injective map $\mathcal{O}_C\{-j\} \to \mathcal{O}_C(-j)$ with cokernel killed by $(\zeta_p - 1)^j$. Since $v((\zeta_p - 1)^j) = \frac{j}{p-1}$ we can divide this map by $p^{\lfloor \frac{j}{p-1} \rfloor}$ to obtain a $\Gamma_K$-equivariant injective map $\mathcal{O}_C\{-j\} \to \mathcal{O}_C(-j)$ with cokernel killed by $p$.

The claimed bound now follows from applying the long exact sequence in cohomology together with Theorem 4.0.5. ☐

**Lemma 5.4.2.** *For each* $i \geq 0$*, the cohomology* $H^i(\tilde{\Omega}^{h\Gamma_K}_{\geq 1})$ *is killed by* $p^{C_{K^t}(d)}$*, where* $C_{K^t}(d)$ *is a linear function of* $d$ *which only depends on* $K^t$*.*

*Proof.* The filtration of $\mathrm{cof}(\alpha_2)$ above gives rise to a spectral sequence

$$\Omega^j_{\mathfrak{X},\log}(\mathfrak{X}) \otimes_{\mathcal{O}_K} H^i_{\mathrm{cts}}(\Gamma_K, \mathcal{O}_C\{-j\}) \implies H^{i+j}(\mathrm{cof}(\alpha_2)),$$

where the left side is nonzero only for $j = 1, \ldots, d$ and $i \geq 1$. Therefore a power of $p$ which annihilates $H^s(\mathrm{cof}(\alpha_2))$ can be obtained by summing together the powers of $p$ which annihilate $H^i_{\mathrm{cts}}(\Gamma_K, \mathcal{O}_C)$ for $\min\{s-d, 1\} \leq i \leq s-1$. The proposition then follows from Lemma 5.4.1, noting the crude bound $v(j) \leq d$ whenever $j \leq d$. ☐

5.5. **Controlling the décalage functor.** We record here some lemmas regarding the $L\eta$ functor, recalling and extending slightly the results in §6 of [BMS18]. Let $(T, \mathcal{O}_T)$ be a ringed topos, and let $D(\mathcal{O}_T)$ be the derived category of $\mathcal{O}_T$-modules. Let $\mathcal{I} \subset \mathcal{O}_T$ be an invertible ideal sheaf. We use $L\eta_{\mathcal{I}}$ to denote the lax symmetric monoidal functor $D(\mathcal{O}_T) \to D(\mathcal{O}_T)$ as in [BMS18, Corollary 6.5, Proposition 6.7]. This functor has the effect of killing the $\mathcal{I}$-torsion in the cohomology of a complex. Throughout this section, we will say that an object in $D(\mathcal{O}_T)$ is $\mathcal{I}^s$*-torsion* for some $s \geq 1$ if it admits the structure of an $\mathcal{O}_T/\mathcal{I}^s$-module.

The functor $L\eta_{\mathcal{I}}$ commutes with truncations and in particular preserves the subcategories of bounded complexes $D^{\geq 0}(\mathcal{O}_T)$, $D^{\leq d}(\mathcal{O}_T)$, $D^{[0,d]}(\mathcal{O}_T)$.

**Lemma 5.5.1.** *Let* $\mathcal{C}$ *be an object in* $D(\mathcal{O}_T)$*. Then*

(1) *Assume that* $\mathcal{C} \in D^{\geq 0}(\mathcal{O}_T)$ *and that* $H^0(\mathcal{C})$ *is* $\mathcal{I}$*-torsion free. We have a natural map in* $D(\mathcal{O}_T)$*:*

$$a \colon L\eta_{\mathcal{I}}(\mathcal{C}) \to \mathcal{C}.$$

(2) *Assume that* $\mathcal{C} \in D^{\leq d}(\mathcal{O}_T)$*. We have a natural map in* $D(\mathcal{O}_T)$

$$b \colon \mathcal{C} \otimes \mathcal{I}^{\otimes d} \to L\eta_{\mathcal{I}}(\mathcal{C})$$

(3) *Assume that* $\mathcal{C} \in D^{[0,d]}(\mathcal{O}_T)$ *and that* $H^0(\mathcal{C})$ *is* $\mathcal{J}$*-torsion free. The cofibers of* $b \circ (a \otimes \mathcal{J}^{\otimes d})$ *and* $a \circ b$ *are* $\mathcal{J}^{\otimes d}$*-torsion.*

*Proof.* These claims are all part of [BMS18, Lemma 6.9]. □

Applying Lemma 5.5.1 to the ringed topos $(\overline{\mathfrak{X}}_{\text{ét}}, \mathcal{O})$, the invertible ideal $\mathcal{J} = (\zeta_p - 1)$, and the objects $R\overline{\nu}_* \hat{\mathcal{O}}^+$ and $\tau^{\leq d} R\overline{\nu}_* \hat{\mathcal{O}}^+$, we obtain morphisms:

$$\begin{aligned} a\colon L\eta_{\mathcal{J}} R\overline{\nu}_* \hat{\mathcal{O}}^+ &\to R\overline{\nu}_* \hat{\mathcal{O}}^+ \\ \overline{a}\colon L\eta_{\mathcal{J}} \tau^{\leq d} R\overline{\nu}_* \hat{\mathcal{O}}^+ &\to \tau^{\leq d} R\overline{\nu}_* \hat{\mathcal{O}}^+ \\ b\colon \tau^{\leq d} R\overline{\nu}_* \hat{\mathcal{O}}^+ \otimes \mathcal{J}^{\otimes d} &\to L\eta_{\mathcal{J}} R\overline{\nu}_* \hat{\mathcal{O}}^+ \end{aligned}$$

**Lemma 5.5.2.** *The object* $\operatorname{cof}(\overline{a})$ *is* $\mathcal{J}^{2d}$*-torsion.*

*Proof.* Consider the diagram of exact triangles constructed from the octahedral axiom:

$$\begin{array}{ccccc} \operatorname{cof}(b \circ (\overline{a} \otimes \mathcal{J}^{\otimes d})) & \xrightarrow{=} & \operatorname{cof}(b \circ (\overline{a} \otimes \mathcal{J}^{\otimes d})) & \longrightarrow & 0 \\ \downarrow & & \downarrow & & \downarrow \\ \operatorname{cof}(b) & \longrightarrow & \operatorname{cof}(\overline{a}b) & \longrightarrow & \operatorname{cof}(\overline{a}) \\ \downarrow & & \downarrow & & \downarrow = \\ \mathcal{J}^{\otimes d} \otimes \operatorname{cof}(\overline{a})[1] & \longrightarrow & N & \longrightarrow & \operatorname{cof}(\overline{a}). \end{array}$$

Here $N$ is the cofiber of the map $\operatorname{cof}(b \circ (\overline{a} \otimes \mathcal{J}^{\otimes d})) \to \operatorname{cof}(\overline{a}b)$. Since $\operatorname{cof}(b \circ (\overline{a} \otimes \mathcal{J}^{\otimes d}))$ and $\operatorname{cof}(\overline{a}b)$ are both $\mathcal{J}^d$-torsion, $N$ is $\mathcal{J}^{2d}$-torsion. By Lemma 5.5.3 below, the map $\operatorname{cof}(\overline{a}) \to \mathcal{J}^{\otimes d} \otimes \operatorname{cof}(\overline{a})[2]$ is 0, so that

$$N \cong \left(\mathcal{J}^{\otimes d} \otimes \operatorname{cof}(\overline{a})[1]\right) \oplus \operatorname{cof}(\overline{a}).$$

It follows that $\operatorname{cof}(\overline{a})$ is $\mathcal{J}^{2d}$-torsion. □

**Lemma 5.5.3.** *Let* $a_1, a_2, a_3$ *be morphisms in a triangulated category which fit into a diagram:*

$$A \xrightarrow{a_1} B \xrightarrow{a_2} C \xrightarrow{a_3} D.$$

*Then the composition of the canonical maps*

$$\operatorname{cof}(a_3) \to \operatorname{cof}(a_2)[1] \to \operatorname{cof}(a_1)[2]$$

*is 0.*

*Proof.* We have a commutative diagram of fiber sequences:

$$\begin{array}{ccccccc} \operatorname{cof}(a_2 a_1) & \longrightarrow & \operatorname{cof}(a_3 a_2 a_1) & \longrightarrow & \operatorname{cof}(a_3) & \longrightarrow & \operatorname{cof}(a_2 a_1)[1] \\ \downarrow & & \downarrow & & = \downarrow & & \downarrow \\ \operatorname{cof}(a_2) & \longrightarrow & \operatorname{cof}(a_3 a_2) & \longrightarrow & \operatorname{cof}(a_3) & \longrightarrow & \operatorname{cof}(a_2)[1] \end{array}$$

It follows that the canonical map $\operatorname{cof}(a_3) \to \operatorname{cof}(a_2)[1]$ factors through $\operatorname{cof}(a_2 a_1)[1]$, which is exactly the fiber of the canonical map $\operatorname{cof}(a_2)[1] \to \operatorname{cof}(a_1)[2]$. Thus the composition $\operatorname{cof}(a_3) \to \operatorname{cof}(a_1)[2]$ is zero. □

**Lemma 5.5.4.** *The complexes* $\operatorname{cof}(a)$ *and* $\operatorname{cof}(\overline{a})$ *are almost isomorphic. In particular,* $\operatorname{cof}(a)$ *is* $(\zeta_p - 1)^{2d+1}$*-torsion.*

*Proof.* Let

$$t\colon \tau^{\leq d} R\overline{\nu}_* \hat{\mathcal{O}}^+ \to R\overline{\nu}_* \hat{\mathcal{O}}^+$$

be the natural map. Consider the diagram of cofiber sequences:

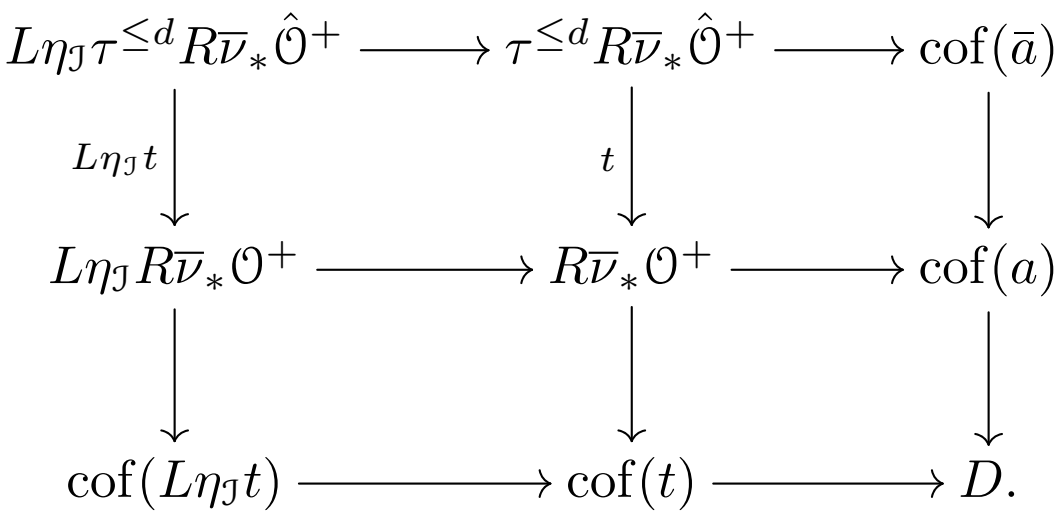

We claim that $t$ and $L\eta_{\mathcal{J}} t$ are almost isomorphisms. We begin with $t$. It suffices to verify that $t$ is an almost isomorphism after evaluating the source and target complexes of sheaves on an open cover. Given an open $U \to \mathfrak{X}'_{\text{ét}}$, $(R\overline{\nu}_* \mathcal{O}^+)(U) \cong R\Gamma((\overline{\nu})^{-1}(U)_{\text{proét}}, \mathcal{O}^+)$ is subject to the constraints imposed by the discussion around [ČK19, Equation 3.3.1], in which $\Delta$ has cohomological dimension $d$, so the cofiber of the map $e$ is almost zero. Applied to our situation, this implies that the cofiber of $t$ is almost zero. A similar argument applying the source and target of $L\eta_{\mathcal{J}} t$ to an open allows us to apply [ČK19, Theorem 3.9] to see that the cofiber of $L\eta_{\mathcal{J}} t$ is almost zero.

Therefore $\mathrm{cof}(L\eta_{\mathcal{J}} t)$ and $\mathrm{cof}(t)$ are both almost zero. This implies that $D$ is almost zero, and this further implies that $\mathrm{cof}(\overline{a}) \to \mathrm{cof}(a)$ is an almost isomorphism. Since an almost zero sheaf of complexes is $\mathcal{J}$-torsion, we have that the cohomology of $\mathrm{cof}(a)$ is $\mathcal{J}^{2d+1}$-torsion. □

**Lemma 5.5.5.** *For $i \geq 0$ we have*

$$H^i(\mathrm{cof}(a)^{h\Gamma_K}) = \begin{cases} 0, & i = 0 \\ p^{2d+1}\textit{-torsion}, & i \geq 1. \end{cases}$$

*Proof.* Since $\mathrm{cof}(a)$ is $(\zeta_p - 1)^{2d+1}$-torsion by Lemma 5.5.4, we (rather crudely) conclude that it is $p^{2d+1}$-torsion. It follows that the homotopy fixed points $\mathrm{cof}(a)^{h\Gamma_K}$ are also $p^{2d+1}$-torsion: Indeed, this can be seen by using the observation that taking homotopy fixed points is lax monoidal, so $\mathrm{cof}(a)^{h\Gamma_K}$ is also a module over $\mathbb{Z}/p^{2d+1}\mathbb{Z}$. This implies the claim about the cohomology groups of $\mathrm{cof}(a^{h\Gamma_K})$. □

5.6. **Tame descent.** Theorem 5.3.4 gives control over the pro-étale cohomology of a rigid-analytic space $X/K$ which admits an affine semistable model over $\mathcal{O}_K$. We now pose the question of whether this control may be obtained for the general case of a smooth rigid-analytic space $X/K$. The trouble we encounter is that $X/K$ does not necessarily have a semistable model; it may be necessary to extend scalars. The best result we know along these lines is the main theorem of [Har03], which states that $X_{\text{ét}}$ admits an open cover by rigid-analytic generic fibers of semistable formal schemes defined over finite extensions $L/K$. If $X$ is not quasi-compact then there may not be a single $L/K$ that suffices for this purpose, and then the constants appearing in Theorem 5.3.4 could a priori accumulate.

In the special case that all the $L/K$ are *tamely ramified* extensions (or, at least, if the wild ramification of the $L/K$ is bounded), then it becomes possible to get uniform control over the pro-étale cohomology of $X$. We will see in the next section that this is the case for the open unit ball over $K$.

We gather a few lemmas about descent through tame extensions. For a Galois extension of local fields $L/K$ with group $G$, let $\mathcal{O}_L \langle G \rangle$ denote the associative $\mathcal{O}_L$-algebra whose underlying module is freely generated by symbols $[g]$, and whose multiplication law is determined by the relations $[gh] = [g][h]$ and $[g]\alpha = g(\alpha)[g]$ for $\alpha \in \mathcal{O}_L$ and $g, h \in G$. Let $\mathrm{Mod}_{\mathcal{O}_L\langle G\rangle}$ denote the category of $\mathcal{O}_L\langle G\rangle$-modules. An object of this category is an $\mathcal{O}_L$-module $A$ endowed with an action of $G$ which acts semilinearly over the action of $G$ on $\mathcal{O}_L$; i.e., a descent datum along $L/K$. Its module of fixed points $A^G$ is an $\mathcal{O}_K$-module.

**Lemma 5.6.1.** *Let $K$ be a local field of characteristic $(0, p)$, and let $L/K$ be a finite tame Galois extension with group $G$. The functor $A \mapsto A^G$ from $\mathcal{O}_L\langle G\rangle$-modules to $\mathcal{O}_K$-modules is exact. As a result, if $C \in D(\mathcal{O}_L\langle G\rangle)$ is a complex whose cohomology groups are all killed by $p^k$ for some $k \geq 1$, then the same is true for $C^{hG} \in D(\mathcal{O}_K)$.*

*Proof.* Since a composition of exact functors is exact, we can reduce to the cases where $L/K$ is unramified or totally tamely ramified.

In the case where $L/K$ is unramified, $\mathcal{O}_L$ is an étale $G$-torsor over $\mathcal{O}_K$, and thus by étale descent, the functor $A \mapsto A^G$ is an equivalence of categories.

In the case where $L/K$ is totally tamely ramified, exactness follows from the fact that the order of $G$ is a unit in $\mathcal{O}_K$. ☐

**Lemma 5.6.2.** *Suppose $X/K$ is a rigid-analytic space. Define $\beta_X$ as the composition:*

$$\mathcal{O}^+(X)[\varepsilon] \overset{unit}{\to} R\Gamma(X_{\text{ét}}, \mathcal{O}^+)[\varepsilon] \overset{\nu_X}{\to} R\Gamma(X_{\text{proét}}, \hat{\mathcal{O}}^+).$$

*Assume there exists a finite tame Galois extension $L/K$ such that $X_L$ is the generic fiber of an affine semistable formal scheme $\mathfrak{X}$ of dimension $d$ over $\mathcal{O}_L$. Then the cohomology groups $H^i(\mathrm{cof}(\beta_X))$ are all killed by a uniform power of $p$ which only depends on $K$.*

*Proof.* Let $G = \mathrm{Gal}(L/K)$. Theorem 5.3.4 applied to the affine semistable formal scheme $\mathfrak{X}$ gives a bound on the cofiber of the $G$-equivariant map

$$\beta_{X_L} : \mathcal{O}^+(X_L)[\varepsilon] \to R\Gamma(X_{L,\text{proét}}, \hat{\mathcal{O}}^+).$$

Here we have used Lemma 5.2.3 to identify $\mathcal{O}(\mathfrak{X})$ with $\mathcal{O}^+(X_L)$. Applying homotopy fixed points under $G$ is exact by Lemma 5.6.1, so the same bounds apply to the cofiber of $\beta_X$. ☐

We now pass to the non-affinoid case. To do this, we consider covers of a rigid-analytic space by affinoids satisfying the hypotheses of Lemma 5.6.2. At this point we will start encountering objects in $D(\mathrm{Solid})$ which are limits of $p$-adic objects but which are no longer $p$-adic themselves.

**Definition 5.6.3.** Let $X/K$ be a smooth rigid analytic space of dimension $d$. Let $\mathcal{U} \coloneqq \{U_i\}_{i\in I}$ be an open cover in $X_{\mathrm{an}}$. For a finite subset $J \subset I$, let $U_J$ be the intersection of the $\{U_i\}_{i\in J}$. We say that $\mathcal{U}$ is *tamely semistable* if for every such $J$, there exists a finite Galois tame extension $L_J/K$ such that $(U_J)_{L_J}$ admits an affine semistable formal model over $\mathcal{O}_{L_J}$.

Such a cover appears when there is already a semistable model for $X$:

**Lemma 5.6.4.** *Let $\mathfrak{X}$ be a quasi-separated semistable formal scheme over $\mathcal{O}_K$ with generic fiber $X$. Let $\{\mathfrak{U}_i\}_{i\in I}$ be an affine open cover of $\mathfrak{X}$. Let $\mathcal{U}$ be the cover of $X$ obtained by pulling this cover back through the reduction map. Then $\mathcal{U}$ is a tamely semistable cover for $X$.*

*Proof.* The quasi-separatedness assumption implies that for all finite non-empty $J \subset I$ the open subset $\mathfrak{U}_J \subset \mathfrak{X}$ is affine. Let

$$\mathrm{red} \colon X \to \mathfrak{X}$$

be the reduction map. For each $J$ we have $U_J \cong \mathrm{red}^{-1}(\mathfrak{U}_J)$. Since semistablity is an étale local condition, we see that $\mathfrak{U}_J$ is an affine semistable model for $U_J$. ☐

**Definition 5.6.5.** Let $X/K$ be a rigid-analytic space. Let $\mathcal{U}$ be an affinoid open cover in $X_{\mathrm{an}}$, with the same notational conventions as in Definition 5.6.3. Define the *condensed Čech complex* by

$$\check{C}(\mathcal{U}, \mathcal{O}^+_{\mathrm{cond}}) \coloneqq \lim_{[n]\in\Delta} \prod_{f\in I^{[n]}} \mathcal{O}^+_{\mathrm{cond}}(U_{\mathrm{im}(f)}).$$

In this expression:

- $\Delta$ is the simplex category, with objects $[n] = \{0, 1, \dots, n\}$ for $n \geq 0$, and morphisms order preserving maps,
- $I^{[n]}$ is the set of functions $f \colon [n] \to I$, and $\mathrm{im}(f)$ is the image of $f$, and
- the limit is computed in $D(\mathrm{Cond}(\mathrm{Ab}))$.

Thus, $\check{C}(\mathcal{U}, \mathcal{O}^+_{\mathrm{cond}})$ is a condensed enhancement of the usual Čech complex $\check{C}(\mathcal{U}, \mathcal{O}^+)$.

*Remark* 5.6.6. Since $\mathcal{O}^+_{\mathrm{cond}}$ (in contrast to $\mathcal{O}_{\mathrm{cond}}$) is not acyclic on affinoids, in general $\check{C}(\mathcal{U}, \mathcal{O}^+_{\mathrm{cond}})$ does not just compute the $\mathcal{O}^+_{\mathrm{cond}}$-cohomology of $X$.

**Lemma 5.6.7.** *Let $\mathfrak{X}$ be a quasi-separated semistable formal scheme over $\mathcal{O}_K$ with generic fiber $X$. Let $\mathcal{U}$ be the affinoid open cover of $X$ obtained by pulling back an affine open cover $\{\mathfrak{U}_i\}_{i\in I}$ of $\mathfrak{X}$. Then $R\Gamma(\mathfrak{X}, \mathcal{O}_{\mathrm{cond}})$ is quasi-isomorphic to $\check{C}(\mathcal{U}, \mathcal{O}^+_{\mathrm{cond}})$.*

*Proof.* By Lemma 3.5.1, $\mathcal{O}_{\mathrm{cond}}$ is acyclic on affine formal schemes. Therefore $R\Gamma(\mathfrak{X}, \mathcal{O}_{\mathrm{cond}})$ is computed by Čech cohomology with respect to the affine cover $\{\mathfrak{U}_i\}_{i\in I}$:

$$R\Gamma(\mathfrak{X}, \mathcal{O}_{\mathrm{cond}}) \cong \lim_{[n]\in\Delta} \prod_{f\in I^{[n]}} \mathcal{O}_{\mathrm{cond}}(\mathfrak{U}_{\mathrm{im}(f)}).$$

For each finite $J \subset I$ we have an isomorphism of topological rings Lemma 5.2.3 $\mathcal{O}(\mathfrak{U}_J) \cong \mathcal{O}^+(U_J)$. Applying the functor $A \mapsto \underline{A}$ from topological abelian groups to condensed abelian groups, we find that $R\Gamma(\mathfrak{X}, \mathcal{O}_{\mathrm{cond}})$ is quasi-isomorphic to $\check{C}(\mathcal{U}, \mathcal{O}^+_{\mathrm{cond}})$. $\square$

The next theorem shows that if $\mathcal{U}$ is a tamely semistable cover of $X$, then the condensed Čech complex associated to $\mathcal{U}$ is a good approximation for the pro-étale cohomology of $X$.

**Theorem 5.6.8.** *Let $X/K$ be a smooth rigid analytic space of dimension $d$, and let $\mathcal{U} \coloneqq \{U_i\}_{i\in I}$ be a tamely semistable cover of $X$. We have a natural map*

$$\alpha_{\mathcal{U}} \colon \check{C}(\mathcal{U}, \mathcal{O}^+_{\mathrm{cond}})[\varepsilon] \to R\Gamma(X_{\mathrm{pro\acute{e}t}}, \hat{\mathcal{O}}^+_{\mathrm{cond}}).$$

*The cohomology groups $H^i(\mathrm{cof}(\alpha_{\mathcal{U}}))$ are each killed by a power of $p$. If the dimension of the nerve of $\mathcal{U}$ is finite, the same power of $p$ kills $H^i(\mathrm{cof}(\alpha_{\mathcal{U}}))$ for all $i$.*

*Proof.* For each finite non-empty set of indices $J \subset I$, we have a natural map

$$\beta_J \colon \mathcal{O}^+(U_J)[\varepsilon] \to R\Gamma(U_{J,\mathrm{pro\acute{e}t}}, \hat{\mathcal{O}}^+)$$

as in Lemma 5.6.2. By that lemma, there exists $n \geq 1$ such that $p^n$ kills $H^i(\mathrm{cof}(\beta_J))$ for all $i, J$.

We consider now the condensed enhancement of $\beta_J$. Note that $U_J$ is affinoid, and as a topological ring, $\mathcal{O}^+(U_J)$ is $p$-adic, in the sense that $\mathcal{O}^+(U_J) \cong \varprojlim \mathcal{O}^+(U_J)/p^n$ as a limit of discrete groups. On the other hand, Lemma 3.7.1 states that $R\Gamma(U_{J,\mathrm{pro\acute{e}t}}, \hat{\mathcal{O}}^+_{\mathrm{cond}})$ is $p$-adic as well. Therefore we can apply the functor $X \mapsto \varprojlim X/p^n$ (where the limit is calculated in $D(\mathrm{Cond}(\mathrm{Ab}))$) to $\beta_J$, we obtain a morphism in $D(\mathrm{Cond}(\mathrm{Ab}))$ which is functorial in $J$:

$$\beta_{J,\mathrm{cond}} \colon \mathcal{O}^+_{\mathrm{cond}}(U_J)[\varepsilon] \to R\Gamma(U_{J,\mathrm{pro\acute{e}t}}, \hat{\mathcal{O}}^+_{\mathrm{cond}}).$$

The $H^i(\mathrm{cof}(\beta_{J,\mathrm{cond}}))$ are $p$-adic solid abelian groups killed by $p^n$ (which implies they are discrete).

Combining the $\beta_{J,\mathrm{cond}}$ together, we find a morphism $\alpha_{\mathcal{U}}$ as in the theorem. Its cofiber is

$$\mathrm{cof}(\alpha_{\mathcal{U}}) \cong \lim_{[n]\in\Delta} \prod_{f\in I^{[n]}} \mathrm{cof}(\beta_{\mathrm{im}(f)}).$$

Taking the associated cosimplicial limit spectral sequence

$$E_1^{s,q} = \prod_{f\in I^{[q]}} H^s(\mathrm{cof}(\beta_{\mathrm{im}(f)})) \implies H^{s+q}(\mathrm{cof}(\alpha_{\mathcal{U}})),$$

we find that $H^i(\mathrm{cof}(\alpha_{\mathcal{U}}))$ is killed by $p^{n(i+1)}$. If the dimension of the nerve of the cover is $\delta$, we find that $H^i(\mathrm{cof}(\alpha_{\mathcal{U}}))$ is killed by $p^{n(\delta+1)}$. □

**Corollary 5.6.9.** *Let $\mathfrak{X}$ be a quasi-separated semistable formal scheme over $\mathcal{O}_K$ with generic fiber $X$. There is a natural morphism in $D(\mathrm{Cond})$:*

$$\alpha_{\mathfrak{X}} \colon R\Gamma(\mathfrak{X}, \mathcal{O}_{\mathrm{cond}})[\varepsilon] \to R\Gamma(X_{\mathrm{pro\acute{e}t}}, \hat{\mathcal{O}}^+_{\mathrm{cond}}),$$

*such that the cohomology groups of $\mathrm{cof}(\alpha_{\mathfrak{X}})$ are $p$-power torsion.*

*Proof.* This is immediate from Lemma 5.6.7 and Theorem 5.6.8. □

Applying $\Gamma_*$ to the morphism in Corollary 5.6.9 and inverting $p$ gives the non-condensed statement appearing in Theorem D.

# 6. Proof of Theorem B

Let us recall the main players in Theorem B: a prime number $p$, an integer $n \geq 1$, the ring of $p$-typical Witt vectors $W = W(\overline{\mathbb{F}}_p)$, and the Lubin–Tate ring $A \cong W[\![u_1, \ldots, u_{n-1}]\!]$, which admits a continuous action of the Morava stabilizer group $\mathbb{G}_n$. Proposition 2.5.1 states that the inclusion $W \hookrightarrow A$ admits a continuous $\mathbb{G}_n$-equivariant additive splitting, say with complement $A^c$. Theorem B is the statement that the continuous cohomology $H^*_{\mathrm{cts}}(\mathbb{G}_n, A^c)$ is $p$-power torsion.

In Section 3.9 we explained how to reduce Theorem B to a statement (Theorem 3.9.4) controlling the pro-étale cohomology of the open ball and Drinfeld's symmetric space. We examine these cases in turn.

## 6.1. Pro-étale cohomology of the open ball.

Let $K$ be a local field of characteristic $(0, p)$. For an integer $d \geq 1$, recall that we had defined the $d$-dimensional rigid-analytic open ball over $K$:

$$B^{\circ,d} = (\mathrm{Spa}\, \mathcal{O}_K[\![T_1, \ldots, T_d]\!]) \setminus \{|p| = 0\}$$

The ball $B^{\circ,d}$ is not quasi-compact. It is exhausted by affinoid (closed) balls of increasing radius:

$$B^{\circ,d} = \varinjlim_{r<1} B^d_r, \tag{6.1.1}$$

where $r$ runs over real numbers in $|p|^{\mathbb{Q}_+}$, and for each $r = |p|^{m/n}$ with $m, n$ relatively prime positive integers, $B^d_r$ is the rational subset defined by the inequality $|T|^n \leq |p|^m$.

The goal of this section is to control $R\Gamma(B^{\circ,d}_{\mathrm{pro\acute{e}t}}, \hat{\mathcal{O}}^+)$. The idea is to apply Theorem 5.6.8 to the cover $\mathcal{U}$ given by a collection of affinoid balls $B^d_r$ which cover $B^{\circ,d}$. A convenient choice of radii is $r_\ell = |p|^{1/\ell}$, where $\ell$ runs over prime numbers $\neq p$. The closed ball $B^d_{r_\ell}$ admits a smooth formal model after passage from $K$ to the tamely ramified extension

$$L_\ell = K(p^{1/\ell}).$$

Indeed, $B^d_{r_\ell}$ is the rigid-analytic generic fiber of the smooth formal scheme $\mathfrak{B}^d_{r_\ell}$, where

$$\mathfrak{B}^d_{r_\ell} = \mathrm{Spf}\, \mathcal{O}_{L_\ell} \left\langle \frac{T_1}{p^{1/\ell}}, \ldots, \frac{T_d}{p^{1/\ell}} \right\rangle.$$

Therefore the cover $\mathcal{U}$ is tamely semistable.

The solid Čech complex $\check{C}(\mathcal{U}, \mathcal{O}^+_{\mathrm{cond}})$ is quasi-isomorphic to $R\varprojlim_\ell H^0(B^d_{r_\ell}, \mathcal{O}^+_{\mathrm{cond}})$, which could only have an $H^0$ and an $H^1$. The $H^0$ is

$$\varprojlim H^0(B^d_{r_\ell}, \mathcal{O}^+_{\mathrm{cond}}) = \mathcal{O}_K[\![T_1, \ldots, T_d]\!],$$

where the latter must be considered as the solid abelian group associated to the topological ring $\mathcal{O}_K[\![T_1, \ldots, T_d]\!]$ (with its $(p, T_1, \ldots, T_d)$-adic topology). We claim there is no $H^1$. In order to prove this, we begin with a formal lemma.

**Lemma 6.1.2.** *If $A_i$ is a cofiltered diagram of discrete abelian groups satisfying the Mittag-Leffler condition and $S$ is a profinite set, then the diagram of abelian groups $C_{\mathrm{cts}}(S, A_i)$ satisfies the Mittag-Leffler condition as well.*

*Proof.* This follows from the fact that for a map of discrete sets $f\colon A \to B$ and a profinite set $S$, we have

$$\mathrm{Im}\,[C_{\mathrm{cts}}(S, A) \to C_{\mathrm{cts}}(S, B)] = C_{\mathrm{cts}}(S, \mathrm{Im}(f)).$$

Indeed, the functions involved here are all locally constant, which reduces us to the obvious case where $S$ is finite. □

**Lemma 6.1.3.** *Let $K$ be a local field of characteristic $(0, p)$ with uniformizer $\varpi$. Let $j \in \mathbb{N}$ and let $0 < r < 1$ and set $M = \lceil \frac{1+j}{\log_{|\varpi|}(r)} \rceil$. Then for every $1 > s \geq r^{\frac{1}{j+1}}$ and $j' \geq j$ the image of the map*

$$H^0(B^d_s, \mathcal{O}^+)/\varpi^{j'} \to H^0(B^d_r, \mathcal{O}^+)/\varpi^j$$

*is the $\mathcal{O}_K/\varpi^j$-submodule spanned by the $T^i$ for multi-indices $i$ with $|i| < M$.*

*Proof.* Denote by $V \subset H^0(B^d_r, \mathcal{O}^+)/\varpi^j$ the $\mathcal{O}_K/\varpi^j$-submodule spanned by $T^i$ for $|i| < M$. Let $b = \sum_i b_i T^i \in H^0(B^d_s, \mathcal{O}^+)$ and define two elements in $H^0(B^d_s, \mathcal{O}^+)$ by

$$b' = \sum_{|i|<M} b_i T^i \quad \text{and} \quad b'' = \sum_{|i|\geq M} b_i T^i.$$

We have $b = b' + b''$. It is enough show that

(1) the image of $b'$ in $H^0(B^d_r, \mathcal{O}^+)/\varpi^j$ belongs to $V$;
(2) and the image of $b''$ in $H^0(B^d_r, \mathcal{O}^+)/\varpi^j$ is zero.

Thus it is enough to show

(1) $|b_i| \leq 1$ for $|i| < M$;
(2) and $|\frac{b_i}{\varpi^j}| r^{|i|} \leq 1$ for $|i| \geq M$.

Since the value group of $K$ is isomorphic to $|\varpi|^{\mathbb{Z}}$, Condition (1) is equivalent to $|b_i| < \frac{1}{|\varpi|}$ for $|i| < M$. Now since $b \in H^0(B^d_s, \mathcal{O}^+)$, we have for all $i$ that $|b_i| \leq s^{-|i|} \leq r^{-\frac{|i|}{j+1}}$. Thus we get

$$|b_i| \leq r^{\frac{-(M-1)}{j+1}}$$

for $|i| < M$ and

$$|\frac{b_i}{\varpi^j}| r^{|i|} = |\varpi|^{-j} |b_i| r^{|i|} \leq |\varpi|^{-j} r^{\frac{j|i|}{j+1}} \leq |\varpi|^{-j} r^{\frac{jM}{j+1}}$$

for $|i| \geq M$. It is thus enough to show

$$r^{\frac{-(M-1)}{j+1}} < \frac{1}{|\varpi|} \quad \text{and} \quad |\varpi|^{-j} r^{\frac{jM}{j+1}} \leq 1.$$

Rearranging these two inequalities, we obtain the following condition:

$$M \geq \frac{j+1}{\log_{|\varpi|}(r)} > M - 1,$$

which is satisfied for for $M = \lceil \frac{1+j}{\log_{|\varpi|}(r)} \rceil$, as needed. □

**Theorem 6.1.4.** *Fix $d \geq 1$. There is a morphism in $D(\mathrm{Solid})$:*

$$\alpha_{B^{d,\circ}} : \underline{\mathcal{O}_K[\![T_1, \dots, T_d]\!]}[\varepsilon] \to R\Gamma(B^{d,\circ}_{\text{proét}}, \hat{\mathcal{O}}^+_{\text{cond}}),$$

*such that for all $i \geq 0$, $H^i(\mathrm{cof}(\alpha_{B^{d,\circ}}))$ is killed by a power of $p$ which does not depend on $i$.*

*Proof.* We apply Theorem 5.6.8 to the tamely semistable cover $\mathcal{U} = \{B^d_{r_\ell}\}$. The condensed Čech complex $\check{C}(\mathcal{U}, \mathcal{O}^+_{\text{cond}})$ is quasi-isomorphic to $R\varprojlim_\ell H^0(B^d_{r_\ell}, \mathcal{O}^+_{\text{cond}})$. By Lemma 6.1.3 and Lemma 6.1.2, the $\mathbb{N}^2$-diagram $H^0(B^d_{r_\ell}, \mathcal{O}^+_{\text{cond}})/(p^j)(S)$ satisfies the Mittag-Leffler condition for each profinite set $S$. We then compute:

$$\begin{aligned}
R\varprojlim_\ell H^0(B^d_{r_\ell}, \mathcal{O}^+_{\text{cond}}) &\cong R\varprojlim_\ell R\varprojlim_j H^0(B^d_{r_\ell}, \mathcal{O}^+_{\text{cond}})/p^j \\
&\cong R\varprojlim_{\ell,j} H^0(B^d_{r_\ell}, \mathcal{O}^+_{\text{cond}})/p^j \\
&\cong \varprojlim_{\ell,j} H^0(B^d_{r_\ell}, \mathcal{O}^+_{\text{cond}})/p^j \\
&\cong \varprojlim_{\ell} H^0(B^d_{r_\ell}, \mathcal{O}^+_{\text{cond}}) \\
&\cong \underline{\mathcal{O}_K[\![T_1, \dots, T_d]\!]}.
\end{aligned}$$

The first and forth isomorphisms follow from the derived $p$-completeness of $H^0(B^d_{r_\ell}, \mathcal{O}^+_{\text{cond}})$. The third isomorphism is a consequence of the Mittag-Leffler condition. □

**Corollary 6.1.5.** *Let $W = W(\overline{\mathbb{F}}_p)$, let $K = W[1/p]$, and let $\mathrm{LT}_K$ be Lubin–Tate space in height $n$, with $A = H^0(\mathrm{LT}_K, \mathcal{O}^+) \cong \mathcal{O}_K[\![u_1, \dots, u_{n-1}]\!]$ the Lubin–Tate ring. We have an isomorphism of graded $K$-vector space objects in $D(\mathrm{Solid})$:*

$$H^*(R\Gamma(\mathrm{LT}_{K,\text{proét}}, \hat{\mathcal{O}}^+)^{h\mathcal{O}_D^\times}) \otimes_W K \cong \Lambda_K(x_1, x_3, \dots, x_{2n-1})[\varepsilon] \oplus ((A^c)^{h\mathcal{O}_D^\times} \otimes_W K)[\varepsilon].$$

*Proof.* $\mathrm{LT}_K$ is a rigid-analytic open ball, so Theorem 6.1.4 applies. After taking homotopy fixed points under $\mathcal{O}_D^\times$, we find a morphism

$$\alpha^{h\mathcal{O}_D^\times}_{\mathrm{LT}_K} : A^{h\mathcal{O}_D^\times}[\varepsilon] \to R\Gamma(\mathrm{LT}_{K,\text{proét}}, \hat{\mathcal{O}}^+_{\text{cond}})^{h\mathcal{O}_D^\times}$$

The cohomology groups of $\mathrm{cof}(\alpha^{h\mathcal{O}_D^\times}_{\mathrm{LT}_K}) \cong \mathrm{cof}(\alpha_{\mathrm{LT}_K})^{h\mathcal{O}_D^\times}$ are all $p$-power torsion. After inverting $p$, we arrive at an isomorphism of $K$-vector space objects in $D(\mathrm{Solid})$:

$$A^{h\mathcal{O}_D^\times} \otimes_W K[\varepsilon] \cong R\Gamma(\mathrm{LT}_{K,\text{proét}}, \hat{\mathcal{O}}^+_{\text{cond}})^{h\mathcal{O}_D^\times} \otimes_W K.$$

Now we apply the decomposition $A = W \oplus A^c$ from Proposition 2.5.1. By Proposition 3.8.1 (after extending scalars to $K$), we have $W^{h\mathcal{O}_D^\times} \otimes_W K \cong K^{h\mathcal{O}_D^\times}$ is an exterior algebra over $K$. □

6.2. **The pro-étale cohomology of Drinfeld's symmetric space.** Let $\mathcal{H}$ be Drinfeld's symmetric space of dimension $n-1$. Then $\mathcal{H}$ admits a semistable formal model $\mathfrak{H}/\mathbb{Z}_p$. For a construction of $\mathfrak{H}$, see [GK05a, §6]. The formal scheme $\mathfrak{H}/\mathbb{Z}_p$ admits a continuous action of $\mathrm{GL}_n(\mathbb{Q}_p)$ which is compatible with the isomorphism $\mathfrak{H}_{\mathbb{Q}_p} \cong \mathcal{H}$. The continuity condition means that the stabilizer $H$ of an open affine subscheme $\mathrm{Spf}\, R \subset \mathfrak{H}$ is open and acts continuously on the $p$-complete algebra $R$. As a result, the coherent cohomology $R\Gamma(\mathfrak{H}, \mathcal{O}_{\text{cond}})$ is an object of $D(\mathrm{Solid}_{\mathrm{GL}_n(\mathbb{Q}_p)})$.

Our goal is to control the pro-étale cohomology of $\mathcal{H}$. Key to that calculation is the following acyclicity result for coherent cohomology on the formal model.

**Theorem 6.2.1.** *The unit map*

$$\mathbb{Z}_p \to R\Gamma(\mathfrak{H}, \mathcal{O}_{\mathrm{cond}})$$

*is an isomorphism in* $D(\mathrm{Solid}_{\mathrm{GL}_n(\mathbb{Q}_p)})$. *(Here* $\mathbb{Z}_p$ *with its* $p$*-adic topology is considered as a solid complex, with trivial* $\mathrm{GL}_n(\mathbb{Q}_p)$*-action.)*

*Remark* 6.2.2. The non-condensed version of the theorem is a special case of a result of Grosse-Klönne [GK05b, Theorem 4.5], who proves the acyclicity of $\mathrm{GL}_n(\mathbb{Q}_p)$-equivariant coherent sheaves on $\mathfrak{H}$ which are "strongly dominant", a condition on weights trivially satisfied by the structure sheaf $\mathcal{O}$. We offer a self-contained version of Grosse-Klönne's proof below; ultimately it is a combinatorial calculation taking place on the special fiber of $\mathfrak{H}$.

The result in [GK05b] also computes the cohomology of each of the sheaves $\Omega^j_{\mathfrak{H},\log}$ in terms of certain lattices in Steinberg representations of $\mathrm{PGL}_{d+1}(\mathbb{Q}_p)$. This would allow us to control the $\hat{\mathcal{O}}^+$-cohomology of $\mathcal{H}^d_{\text{proét},C}$. From this we could ultimately compute the continuous cohomology $H^i_{\mathrm{cts}}(\mathcal{O}^1_D, A)[1/p]$, where $\mathcal{O}^1_D \subset \mathcal{O}^\times_D$ is the subgroup of elements of reduced norm 1. We do not pursue these computations here.

*Proof of Theorem 6.2.1.* Let $\mathfrak{H}_s$ denote the special fiber of $\mathfrak{H} \to \operatorname{Spf} \mathbb{Z}_p$. Since $R\Gamma(\mathfrak{H}, \mathcal{O}_{\mathrm{cond}})$ is $p$-adically complete and its $H^0$ is torsion-free (since $\mathfrak{H}$ is semistable over $\operatorname{Spf}\mathbb{Z}_p$, it is also flat), it suffices to show that $R\Gamma(\mathfrak{H}_s, \mathcal{O}_{\mathrm{cond}}) \cong \mathbf{F}_p$.

We need the following facts about $\mathfrak{H}_s$, and in particular its relation to the Bruhat–Tits building BT for $\mathrm{PGL}_n(\mathbb{Q}_p)$ (see [GK05a, §6]):

(1) BT is the simplicial complex whose $q$-simplices $\mathrm{BT}_q$ are the homothety classes of flags of $\mathbb{Z}_p$-lattices in $\mathbb{Q}_p^n$:
$$L_0 \subsetneq L_1 \subsetneq \cdots \subsetneq L_q \subsetneq p^{-1}L_0.$$
In particular BT is $(n-1)$-dimensional.
(2) BT can be expressed as an increasing union of finite contractible subcomplexes.
(3) There is a bijection $\sigma \mapsto Z_\sigma$ between the vertices $\mathrm{BT}_0$ of BT and irreducible components of $\mathfrak{H}_s$.
(4) A collection of vertices $\sigma$ in $\mathrm{BT}_0$ are the vertices of a simplex in BT if and only if the $Z_\sigma$ have nontrivial intersection. In that case, the $Z_\sigma$ meet transversely.
(5) For each $q$-simplex $\tau$, let $Z_\tau$ be the intersection of $Z_\sigma$ over all vertices $\sigma$ of $\tau$. Then $Z_\tau$ is a rational projective variety. More precisely, it is the result of a sequence of blow-ups of $(n-1-q)$-dimensional projective space over $\mathbf{F}_p$, performed along nonsingular subvarieties.

For $\tau \in \mathrm{BT}_q$, we let $\tau_0$ be its set of vertices, and we let $U_\tau$ be the open subset of $\cup_{\sigma\in\tau_0} Z_\sigma$ consisting of points not contained in $Z_{\sigma'}$ for any $\sigma' \notin \tau$. Then the $U_\tau$ are a collection of (generally singular) affine varieties which is closed under intersection.

The $U_\tau$ for $\tau \in \mathrm{BT}_{n-1}$ form an affine cover $\mathfrak{H}_s$, so $R\Gamma(\mathfrak{H}_s, \mathcal{O}_{\mathrm{cond}})$ is computed by the condensed Čech complex

$$\lim_{[q]\in\Delta} \prod_{f\colon [q]\to \mathrm{BT}_{n-1}} \mathcal{O}(U_f), \tag{6.2.3}$$

where $U_f = \cap_{i\in[q]} U_{f(i)}$. Each $\mathcal{O}(U_f)$ is taken to be discrete, and the limit is computed in $D(\mathrm{Solid})$.

For $q > 0$ and $\tau \in \mathrm{BT}_q$, the affine variety $U_\tau$ is the union of irreducible components $Z_\sigma \cap U_\tau$ (for $\sigma$ a vertex of $\tau$) which meet transversely. Each $Z_\sigma \cap U_\tau$ is an affine open in a rational variety $Z_\sigma$, and so it is isomorphic to an open subset of $(n-1)$-dimensional affine space. Therefore there is an isomorphism of $U_\tau$ onto an open subset of $\operatorname{Spec} B$, where

$$B = \mathbf{F}_p[x_1, \dots, x_n]/(x_1 x_2 \cdots x_q),$$

which carries $Z_\sigma \cap U_\tau$ into $\operatorname{Spec} B/(x_i)$ for the corresponding $1 \leq i \leq q$. Thus, the map $B \to \mathcal{O}(U_\tau)$ is flat and we have $\mathcal{O}(Z_\sigma \cap U_\tau) \simeq \mathcal{O}(U_\tau)/(x_i)$.

The resolution of $U_\tau$ by the disjoint union of the $Z_\sigma \cap U_\tau$ gives a quasi-isomorphism:

$$\mathcal{O}(U_\tau) \cong \lim_{[r]\in\Delta} \prod_{\sigma\in \mathrm{BT}_r} \mathcal{O}(Z_\sigma \cap U_\tau). \tag{6.2.4}$$

This follows from the standard observation that the complex of $B$-modules

$$B \to \bigoplus_{1\leq i\leq q} B/(x_i) \to \bigoplus_{1\leq i<j\leq q} B/(x_i,x_j) \to \cdots$$

is exact; base-changing along the flat map $B \to \mathcal{O}(U_\tau)$ gives the isomorphism in (6.2.4).

Each $U_f$ appearing in (6.2.3) equals $U_\tau$ for $\tau = \cap_{i\in[q]} f(i)$. Substituting (6.2.4) into (6.2.3) and rearranging the order of limits shows that $R\Gamma(\mathfrak{H}_s, \mathcal{O}_{\mathrm{cond}})$ is quasi-isomorphic to

$$\lim_{[r]\in\Delta} \prod_{\sigma\in \mathrm{BT}_r} \lim_{[q]\in\Delta} \prod_{f\colon [q]\to \mathrm{BT}_{n-1}} \mathcal{O}(Z_\sigma \cap U_f). \tag{6.2.5}$$

For a given $\sigma \in \mathrm{BT}_q$, the projective variety $Z_\sigma$ admits a cover by affine varieties $Z_\sigma \cap U_\tau$, where $\tau$ runs through $\mathrm{BT}_{n-1}$. Therefore the inner limit in (6.2.5) computes $R\Gamma(Z_\sigma, \mathcal{O})$, and we find:

$$R\Gamma(\mathfrak{H}_s, \mathcal{O}_{\mathrm{cond}}) \cong \lim_{[r]\in\Delta} \prod_{\sigma\in \mathrm{BT}_r} R\Gamma(Z_\sigma, \mathcal{O})$$

We are left with the task of computing the cohomology $R\Gamma(Z_\sigma, \mathcal{O})$ of the nonsingular projective variety $Z_\sigma$. Property (5) above states that $Z_\sigma$ is constructed by a sequence of blow-ups of projective space along nonsingular subvarieties. For any projective space $\mathbf{P}$ over $\kappa$ we have $R\Gamma(\mathbf{P}, \mathcal{O}) \cong \mathbf{F}_p$. We now use the general fact [Hir64, §7, Corollary 2] that if $f\colon X' \to X$ is a blow-up of a nonsingular variety $X$ along a nonsingular subvariety, then $Rf_*\mathcal{O}_{X'} \cong \mathcal{O}_X$ as complexes, so $R\Gamma(Z_\sigma, \mathcal{O}) \cong \mathbf{F}_p$ as well.

Therefore

$$R\Gamma(\mathfrak{H}_s, \mathcal{O}_{\mathrm{cond}}) \cong \lim_{[r]\in\Delta} \prod_{\sigma\in \mathrm{BT}_r} \mathbf{F}_p$$

is quasi-isomorphic to the (condensed!) $\mathbf{F}_p$-cohomology of BT itself. Since BT is exhausted by finite contractible subcomplexes $\mathrm{BT}^i \subset \mathrm{BT}$, its cohomology is $\varinjlim_i \mathbf{F}_p \cong \mathbf{F}_p$. □

**Theorem 6.2.6.** *There is a morphism of algebra objects in* $D(\mathrm{Solid}_{\mathrm{GL}_n(\mathbb{Q}_p)})$*:*

$$\alpha_{\mathcal{H}}\colon \mathbb{Z}_p[\varepsilon] \to R\Gamma(\mathcal{H}_{\mathrm{pro\acute{e}t}}, \hat{\mathcal{O}}^+_{\mathrm{cond}})$$

*such that the* $H^i(\mathrm{cof}(\alpha))$ *are killed by a uniform power of* $p$ *for any* $i$*. Here, the* $\mathrm{GL}_n(\mathbb{Q}_p)$ *acts trivially on* $\mathbb{Z}_p[\varepsilon]$*.*

*Proof.* Combine Corollary 5.6.9 with Theorem 6.2.1. □

**Corollary 6.2.7.** *We have an isomorphism of graded* $\mathbb{Q}_p$*-algebras in* $D(\mathrm{Solid})$*:*

$$H^*([\mathcal{H}^\diamond/\mathrm{GL}_n(\mathbb{Z}_p)]_{\mathrm{pro\acute{e}t}}, \hat{\mathcal{O}}^+) \otimes \mathbb{Q}_p \cong \Lambda_{\mathbb{Q}_p}(y_1, y_3, \ldots, y_{2n-1})[\varepsilon]$$

*Proof.* After taking homotopy fixed points under $\mathrm{GL}_n(\mathbb{Z}_p)$ in Theorem 6.2.6, we find a morphism

$$\alpha_{\mathcal{H}}^{h\,\mathrm{GL}_n(\mathbb{Z}_p)}\colon \mathbb{Z}_p^{h\,\mathrm{GL}_n(\mathbb{Z}_p)}[\varepsilon] \to H^*([\mathcal{H}^\diamond/\mathrm{GL}_n(\mathbb{Z}_p)]_{\mathrm{pro\acute{e}t}}, \hat{\mathcal{O}}^+).$$

By Proposition 3.8.1, the cohomology groups of $\mathrm{cof}(\alpha_{\mathcal{H}}^{h\,\mathrm{GL}_n(\mathbb{Z}_p)}) \cong \mathrm{cof}(\alpha_{\mathcal{H}})^{h\,\mathrm{GL}_n(\mathbb{Z}_p)}$ are all $p$-power torsion. After inverting $p$, we arrive at an isomorphism of $\mathbb{Q}_p$-vector space objects in $D(\mathrm{Solid})$:

$$\mathbb{Q}_p^{h\,\mathrm{GL}_n(\mathbb{Z}_p)}[\varepsilon] \cong H^*([\mathcal{H}^\diamond/\mathrm{GL}_n(\mathbb{Z}_p)]_{\mathrm{pro\acute{e}t}}, \hat{\mathcal{O}}^+) \otimes \mathbb{Q}_p$$

By Proposition 3.8.1, the continuous cohomology of $\mathrm{GL}_n(\mathbb{Z}_p)$ is the stated exterior algebra. □

6.3. **Conclusion of the proof.** We now complete the proofs of Theorem B and of Theorem A.

**Theorem 6.3.1.** *For every $i$, the group $H^i_{\mathrm{cts}}(\mathbb{G}_n, A^c)$ is $p$-power torsion.*

*Proof.* The isomorphism between the two towers, as in Theorem 3.9.1 gives an isomorphism of graded $\mathbb{Q}_p$-algebras:

$$H^*([\mathrm{LT}_K^\diamond/\mathbb{G}_n]_{\text{proét}}, \hat{\mathcal{O}}^+_{\text{cond}}) \otimes \mathbb{Q}_p \cong H^*([\mathcal{H}^\diamond/\,\mathrm{GL}_n(\mathbb{Z}_p)]_{\text{proét}}, \hat{\mathcal{O}}^+_{\text{cond}}) \otimes \mathbb{Q}_p.$$

Combining Corollary 6.1.5 and Corollary 6.2.7 and ignoring the condensed structure, we find an isomorphism of graded $\mathbb{Q}_p$-vector spaces:

$$H^*_{\mathrm{cts}}(G_n, A) \otimes \mathbb{Q}_p[\varepsilon] \cong \Lambda_K(y_1, y_3, \ldots, y_{2n-1})[\varepsilon]$$

Using the decomposition $A = W \oplus A^c$ and the identification of $H^*_{\mathrm{cts}}(\mathbb{G}_n, W)[1/p]$ from Lemma 3.8.5, we find:

$$\Lambda_{\mathbb{Q}_p}(y_1, y_3, \ldots, y_{2n-1})[\varepsilon] \oplus H^*_{\mathrm{cts}}(\mathbb{G}_n, A^c)[1/p] \cong \Lambda_{\mathbb{Q}_p}(x_1, x_3, \ldots, x_{2n-1})[\varepsilon].$$

By comparing dimensions degree-wise on either side, we conclude that $H^*_{\mathrm{cts}}(\mathbb{G}_n, A^c)[1/p] = 0$. □

*Proof of Theorem A.* This is now immediate from Theorem B and Proposition 2.6.3. □

Max Planck Institute for Mathematics, Vivatsgasse 7, 53111 Bonn, Germany
*Email address*: tbarthel@mpim-bonn.mpg.de

The Hebrew University of Jerusalem
*Email address*: tomer.schlank@mail.huji.ac.il

University of Kentucky
*Email address*: nat.j.stapleton@uky.edu

Boston University Department of Mathematics and Statistics, 665 Commonwealth Avenue, Boston, MA, USA
*Email address*: jsweinst@bu.edu